\documentclass[12pt]{amsart}

\usepackage{verbatim}

\newtheorem{theorem}{Theorem}[section]

\newtheorem{corollary}[theorem]{Corollary}

\newtheorem{lemma}[theorem]{Lemma}

\newtheorem{proposition}[theorem]{Proposition}

\theoremstyle{definition}

\newtheorem{definition}[theorem]{Definition}

\newtheorem{remark}[theorem]{Remark}

\newtheorem{example}[theorem]{Example}

\theoremstyle{parrafo}

\newtheorem{parrafo}[theorem]{{\!}}

\numberwithin{equation}{theorem}

\newcommand{\nat}{\mathbb N}

\newcommand{\ent}{\mathbb Z}

\newcommand{\cali}{{\mathcal {I}}}

\newcommand{\calo}{{\mathcal {O}}}

\newcommand{\Sing}{\mbox{Sing\ }}

\title[]{Hypersurface singularities in positive
characteristic.}

\author{Orlando Villamayor U.}

\address{Dpto. Matem\'aticas,  Universidad
Aut\'onoma de Madrid, Canto Blanco 28049 Madrid, Spain.}
\email{villamayor@uam.es}
\thanks{.}


\date{Oct. 2006} \dedicatory{} \commby{}

\begin{document}
\maketitle
\begin{abstract}

The paper is motivated on the open problem of resolution of singularities in positive characteristic. The aim is to present a
form of induction which is different from that used by Hironaka.
In characteristic zero induction is formulated by restriction to smooth hypersurfaces (hypersurfaces of maximal contact). Our alternative approach, introduced here, replaces restrictions to smooth sub-schemes by generic projections on smooth schemes of smaller dimension. 
We introduce a generalization of the discriminant, and our result makes use of the elimination theory. In the case of fields of characteristic zero, elimination gives exactly the same information as the form of induction used by Hironaka.

The properties of this new form of elimination remain weaker in
positive characteristic, than it does in characteristic zero, when
it comes to resolution of singularities. But it opens the way to new invariants for this open problem.
%problems. We also discuss the behavior of invariants, attached to
%a singularity at a point, by using this form of elimination.

\end{abstract}
%%% ----------------------------------------------------------------------

\tableofcontents

%%%%%%%%%%%%%%%%%%%%%%%%%%%%%%%%%%%%%%%%%%%%%%%%%%%%%%%%%%%%%%%%%

%%%%%%%%%%%%%%%%%%%%%%%%%%%%%%%%%%%%%%%%%%%%%%%%%%%%%%%%%%%%%%%%%

%%%%%%%%%%%%%%%%%%%% PART: INTRODUCTION  %%%%%%%%%%%%%%%%%%%%%%%%

%%%%%%%%%%%%%%%%%%%%%%%%%%%%%%%%%%%%%%%%%%%%%%%%%%%%%%%%%%%%%%%%%

%%%%%%%%%%%%%%%%%%%%%%%%%%%%%%%%%%%%%%%%%%%%%%%%%%%%%%%%%%%%%%%%%

\part{Introduction.}

\label{introduction}

%%% ----------------------------------------------------------------------

%%% ----------------------------------------------------------------------
 Let $(S,M)$ denote a local ring, and fix a monic polynomial
$f(Z)$, of degree $b$ in $S[Z]$. This defines a finite ring
extension $S\subset S[Z]/\langle f(Z) \rangle$, and hence a finite
morphism 
\begin{equation}\label{lpeq}
Spec( S[Z]/\langle f(Z)
\rangle)\stackrel{\sigma}{\longrightarrow}Spec(S).
\end{equation}

It seems clear that a lot of information of this finite morphism
is encoded in the coefficients of the monic polynomial $f(Z)$. Higher differential operators
in the variable
$Z$, applied to $f(Z)$, also relate to ramification theoretical methods and to multiplicity theory.

The clarification of this point, initiated in \cite{VV2}, is the
objective of Section 1; where we present the fundamental results
on which our development will rely. The tools developed in this
first section will open the way for explicit computation in the
last section.

Briefly speaking, in the first section we produce a {\em
universal} morphism $\pi$, on which the permutation group acts;
and our invariants arise as an invariant subring by this group
action. The discriminant is an element of the invariant subring.

The Weierstrass Preparation Theorem enables us to reduce the local study of hypersurfaces, embedded in a smooth d-dimensional
scheme, to the study of finite (ramified) covers
defined in terms of one variable $Z$. Namely by $f(Z)\in S[Z]$ as
above, where $S$ is a local regular ring at a point of a smooth
$(d-1)$-dimensional scheme, as we indicate below.

We recall now how higher order differential
operators relate to multiplicity theory, in order to
motivate the development in the coming sections.

Let $V$ be a smooth scheme over a perfect field $k$, and let $J\subset
\calo_V$ be a non-zero sheaf of ideals. For example take $J$ to be
the sheaf of ideals defining a hypersurface. Define a function,
say
$$ord_J :V\to \mathbb{Z},$$ where $ord_J(x)$ denotes the order of
$J_x$ at the local regular ring $\calo_{V,x}$. Let $b$ denote the
biggest value achieved by this function (the biggest order of
$J$). The pair $(J,b)$ is {\em the object of interest} in resolution of singularities and in
Log-principalization of ideals. There is a closed set attached to
this pair in $V$, namely the set of points where $J$ has order
$b$. So if $J$ is a locally principal ideal (defining a
hypersurface), the closed set is the set of points of multiplicity
b at the hypersurface.

For every non-negative integer $s$, the sheaf of $k$-linear
differential operators of order $s$, say $Diff^s_k$, is coherent
and locally free over $V$. There is a natural identification, say
$Diff^0_k=\calo_V$, and for each $s\geq 0$ there is a natural
inclusions $Diff^s_k\subset Diff^{s+1}_k$.

If $U$ is an affine open set in $V$, each $D\in Diff^s_k(U)$ is a
differential operator: $D:\calo_V(U)\to \calo_V(U)$. We define an
extension of the sheaf of ideals $J\subset \calo_V$, say
$Diff^s_k(J)$, so that over the affine open set $U$,
$Diff^s_k(J)(U)$ is the extension of $J(U)$ obtained by adding all
elements $ D(f)$, for all $D\in Diff^s_k(U)$ and $f\in J(U)$.
One can check that $Diff^0(J)=J$, and that $Diff^s(J)\subset Diff^{s+1}(J)$ as sheaves of
ideals in $\calo_V$. Let $V(I)\subset V$ denote the closed set
defined by an ideal $I\subset \calo_V$. So $$V(J)\supset
V(Diff^1(J))\supset \dots \supset V(Diff^{s-1}(J)) \supset
V(Diff^{s}(J))\dots $$

It is simple to check that the order of $J$ at the local
regular ring $\calo_{V,x}$ is $\geq s$ if and only if $x\in
V(Diff^{s-1}(J))$.

The previous observations say that $ord_J :V\to \mathbb{Z}$ is an
upper-semi-continuous function, and that the highest order of $J$
(at points $x\in V$) is $b$ if $V(Diff^{b}(J))=\emptyset$ and
$V(Diff^{b-1}(J))\neq \emptyset$. There is a notion of
transformation of pairs $(J,b)$, defined by monoidal
transformations. Let
$$V  \stackrel{\pi }{\longleftarrow} V_1 \supset H$$
denote the blow up of $V$ at a smooth sub-scheme $Y$, where
$\pi^{-1}(Y)=H $ is the exceptional hypersurface. If $Y\subset
V(Diff^{b-1}(J)))$ we say that $\pi$ is $b$-permissible. In this
case set
$$J\calo_{V_1}=I(H)^b J_1,$$ where $I(H)$ is the sheaf of
functions vanishing along the exceptional hypersurface $H$.

If $J$ has highest order $b$, and if $\pi$ is $b$-permissible, $J_1$ also has at most order $b$ at points
of $V_1$ (i.e., $V(Diff^{b}(J_1))=\emptyset)$. If, in addition,
$J_1$ has no points of order $b$, then we say that $\pi$ defines a
{\em $b$-simplification} of $J$.

If $V(Diff^{b-1}(J_1))\neq \emptyset$, let $V_1\stackrel{\pi_1
}{\longleftarrow} V_2$ denote a monoidal transformation with
center $Y_1\subset V(Diff^{b-1}(J_1))$. We say that $\pi_1$ is
$b$-permissible, and set
$$J_1\calo_{V_2}=I(H_1)^b J_2.$$
So again $J_2$ has at most points of order $b$. If it does, define
a $b$-permissible transformation at some smooth center $Y_2\subset
V(Diff^{b-1}(J_2)))$. 

For $J$ and $b$ as before, we define, by
iteration, a $b$-permissible sequence:
$$V\stackrel{\pi }{\longleftarrow} V_1 \stackrel{\pi_1 }{\longleftarrow}
V_2\stackrel{\pi_2 }{\longleftarrow} \dots  V_{n-1}\stackrel{\pi_n
}{\longleftarrow}  V_{n}, $$ and a factorization
$J_{n-1}\calo_{V_{n}}=I(H_n)^b J_n.$
Let $H_i\subset V_n$ denote the strict transform of the
exceptional hypersurface $H_i \subset V_{i-1}$, $1\leq i \leq n$. Note that:

1) $\{H,H_1,\dots ,H_{n-1}\}$ are components of the exceptional
locus of $V \leftarrow V_n$.

2) The total transform of $J$ relates to $J_n$ by an expression of
the form: $$J\calo_{V_n}=I(H)^{a_0}I(H_1)^{a_1}\cdots
I(H_{n-1})^{a_0}J_{n}.$$

We say that this $b$-permissible sequence defines a {\em
$b$-simplication} of $J\subset \calo_V$ if $\cup H_i$ has normal
crossings, and $V(Diff^{b-1}(J_n))=\emptyset$ (i.e., $J_n$ has
order at most $b-1$ at $V_n$).

When $k$ is a field of {\em characteristic zero}, and $b$ is the
highest order of a sheaf of ideals $J\subset \calo_V$, Hironaka
proves that there is a $b$-simplification. A Log-principalization
of $J$ is achieved when $J_n$ has at most order zero (i.e., when
$J_n=\calo_{V_n}$).

The key point for $b$-simplification, already used in Hironaka's
proof, is a form of induction. In fact, Hironaka proves
$b$-simplification by induction on the dimension of the ambient
space $V$. To simplify matters, assume that $J$ is locally
principal, and let $b$ denote the highest order of $J$ along
points in $V$, which we take now to be a smooth scheme over a field of characteristic
zero. Let
$$ \{ord_J \geq b\}$$ denote the closed set $\{ x \in V /
ord_J(x)\geq b\}$ (or say $=b$).

Fix a closed point $x \in \{ord_J \geq b\}$, and a regular system
of parameters $\{x_1,x_2,\dots ,x_n\}$ for $\calo_{V,x}$. For every
$\alpha =(\alpha_1,\dots, \alpha_n) \in \mathbb{N}^n$, set
$|\alpha|=\alpha_1+\cdots + \alpha_n$, and
\begin{equation*}
\Delta^{\alpha}=\left( \frac{1}{\alpha_1!}\cdots
\frac{1}{\alpha_n!}\right)\frac{\partial^{\alpha_1}}{\partial^{\alpha_1}x_1
}\cdots \frac{\partial^{\alpha_n}}{\partial^{\alpha_n}x_n }.
\end{equation*}
If $J_x$ is locally generated by $f\in \calo_{V,x}$, then $f$ has
order $b$ at $\calo_{V,x}$, and
$$(Diff^{b-1}(J))_x= \langle f, \Delta^{\alpha}(f) /  0\leq |\alpha|<b \rangle.$$
The relevant property will be that the order of $(Diff^{b-1}(J))_x$ at
$\calo_{V,x}$ is one. This holds when $k$ is a field of
characteristic zero.

Recall that, locally at x,  $ V(Diff^{b-1}(<f>))=\{ ord_J \geq b\}$ locally. One way to check that $(Diff^{b-1}(J))_x$ has order one at
$\calo_{V,x}$, is to check this
 at the completion $\hat{\calo}_{V,x}$, say
$R=k'[[x_1,..,x_n]]$. We may choose the system of parameters so
that, for a suitable unit $u$:
$$ uf=f_1=Z^b+a_1Z^{b-1}+\dots +a_b\in S[Z],$$
where $S=k'[[x_1,..,x_{n-1}]]$, and $ Z=x_n$. As $k$ is a field of
characteristic zero, $S[Z]=S[Z_1]$, where $Z_1=Z+\frac{1}{b}a_1$ ,
and we obtain a Tschirnhausen polynomial :
$$f_1=Z_1^b+a_2'Z_1^{b-2}+\dots +a_b' $$  $ (a'_1=0).$
Then:

A)  $Z_1\in Diff^{b-1}(f)$ (in fact
$\frac{\partial^{b-1}f}{\partial^{b-1}Z}\in Diff^{b-1}(f)$). In
particular the ideal $Diff^{b-1}(f)$ has order one at $x$, and the
closed set $\{ ord_J \geq b\}$ is locally included in a smooth
scheme of dimension $n-1$.

B)({\em Elimination of one variable.}) $\{ord \ f\geq b\}(\subset
V(Z_1)) $ can be described as
$$\{ord \ f\geq b\}= \cap_{2 \leq i \leq b}\{ord \ a_i'\geq i\}.$$

C) ({\em Stability: Elimination and monoidal transformations.}) Both A), and the description
in B), are preserved by  $b$-permissible 
transformations.

We will not go into the  details of A),  B) and C). But let us
point out that there is an elimination of one variable in (B). In fact the
closed set $\{ord \ f\geq b\}$ defined in terms of $f$, is also
described as $\cap_{2 \leq i \leq b}\{ord \ a_i'\geq i\},$ where
now the $a'_i$ involve one variable less.

As it was previously indicated, A),B), and C), together, conform the essential
reason and argument in resolution of singularities in
characteristic zero . They rely entirely
on the hypothesis of characteristic zero. For instance A) does not
hold over fields of positive characteristic; so there is no way to
formulate this form of induction over arbitrary fields (see also \cite{HHauser}).
\vskip 1cm

{\bf On this paper.} The objective of these notes is to report on a different approach
to the elimination in B), that can be formulated over arbitrary fields,
which we discuss below (see also \cite{VV3}). Here the restriction to a smooth hypersurface in A) will be replaced by a projection, as in the Weierstrass Preparation Theorem.

In Section 1 we fix a finite morphism $Spec( S[Z]/\langle f(Z)
\rangle)\stackrel{\sigma}{\longrightarrow}Spec(S)$, defined 
by $f(Z)=Z^b+a_1Z^{b-1}+\cdots +a_b\in S[Z]$, and study a closed set, say $F$, in  
$Spec(S)$ over which the finite map is completely ramified, as we indicate below. 
The equations defining this closed set  will be given by polynomial expressions on the coefficients. Among these equations is the discriminant of $f(Z)$.

This will lead us to the study of universal polynomials that can be evaluated on the coefficients of $f(Z)$. These universal polynomials will be obtain via invariant theory, and they are to be thought of as generalized discriminants.
% (as the discriminant  is among these polynomials).

When evaluated on the coefficients of $f(Z)$, these generalized discriminants define $F \subset Spec(S)$. The main result in Section 1 is Theorem 
\ref{ramlocus} which characterizes $F$ in $Spec(S)$ in terms of ramification theory. It is the biggest closed subset so that there is 
a natural bijections $\sigma: \sigma^{-1}(F) \to F$. More precisely, it is the set of points over which the geometric fiber of $\sigma$ is a unique point.

The connexion of this Theorem with multiplicity theory relies on the fact that the $b$-fold points of the hypersurface defined by $f(Z)$ in $Spec( S[Z])$ are included in $\sigma^{-1}(F)$. This property will ultimately ensure a form of compatibility of elimination with monoidal transformations, in the sense of C), addressed in the last Section.

Generalized discriminants are also endowed with a natural weight, a matter to be carefully discussed in Section 1. 
For example, if $D(a_1, a_2,  \dots , a_b)$ denotes the discriminant of $f(Z)$, then it is the evaluation on the $a_i$'s of a universal polynomial, say $D(Y_1, Y_2,  \dots , Y_b)$, which is weighted homogeneous. Equations, as these, endowed with a weight, lead us to the notion of Rees algebras that we discuss below. Appendices 1) and 2) in Section 1 can be avoided in a first look at the paper.

In Section 2 and Section 3 we study the notion of Rees algebras. 
Here we fix $V$, a smooth scheme over a perfect field $k$. A Rees algebras will be of form 
$\mathcal G=\bigoplus_{k\geq 0} I_k W^k$,
$I_0=\calo_V$, and each $I_k$ is an ideal in $\calo_V$. It is a graded subring of $\calo_V[W]$.
 
There are two reasons that justify our attention on these algebras:

i) There is a notion of {\em resolution of Rees algebras}. Moreover, the problem of resolution of singularities reduces to that of resolution of Rees algebras.

ii) Rees algebras can be naturally enriched by the action of higher order differentials.

As for ii) let us indicate simply that these algebras come with a grading. And it is with this grading that one can naturally extend them to new Rees algebras enriched with differential operators. To be precise, each  $\mathcal G$ can be extended to say 
$\mathcal G\subset \mathcal G '$ where now $\mathcal G '$ is what we call a differential Rees algebra, or simply a {\em Diff-algebra}. These are Rees algebras which are, in some sense, compatible with differential operators on $V$. It will be shown that there is a smallest Diff-algebra containing $\mathcal G$, say $\mathcal G\subset Diff(\mathcal G)$. This natural extension has fundamental properties.
It is shown that the construction of a resolution of  $\mathcal G$ is equivalent to a resolution of  $Diff(\mathcal G)$. So, for the purpose 
of resolution of singularities, we may always assume that we start with a differential Rees algebra.

As a Rees algebra over $V$ is a graded subring of $\calo_V[W]$, its integral closure is also a Rees algebra. If two Rees algebras over $V$ have the same integral closure, then it is easy to check that the construction of a resolution of one is equivalent to a resolution of the other. For this reason it is natural to expect that invariants that aim to the definition of constrictive resolutions of Rees algebras should not distinguish one from the other. 
In \cite{VV1} it is proved that if $\mathcal G_1$ and $\mathcal G_2$ are Rees algebras over $V$ with the same integral closure, then 
$Diff(\mathcal G_1)$ and $Diff(\mathcal G_2)$ also have the same integral closure.
Due to the importance of this property, to be used in the further development, 
these results are discussed in Sections 2 and 3.

In Section 4 we develop our form of elimination of one variable,
which is formulated at first in the context of differential Rees algebras. Let
$\mathcal{G}$ be a Diff-algebra over a smooth scheme $V$ of
dimension $d$, and fix $x\in Sing(\mathcal{G})$. We consider here
a suitable smooth morphisms $\pi: V\to V^{(1)}$, defined at an \'etale
neighborhood of $x\in V$,  where $V^{(1)}$ is smooth of dimension
$d-1$. Then a new differential Rees algebra, say
$\mathcal{R}_{\mathcal{G}}$, will be defined over the smooth $d-1$-dimensional
scheme $V^{(1)}$. So $\mathcal{R}_{\mathcal{G}}\subset
\calo_{V^{(1)}}[W]$ is defined in terms of $\pi: V\to V^{(1)}$ and
$\mathcal{G}$. It has the property that $Sing(\mathcal{G})
(\subset V)$ can be identified with
$Sing(\mathcal{\mathcal{R}_{\mathcal{G}}})$ in $V^{(1)}$ via
$\pi$. Here $\mathcal{R}_{\mathcal{G}}$ is called the elimination
algebra.

If we take  $\mathcal{G}$ to be the Diff-algebra spanned by
$\calo_V[fW^b]$, $Sing(\mathcal{G})$ is the set of points of
multiplicity $b$ in the hypersurface defined by $f$, which can be
identified with $Sing(\mathcal{\mathcal{R}_{\mathcal{G}}})$,
defined now in a $(d-1)$-dimensional scheme. This is our approach
to B).

The main properties of elimination algebras are collected in Theorem \ref{4th10}.

Section 5 is devoted to one of the main results in this paper, Theorem \ref{th55}. It is shown there that the main invariant used for resolution of singularities in characteristic zero has a natural extension 
over perfect fields. Yet this invariant is new and has never been treated previously in the study of singularities over fields of positive characteristic.
%
%a we relate elimination with the main invariant used in
%resolution problems ;  these properties will lead us to further
%applications of these results. 
%We also address here the behavior of
%this form of elimination with other known local invariants of
%singularities over arbitrary fields.

The notion of {\em stability of elimination} (see (C)
above) is addressed in Section 6. Results remain stronger over fields of characteristic
zero, where we provide an alternative approach to induction for
desingularization theorems. This form of elimination opens the way
to new questions on resolution problems, also in characteristic
zero; some suggestive properties in this sense are discussed. This
section also includes some explicit calculation of elimination through
examples.

A first look at Definition \ref{defmathcalR} and Theorem
\ref{4th10} is suggested right after reading Section 1, omitting
Appendices 1 and 2 of this first section.

Differential Rees algebras appear in Wlodarczyk's work \cite{WL},
and
 play a central role in Koll\'{a}r's presentation in
(\cite{kollar}), particularly with his notion of {\em tuned
ideals}. Hironaka studies the relation of differential Rees
algebras with integral closure of Rees algebras in
 \cite{Hironaka77},\cite{Hironaka03},\cite{Hironaka05}, in
 connection with the theory of infinitely closed points. These
 connections, and various other aspects, are studied in detail within Kawanoue's recent paper
 \cite{kaw}. See also \cite{VV1}, and \cite{VV3}.

I profited from discussions with Ana Bravo, Vincent Cossart, Marco
Farinati, and Monique Lejeune.

%%%%%%%%%%%%%%%%%%%%%%%%%%%%%%%%%%%%%%%%%%%%%%%%%%%%%%%%%%%%%%%%%%%%%%%%%

%%%%%%%%%%%%%%%%%%%%%%%%%%%%%%%%%%%%%%%%%%%%%%%%%%%%%%%%%%%%%%%%%

%%%%%%%%%%%%%%%%%%%%%%%%%%%%%%%%%%%%%%%%%%%%%%%%%%%%%%%%%%%%%%%%%

\section{ Finite covers: multiplicity  and ramification.}

The aim in this paper is to introduce invariants of singularities over perfect fields. In particular we will refine the notion of multiplicity at points of a hypersurface embedded in a smooth scheme of dimension $d$. Using the Weierstrass Preparation Theorem one can assume that the equation defining the hypersurface is given, locally, by a monic polynomial, in one variable, and coefficients in a smooth scheme of dimension $d-1$.
In this section we study relations of higher order differential
with elimination theory. 

Here we fix a positive integer $b$, and focus on a monic polynomial $f(Z)\in S[Z]$, where  $S$ is a $k$-algebra. The objective is to study the ramification of the finite 
map in (\ref{lpeq}).
%The main example is the universal polynomial of degree $b$ of 
Take $F(Z)=(Z-Y_1)\cdot
(Z-Y_2)\cdots (Z-Y_b)$, in the polynomial ring in $b+1$ variables
$k[Y_1,\dots ,Y_b, Z]$. Consider the $k$-algebra
$S=k[Y_1,\dots ,Y_b]^{\mathbb{S}_b}$, ring of invariants where
$\mathbb{S}_b$ is the group of permutation of $b$ elements. Then
$F(Z)\in S[Z]$ is monic in this $k$-algebra, and it is the {\em
universal} monic polynomial of degree $b$ within the class of $k$-algebras.

The theory of elimination (of the variable $Z$) leads us to the study
of the action of $\mathbb{S}_b$ on a subring of $k[Y_1,\dots
,Y_b]$; namely on the $k$-algebra generated by all differences
$Y_j-Y_k$. This provides a new ring of invariants which is the
object of interest of Section 1, and a basic tool to be considered
throughout the paper. Note that the discriminant of $F(Z)$ is an
element of this invariant ring. 

The section begins with the
relation of this invariant ring with higher order differentials
applied to $F(Z)$ (Remark \ref{rkjun06}); and the main result in this section is Theorem \ref{ramlocus}. 

Appendices 1) and 2) can be avoided in a first look at the paper.

%
%Other invariant rings will arise when we consider more then one
%polynomial. For example the resultant of two polynomials $F(Z)$
%and $G(Z)$ is a natural concept in elimination. Generalization to
%more then one polynomial, which begins in \ref{rr1}, will give
%rise to different rings of invariants, and inclusions among them.
%Of special interest, due to applications to resolution problems,
%will be the study of {\em elimination algebras} and of finite
%extensions among them (Def \ref{1def29} and \ref{uni660}).

%

\begin{parrafo} \label{parcero}{\em Let $(S,M)$ denote a local ring, and fix a monic polynomial
$f(Z)$, of degree $b$ in $S[Z]$. This defines a finite ring
extension $S\subset S[Z]/\langle f(Z) \rangle$, and hence a finite
morphism $$Spec( S[Z]/\langle f(Z)
\rangle)\stackrel{\pi}{\longrightarrow}Spec(S)$$

This finite morphism is said to be purely ramified at 
$P\in Spec(S)$, when the geometric fiber at $P$ has a unique
point. Equivalently, set $\overline{k(P)}$ an algebraically closed field
extension of the residue field $k(P)$, then the morphism is purely
ramified at $P$ if and only if the class of $f(Z)$ in
$\overline{k(P)}[Z]$ has a unique root.

We begin by describing the set of prime ideals in $Spec(S)$ for
which the finite extension is purely ramified. Our arguments, for
this and further properties of this finite morphism, focus on two
observations. Since the finite morphism is determined by the monic
polynomial $f(Z)\in S[Z]$, many properties of the morphism should
be encoded in the coefficients. On the other hand, changes of
variables of the form $Z_1=Z-s$, $s\in S$ do not affect the finite
extension $S\subset S[Z]/\langle f(Z) \rangle$.}
\end{parrafo}

\begin{definition} Consider the homomorphism of $S$-algebras, say
$Tay:S[Z]\to S[Z,T]$, defined by setting $Tay(Z)=Z+T$. For each
$F(Z)\in S[Z]$, $$Tay(F(Z))=\sum_{k\geq 0} g_k(Z)T^k.$$

This defines, for each index $r \geq 0 $,  operators $\Delta^r :
S[Z]\to S[Z]$, by setting $\Delta^r(F(Z))=g_r(Z)$, or say
$$Tay(F(Z))=\sum \Delta^r(F(Z))T^r.$$ The $\Delta^r$ are
differential operators of order $r$. These are $S$-linear
operators; $\Delta^n(Z^{n})=1$ and $\Delta^r(Z^{n})=0$ for $r>n$.

Note that the morphism $Tay$, and the operators $\Delta^r $, are
compatible with change of the base ring $S$. The following result is included simply to illustrate their link with ramification.
% and also with changes
%of variable, in $S[Z]$, of the form $Z_1=Z-s$, $s\in S$.
\end{definition}
\begin{lemma}\label{lemma1} Fix $f(Z)\in K[Z]$, monic of degree $b$, where $K$ is
an algebraically closed field. The following are equivalent:

1) $f(Z)=(Z-\alpha)^b$ for some $\alpha \in K$.

2) The class of $\Delta^k (f(Z))$ in $K[Z]/\langle f(Z) \rangle$
is nilpotent for all integer $0\leq k<b$.

\end{lemma}

\proof That 1) implies 2) is clear.

Let $f(Z)=(Z-\alpha_1)^{\beta_1}\cdot (Z-\alpha_2)^{\beta_2}\cdots
(Z-\alpha_r)^{\beta_r}$ be the expression of the monic polynomial
in terms of its $r$ different roots; so $\sum \beta_i=b$. We prove
that 2) implies 1) by showing that the class of the ideal $\langle
\Delta^k(f(Z)), 0\leq k\leq b-1 \rangle \subset K[Z]$ is nilpotent
in $K[Z]/\langle f(Z) \rangle$ only when $r=1$. Assume that $r\geq
2$, and set $G(Z)=(Z-\alpha_2)^{\beta_2}\cdots
(Z-\alpha_r)^{\beta_r}$, so $f(Z)=(Z-\alpha_1)^{\beta_1}\cdot
G(Z)$, and $\beta_1<b$.
Now $$Tay(f(Z))=Tay((Z-\alpha_1)^{\beta_1})\cdot Tay(G(Z))\in
K[Z,T],$$ and $Tay((Z-\alpha_1)^{\beta_1})$ is a monic polynomial
of degree $\beta_1$ in $T$.

On the other hand, at $(K[Z]/\langle Z-\alpha_1 \rangle )[T]$
($=K[T]$), $Tay((Z-\alpha_1)^{\beta_1})=  T^{\beta_1}$  and
$$Tay(G(Z))= G(\alpha_1)+ \mbox{terms of degree } \geq 1 \mbox{ in
} T,$$ where $G(\alpha_1)$ (the class of $G(Z)$ in $K[Z]/\langle
Z-\alpha_1 \rangle $) is non-zero. This shows that
$\Delta^{\beta_1}(f(Z))\notin \langle Z-\alpha_1 \rangle $, and
hence the class of $\Delta^{\beta_1}(F(Z))$ in the ring
$K[Z]/\langle f(Z) \rangle$ is not nilpotent.
\endproof
\begin{parrafo}\label{casas}{\rm The Lemma relates with a stronger
question raised by E. Casas-Alvero: if either $f(Z)=(Z-\alpha)^b$
for some $\alpha \in K$, or $f(Z)$ and $\Delta^{\alpha}(f(Z))$ are
coprime in $K[Z]$ for some $ 1\leq \alpha \leq b-1$.}
\end{parrafo}
\begin{parrafo}{\bf On the general strategy.}\label{ontges} {\em
Fix a ring $k$ and consider $S$ in the class of
$k$-algebras. We recall that there is a universal monic polynomial of degree
$b$ within this class. In our further applications $k$ will be a field, but to some extend the form of elimination we discuss here works over arbitrary rings, in fact over the integers.

 Let $k[Y_1,\dots ,Y_b]$ be a polynomial ring over $k$.
We will denote by $R_b$ the ring of symmetric polynomials in $b$
variables with coefficients in $k$. Here $s_{b,1},\ldots,s_{b,b}$
will denote the symmetric polynomials, where each $s_{b,i}$ is
homogeneous of degree $i$ in $k[Y_1,\dots ,Y_b]$. Therefore
$R_b=k[s_{b,1},\ldots,s_{b,b}]$ is a graded subring in
$k[Y_1,\dots ,Y_b]$. Set $$F_b(Z)=(Z-Y_1)\cdot (Z-Y_2)\cdots
(Z-Y_b)\in k[Y_1,\dots ,Y_b][Z],$$  the generic polynomial of
degree $b$. Recall that
$$F_b(Z)=Z^b-s_{b,1}Z^{b-1}+\dots +(-1)^bs_{b,b}\in R_b[Z],$$ and
note that every monic polynomial of degree $b$, say $f(Z)=Z^b+a_1Z^{b-1}+\dots +a_b\in S[Z],$ over a $k$-algebra $S$,
is obtained from $F_b(Z)\in R_b[Z]$ by a unique morphism $R_b\to S$.

Let ${\mathbb{S}_b}$ denote the symmetric group acting on
$k[Y_1,\dots ,Y_b]$ in the usual manner, so that
\begin{equation}\label{lrbeq}
R_b=k[s_{b,1},\ldots,s_{b,b}]=k[Y_1,\dots ,Y_b]^{{\mathbb{S}_b}}.
\end{equation}

We shall show, via Galois correspondence, that there are natural
$R_b$-isomorphisms, or say identifications: $$R_b[Z]/\langle F_b(Z) \rangle = R_b[Y_j],$$
for every index j, where $R_b[Y_j]$ is a subring of $k[Y_1,\dots
,Y_b]$ (\ref{rk3}, 3)).

Many properties of a finite extension $S \subset
S[Z]/\langle f(Z) \rangle $, defined by a monic polynomial
$f(Z)\in S[Z]$, can be expressed in terms of the
coefficients, and are independent of changes of the form
$Z_1=Z-s$, $s\in S$. Such is the case of the ramification locus of
the induced finite morphism. In fact, the ramification is described by the discriminant, which
is invariant by those changes of $Z$. 

Changes of variables in $S[Z]$ are more general, namely of the form $Z_1=uZ-s$, where $s\in S$ and $u$ is a unit of $S$. We indicate below why it suffices, for our purpose, to restrict to changes of the form $Z_1=Z-s$, at the universal level, in order to obtain, on each specific $S[Z]$, invariants for arbitrary changes $Z_1=uZ-s$.

Note that the action of ${\mathbb{S}_b}$ on $k[Y_1,\dots ,Y_b]$
restricts to an action on the subring $$k[Y_2-Y_1,\dots
,Y_b-Y_1]=k[Z_{i,j}]$$ where $Z_{i,j}=Y_i-Y_j$, $1\leq i,j\leq b$.
In what follows we denote by
\begin{equation}\label{eqerrebarra}
\overline{R}_b=k[Y_2-Y_1,\dots ,Y_b-Y_1]^{{\mathbb{S}_b}},
\end{equation}
the subring of invariants by this action. So $$\overline{R}_b=R_b
\cap k[Y_2-Y_1,\dots ,Y_b-Y_1](\subset k[Y_1,\dots ,Y_b]).$$

Elements of $\overline{R}_b$, are elements in $R_b$, that provide,
for every monic polynomial $f(Z)\in S[Z]$ of degree $b$, equations
on the coefficients which are independent of changes of the form
$Z_1=Z-s$, $s\in S$. Our aim is to study generators of
$\overline{R}_b$, and also its weighted structure as subring of
the graded ring $k[Y_1,Y_2, \dots ,Y_b]$. It is also a graded
subring of $R_b$.  Here $R_b$ is mapped to $S$, so that $F_b(Z)$ defines
$f(Z)=Z^b+a_1Z^{b-1}+\dots +a_b\in S[Z]$ as
above. In particular, an homogenous element $H$ of degree, say
$r$, in $\overline{R}_b$, maps to an element, say $h$, in $S$, which is a polynomial expression on the coefficients $a_i$ of $f(Z)$.
However, in our development we will want to recall the degree of
$H$. To this end, in \ref{lavw}, we will add a dummy variable $W$
, and we assign to $H$ the element $h\cdot W^r$ in the ring
$S[W]$. The algebra $\overline{R}_b$ is generated by finitely many homogeneous elements in $R_b$, say
$\overline{R}_b=k[H_1,\dots ,H_s]$, where each $H_i$ is
homogeneous of degree, say $r_i$. We will consider, by a change of the base
ring $\overline{R}_b\to S$, the $S$ subalgebra of $S[W]$ generated by $\{ h_1\cdot
W^{r_1},\dots ,h_s\cdot W^{r_s}\}$, namely 
$$S[h_1\cdot W^{r_1},\dots , h_s\cdot W^{r_s}].$$
This is a Rees algebra, a direct sum of ideals in $S$, say 
$\bigoplus_{r\geq 0} I_r W^r$. Now for each index $r$, the ideal 
$I_r$ is generated by polynomials on the coefficients of  $f(Z)$.
They are clearly invariants by a change of variable $Z_1=Z-s$, $s\in S$, as each $h_i$ is invariant by such change. But the point is that each 
ideal $I_r$ will be independent of any change $Z_1=uZ-s$ in $S[Z]$.
This last property will rely on the fact that the previous (universal) $H_i$ are weighted homogeneous on the symmetric functions $s_{b,i}$.

Therefore the grading of these invariant rings play a central role in our further discussion. $\overline{R}_b$ is just one
example of a graded ring that will arise. A first look into 
\ref{parraweightord} and the formulation of Theorem \ref{ramlocus}  is suggested to understand the connection
of these invariant rings with the notion of ramification.

}
\end{parrafo}

\begin{remark}\label{rk3}

1) Notice that $$\mbox{ dim}\left(k[Y_1-Y_2, Y_1-Y_3,\dots,
Y_1-Y_b]\right)=\mbox{dim }\left((k[Y_1-Y_2, Y_1-Y_3,\dots,
Y_1-Y_b])^{\mathbb{S}_b}\right)=b-1.$$

2) Identify $\mathbb{S}_{b-1}$ with the subgroup of permutation in
$\mathbb{S}_b$ fixing $Y_1$, so $$(k[Y_1-Y_2, Y_1-Y_3,\dots,
Y_1-Y_b])^{\mathbb{S}_b}=\left((k[Y_1-Y_2, Y_1-Y_3,\dots,
Y_1-Y_b])^{ \mathbb{S}_{b-1}}\right)^{ \mathbb{S}_b}.$$

3) $(k[Y_1-Y_2, Y_1-Y_3,\dots, Y_1-Y_b])^{ \mathbb{S}_{b-1}}\subset
k[Y_1,\dots ,Y_b]^{ \mathbb{S}_{b-1}}=k[s_{b,1},\ldots,s_{b,b}][Y_1].$

We check that $k[Y_1,\dots ,Y_b]^{
\mathbb{S}_{b-1}}=k[s_{b,1},\ldots,s_{b,b}][Y_1]$. The inclusion
$$k[Y_1,\dots ,Y_b]^{ \mathbb{S}_{b-1}}\supset
k[s_{b,1},\ldots,s_{b,b}][Y_1]$$ is clear. Let
$t_{b-1,1},\ldots,t_{b-1,b-1}$ denote the symmetric polynomials in
$b-1$ variables $Y_2,\dots,Y_b$, and note that $k[Y_1,\dots
,Y_b]^{ \mathbb{S}_{b-1}}=k[Y_1,t_{b-1,1},\ldots,t_{b-1,b-1}]$.
For the other inclusion check that

$s_{b,1}=Y_1+t_{b-1,1}$; $s_{b,2}=Y_1t_{b-1,1}+t_{b-1,2}$; \dots
$s_{b,b-1}=Y_1t_{b-1,b-2}+t_{b-1,b-1}$; and
$s_{b,b}=Y_1t_{b-1,b-1}$.

\bigskip

4) Notice that $$(k[Y_1-Y_2, Y_1-Y_3,\dots, Y_1-Y_b])^{ \mathbb{S}_{b-1}}= k[
t_{b-1,1}(Y_1-Y_j),\ldots,t_{b-1,b-1}(Y_1-Y_j)],$$ where, as
before, $t_{b-1,1},\ldots,t_{b-1,b-1}$ denote the symmetric
polynomials in $b-1$ variables, evaluated here at the elements
$\{Y_1-Y_2, Y_1-Y_3,\dots, Y_1-Y_b\}$.
\end{remark}

\begin{parrafo} \label{par04} {\em Define the morphism $Tay:R_b[Z]=k[s_{b,1},\ldots,s_{b,b}][Z] \to
R_b[Z,T]$, and $\Delta^{(\alpha)}:R_b[Z] \to R_b[Z]$, as usual, by
setting $Tay(F_b(Z))=F_b(Z+T)=\sum
\Delta^{(\alpha)}(F_b(Z))T^{\alpha}.$

In what follows recall the natural identification
$$k[s_{b,1},\ldots,s_{b,b}][Y_1]=k[s_{b,1},\ldots,s_{b,b}][Z]/\langle
F_b(Z) \rangle ,$$ and that the $\Delta^{(\alpha)}$ operators are,
in a natural sense, compatible with change of base rings $R_b \to
S$ within the class of algebras over $k$.}
\end{parrafo}
\begin{remark}\label{rkjun06}

Since $F_b(Z)=(Z-Y_1)\cdot (Z-Y_2)\cdots (Z-Y_b)\in k[Y_1,\dots
,Y_b][Z],$ $$F_b(T+Z)=(T+(Z-Y_1))\cdot (T+(Z-Y_2))\cdots
(T+(Z-Y_b))\in k[Y_1,\dots ,Y_b][Z,T].$$

Let $F^{(\alpha)}(Z)$ denote the element
$\Delta^{(\alpha)}(F_b(Z))$. The coefficients of this polynomial
in the variable $T$, are the symmetric polynomials evaluated on
the elements $Z-Y_j$, $1\leq j\leq b$. So
\begin{equation}\label{eqjuni06}
F_b^{(e)}(Z)=(-1)^{b-e}s_{b,b-e}(Z-Y_1, Z-Y_2,\dots , Z-Y_b).
\end{equation}
Here we extend the action of  $ \mathbb{S}_b$, acting on the
variables $Y_j$, setting $\sigma(Z)=Z$. Note that
$$k[Z-Y_1,\dots , Z-Y_b]^{ \mathbb{S}_b}=k[\{F_b^{(e)}(Z), e=0,1,\dots, b-1\}],$$
and that each $F_b^{(e)}(Z)$ is homogeneous. Observe that
$(Z-Y_1)-(Z-Y_j)=Y_j-Y_1$, so $\overline{R}_b$
 is also a subring of
$k[\{F_b^{(e)}(Z), e=0,1,\dots, b-1\}]$ (\ref{eqerrebarra}).
\end{remark}

\begin{lemma}\label{lemma2} Let the setting be as in \ref{rk3}, 4). Then $$k[
t_{b-1,1}(Y_1-Y_j),\ldots,t_{b-1,b-1}(Y_1-Y_j)]=k[F_b^{(1)}(Y_1),\dots
, F_b^{(b-1)}(Y_1)](\subset k[s_{b,1},\ldots,s_{b,b}][Y_1]).$$ In
fact: $F_b^{(e)}(Y_1)=(-1)^{b-e}t_{b-1,b-e}(Y_1-Y_j), \hskip 0.2cm
1\leq e \leq b-1.$
\end{lemma}
\proof This is a consequence of \ref{eqjuni06}. Note also that the
equality
$$s_{b,b-e}(Y_1-Y_1, Y_1-Y_2,\dots
, Y_1-Y_b)=t_{b-1,b-e}(Y_1-Y_2,\dots , Y_1-Y_b)$$ follows from the
definition of symmetric polynomials.
\endproof
\begin{corollary}\label{cor1}

1) $$(k[Y_1-Y_2, Y_1-Y_3,\dots, Y_1-Y_b])^{\mathbb{S}_{b-1}}=
k[F_b^{(1)}(Y_1),\dots , F_b^{(b-1)}(Y_1)](\subset
k[s_{b,1},\ldots,s_{b,b}][Y_1]).$$

2) $\overline{R}_b=k[F_b^{(1)}(Y_1),\dots , F_b^{(b-1)}(Y_1)]\cap
k[s_{b,1},\ldots,s_{b,b}],$
as subrings of
$k[s_{b,1},\ldots,s_{b,b}][Y_1]$(\ref{eqerrebarra}).

\end{corollary}
\begin{parrafo}\label{rk07} {\em  Let $k[Y_2-Y_1,\dots , Y_b-Y_1]$ be graded as subring of $k[Y_1,\dots, Y_b]$ with the usual grade.
Since the action of $ \mathbb{S}_b$ preserves degrees, it follows
that $\overline{R}_b$ and $k[F_b^{(1)}(Y_1),\dots ,
F_b^{(b-1)}(Y_1)]$ are graded subrings of $k[Y_2-Y_1,\dots ,
Y_b-Y_1]$. So we may assume that $\overline{R}_b=k[G_{b,1},\dots
,G_{b,r_b}]$, where the generators $G_{b,i}$ are homogeneous
polynomials. Furthermore,
$$\overline{R}_b=k[F_b^{(1)}(Y_1),\dots , F_b^{(b-1)}(Y_1)]\cap
k[s_{b,1},\ldots,s_{b,b}]= k[F_b^{(1)}(Y_1),\dots ,
F_b^{(b-1)}(Y_1)]^{ \mathbb{S}_b}.$$ So $\overline{R}_b\subset
k[F_b^{(1)}(Y_1),\dots , F_b^{(b-1)}(Y_1)]$ is a finite extension
of graded rings. Therefore $$\langle G_{b,1},\dots ,G_{b,r_b}
\rangle  \subset \langle F_b^{(1)}(Y_1),\dots ,
F_b^{(b-1)}(Y_1)\rangle $$ are homogeneous  ideals in
$k[F_b^{(1)}(Y_1),\dots , F_b^{(b-1)}(Y_1)]$, and
$$\sqrt{\langle G_{b,1},\dots ,G_{b,r_b} \rangle }=\langle
F_b^{(1)}(Y_1),\dots , F_b^{(b-1)}(Y_1) \rangle $$ in the ring
$k[F_b^{(1)}(Y_1),\dots , F_b^{(b-1)}(Y_1)]$ (see \ref{grocal}).
This, together with the inclusion $$k[F_b^{(1)}(Y_1),\dots ,
F_b^{(b-1)}(Y_1)]\subset k[s_{b,1},\ldots,s_{b,b}][Y_1]$$ in
(\ref{cor1}), 1), show that

$$\sqrt{\langle G_{b,1},\dots
,G_{b,r_b} \rangle}=\sqrt{ \langle F_b^{(1)}(Y_1),\dots ,
F_b^{(b-1)}(Y_1)\rangle}$$ as ideals in
$k[s_{b,1},\ldots,s_{b,b}][Y_1]$ (see \ref{rkperma}). }
\end{parrafo}

\begin{remark}\label{rkperma}
Fix a ring $A$ and two ideals $I_1$, $I_2$. If $I_1 \subset I_2$,
then for every ring homomorphism $A \to B$, $I_1B \subset I_2B$. A
similar property holds if $\sqrt{I_1}=\sqrt{I_2}$ in $A$. We set
$R_b=k[s_{b,1},\ldots,s_{b,b}]=k[Y_1,\dots ,Y_b]^{ \mathbb{S}_b},$
and the universal polynomial of degree $b$:
$$F_b(Z)=Z^b-s_{b,1}Z^{b-1}+\dots +(-1)^bs_{b,b}\in R_b[Z].$$

Any monic polynomial of degree $b$, say $f(Z)=Z^b-a_1Z^{b-1}+\dots
+(-1)^ba_b\in R[Z],$
 over a $k$ algebra $R$,
is obtained from $F_b(Z)\in R_b[Z]$, by the $k$-algebra
homomorphism $R_b \to R$ defined by mapping $s_{b,i}$ to $a_i$. In
particular, there is a natural homomorphism $$R_b[Z]/\langle
F_b(Z) \rangle  \to R[Z]/\langle f(Z)\rangle $$ defined by the base
change $R_b \to R$ , and the result in \ref{rk07} ensures that $$\sqrt{\langle
G_{b,1}(a_1,\dots, a_b),\dots ,G_{b,r_b}(a_1,\dots, a_b)\rangle
}=\sqrt{\langle \overline{f^{(1)}(Z)},\dots ,
\overline{f^{(b-1)}(Z) \rangle }}$$ as ideals in $R[Z]/\langle
f(Z)\rangle$.

\end{remark}
\begin{remark}\label{grocal} Let $F_1\subset F_2$ be an inclusion
of finitely generated $\mathbb{N}$-graded algebras over a field
$k$. Let $max_1$ and $max_2$ denote the irrelevant maximal ideals
of $F_1$ and $F_2$ respectively. Note that if $F_1\subset F_2$ is
a finite extension, then $\sqrt{max_1F_2}=max_2$.
\end{remark}
\begin{remark}\label{gradlocal}
Fix a finitely $\mathbb{N}$-graded algebra, say $k[H_1, \dots ,
H_r]$, where each $H_i$ is homogeneous of degree $n_i$, and a
local regular $k$-algebra $(S,M)$. We say that an homomorphism of
$k$-algebras $ \phi:k[H_1, \dots , H_r] \to S$ {\em preserves
degrees}, if, for every homogeneous element $H\in k[H_1, \dots ,
H_r]$ of degree $n$ : $$\nu_S(\phi(H))\geq n,$$ where $\nu_S$
denotes the order at the local regular ring.

This property holds if and only if $\nu_S(\phi(H_i))\geq n_i$,
$1\leq i \leq r$. This follows from the fact that every
homogeneous element $H$, of degree $n$, can be expressed as
$H=G(H_1,\dots ,H_r)$, where $G(Z_1,\dots ,Z_r)\in k[Z_1,\dots
,Z_r]$ is weighted homogeneous of degree $n$, where each $Z_i$ is
considered with degree $n_i$.

\end{remark}

\begin{parrafo}\label{parraweightord} {\em Set $\overline{R}_b(=k[Y_1-Y_2, Y_1-Y_3,\dots,
Y_1-Y_b])^{ \mathbb{S}_{b}})=k[G_{b,1},\dots ,G_{b,r_b}]$; and
assume that each generator $G_{b,i}$ is homogeneous of degree
$n_i$, as polynomial in $k[Y_1, \dots , Y_b]$. So an homogeneous
element $H\in \overline{R}_b$, of degree $m$ in $k[Y_1, \dots , Y_b]$, is also a weighted
homogeneous polynomial of degree $m$ in the $G_{b,i}$'s, provided
each $G_{b,i}$ is given weight $n_i$.

The same argument applies for the inclusion $\overline{R}_b\subset
k[s_1,\dots, s_b],$ so $H$ is also a weighted homogeneous
polynomial of degree $m$ in the $s_i$'s, provided each $s_{i}$ is
given weight $i$.

Let $(S,M)$ denote, as above, a local regular ring; and recall
that $\phi:\overline{R}_b \to S $ preserves degrees if and only if
$\nu_S(\phi(G_{b,i}))\geq n_i,$ for $1\leq i\leq r_b$ (
\ref{gradlocal} ). Consider a regular $k$-algebra $R$. A monic
polynomial $f(Z)=Z^b+a_1Z^{b-1}+\ldots+a_{b-1}Z+a_b\in R[Z]$
defines a hypersurface in the regular scheme $Spec(R[Z])$, and a
finite morphism
$$Spec(R[Z]/\langle f(Z) \rangle )\stackrel{\pi}{\rightarrow}Spec(R).$$
Let $Q$ be a point of this hypersurface. Set $P=\pi(Q)$, and
$S=R_P$ (local regular ring).

We claim that if the hypersurface has multiplicity $b$ at $Q$,
then
 $\nu_S(G_{b,i}(a_1,a_2,\dots, a_b))\geq n_i,$ for $1\leq
i\leq r_b$. In order to prove this claim we use Zariski's
multiplicity formula (see \cite{ZS}, Corollary 1, page 299). It
follows that if $Q$ is a $b$-fold point of this hypersurface,
there is a suitable change of coordinate $Z_1=Z- s$, $s\in S$, so
that: $f(Z)=Z_1^b+c_1Z_1^{b-1}+\ldots+c_{b-1}Z_1+c_b\in S[Z]\hskip
0.1cm \mbox{, and } \nu_S(c_i)\geq i.$ In fact, the multiplicity
formula, applied to the finite extension $S\subset B=S[Z]/\langle
f(Z)\rangle$, ensures that $B$ is local (i.e., $B=B_Q$) and the
residue fields of the local rings $B$ and $S$ are the same. In
particular, the class of $Z$ in the residue field of $B$ is the
class of some element $s\in S$.

 Recall that $k[s_{b,1},\ldots,s_{b,b}]$ is a graded
subalgebra of $k[Y_1,\dots ,Y_b]$, and each $s_{b,i}$ is
homogeneous of degree $i$. The morphism
$k[s_{b,1},\ldots,s_{b,b}]\to S$, defined by mapping
$(-1)^is_{b,i}$ to $c_i$, maps $G_{b,i}(s_{b,1},\ldots,s_{b,b}) $
to $G_{b,i}(c_1,\ldots,c_b)$. As this morphism
$k[s_{b,1},\ldots,s_{b,b}]\to S$ preserves degrees, it follows
that $\nu_S(G_{b,i}(c_1,\ldots,c_b)) \geq n_i.$

On the other hand
$G_{b,i}(c_1,\ldots,c_b)=G_{b,i}(a_1,\ldots,a_b)$ since these
functions are invariant by these changes of the coordinate $Z$,
hence $\nu_S(G_{b,i}(a_1,\ldots,a_b)) \geq n_i.$ (\ref{casas})}
\end{parrafo}
\begin{theorem}
\label{ramlocus} Let $R$ be a $k$ algebra,
$f(Z)=Z^b+a_1Z^{b-1}+\ldots+a_{b-1}Z+a_b\in R[Z]$, and set
$Spec(R[Z]/\langle f(Z) \rangle
)\stackrel{\pi}{\rightarrow}Spec(R)$. Then:

i)$V(\langle G_{b,1}(a_1,\ldots,a_b),\ldots,
G_{b,r_b}(a_1,\ldots,a_b) \rangle) $ is the set of points in $
\mbox{Spec}(R)$ where the finite morphism is purely ramified
(\ref{parcero}).

ii) If $R$ is regular, and $Q\in V(\langle f(Z)\rangle)$ is a
point of multiplicity $b$ of this hypersurface in $Spec(R[Z])$,
then
$$\nu_S(G_{b,i}(a_1,\ldots,a_b)) \geq n_i$$ for $1\leq i\leq r_b$,
where $S=R_P$, $P=\pi(Q)$.

\end{theorem}
\proof

i) Note that $ \mbox{Spec}\left(S[Z]/\langle f\rangle\right) \to
\mbox{Spec}(S)$ arises from $ \mbox{Spec}\left(R_b[Z]/\langle
F_b(Z)\rangle\right) \to \mbox{Spec}(R_b)$ by the change of base
rings  $\phi: R_b\to S$, where $\phi$ is a $k$ algebra morphism, and
$\phi((-1)^is_{b,i})=a_i$. So, as it was indicated in \ref{rkperma}, it
suffices to prove the claim for the universal case.

In \ref{rk07} we show that $\langle G_{b,1},\dots ,G_{b,r_b}
\rangle$ and $\langle F_b^{(1)}(Y_1),\dots , F_b^{(b-1)}(Y_1)
\rangle$ have the same radical ideal,
 as ideals in
$R_b[Z]/\langle F_b(Z)\rangle$ (
$=k[s_{b,1},\ldots,s_{b,b}][Y_1]$). Fix a prime $P\subset
k[s_{b,1},\ldots,s_{b,b}]$ and set $\{Q_1,\dots,Q_s\}$ the primes
in $R_b[Z]/\langle F_b(Z)\rangle$, over $P$. Let $K$ be an
algebraic closure of the residue field of $(R_b)_P$, and argue as
in Lemma \ref{lemma1}. If $P$ contains the ideal $\langle
G_{b,1},\dots ,G_{b,r_b} \rangle$ in $R_b$, then any $Q_i$
contains $ \langle F_b^{(1)}(Y_1),\dots , F_b^{(b-1)}(Y_1)
\rangle$
 as ideals in
$R_b[Z]/\langle F_b(Z)\rangle$, so Lemma \ref{lemma1} asserts that
the fiber over $P$ is purely ramified. Conversely, if $P$ does not
contain the ideal $\langle G_{b,1},\dots ,G_{b,r_b} \rangle$ in
$R_b$, then $ \langle F_b^{(1)}(Y_1),\dots , F_b^{(b-1)}(Y_1)
\rangle$ (in $R_b[Z]/\langle F_b(Z)\rangle$) is not contained in
any $Q_i$. In this case Lemma \ref{lemma1} asserts that the
morphism is not purely ramified at $P$.

ii) Proved in \ref{parraweightord}.
\endproof
\centerline{\bf Appendix 1: On normality and graded structure of
$\overline{R_b}$.}

\medskip

 In
the coming sections we will study invariants that arise from the
graded structure of the rings $\overline{R_b}=k[Y_2-Y_1,\dots
,Y_b-Y_1]^{ \mathbb{S}_b} $ and $k[Z-Y_1,\dots ,Z-Y_b]^{
\mathbb{S}_b} $.

We now introduce a graded subring of $\overline{R_b}$, whose
integral closure is $\overline{R_b}$. This subring will be particularly
useful for explicit computation of invariants (see examples of the
last Section 6). This appendix can be omitted in a first look at the
paper.

\begin{remark}\label{rk09}
We have studied graded subrings of $k[Y_1,\dots ,Y_b,Z]$, which we
consider with the usual grade. In particular an element of a
graded subring is homogeneous if and only if it is homogeneous in
$k[Y_1,\dots ,Y_b,Z]$. A subring of this ring is graded when it is generated by
homogeneous elements. Since
$\overline{R_b}$ is the subring of $ \mathbb{S}_b$ invariants in
$k[Y_2-Y_1,\dots ,Y_b-Y_1]$,
$$\overline{R_b}=k[Y_2-Y_1,\dots ,Y_b-Y_1]\cap R_b,$$ as subrings
of $k[Y_1,\dots ,Y_b]$ (see \ref{lrbeq}). In particular $\overline{R_b}$ is an
intersection of normal rings, so it is normal.
\end{remark}

\begin{parrafo}\label{par010}{\em
 The ring $k[s_{b,1},\ldots,s_{b,b}][Y_1]$ is a free module of
rank $b$ over $k[s_{b,1},\ldots,s_{b,b}]$. Multiplication by an
element $\Theta \in k[s_{b,1},\ldots,s_{b,b}][Y_1]$ defines an
endomorphism, say $\phi_{\Theta}$, with characteristic polynomial,
say $$\psi_{\Theta}(V)=V^b+h_1V^{b-1}+\dots +h_b\in
k[s_{b,1},\ldots,s_{b,b}][V].$$}
\end{parrafo}
\begin{lemma}\label{lemma011}
We claim that, if $\Theta= a_0+a_1Y_1+\dots a_{b-1}Y_1^{b-1}$,
$$\psi_{\Theta}(V)=\prod_{1\leq j \leq b} (V-(a_0+a_1Y_j+\dots
a_{b-1}Y_j^{b-1})).$$ (i.e., the coefficients $h_i$ are, up to
sign, the symmetric polynomials evaluated at the elements
$a_0+a_1Y_j+\dots a_{b-1}Y_j^{b-1}$).
\end{lemma}
\proof

The proof follows from the observations:

1) $\psi_{\Theta}(\Theta)=0.$

2) There is an isomorphism of  $k[s_{b,1},\ldots,s_{b,b}]$-modules, say $$\beta_j : k[s_{b,1},\ldots,s_{b,b}][Y_1]\to
k[s_{b,1},\ldots,s_{b,b}][Y_j]\hskip 1cm \beta_j(Y_1)=Y_j.$$ So
the characteristic polynomial of $\Theta$ in
$k[s_{b,1},\ldots,s_{b,b}][Y_1]$, is the same as that of
$\beta_j(\Theta)$ in $k[s_{b,1},\ldots,s_{b,b}][Y_j]$.
\endproof

Note that $ \mathbb{S}_b$ is the Galois group of the extension
$k[s_{b,1},\ldots,s_{b,b}]\subset k[Y_1,\dots ,Y_b]$. An element
$\sigma \in\mathbb{S}_b$, such that $\sigma(Y_1)=Y_j$, induces
$\beta_j$ when restricted to  $k[s_{b,1},\ldots,s_{b,b}][Y_1]$.

\begin{corollary} \label{cor012} If  $\Theta= a_0+a_1Y_1+\dots a_{b-1}Y_1^{b-1}$
is weighted homogeneous in $k[s_{b,1},\ldots,s_{b,b}][Y_1]$ (i.e., if
$\Theta$ is homogeneous as element of $k[Y_1,\dots ,Y_b]$), then
the coefficients of the characteristic polynomial
$\psi_{\Theta}(V)$ are weighted homogeneous in
$k[s_{b,1},\ldots,s_{b,b}]$ (i.e., are homogeneous in $k[Y_1,\dots
,Y_b]$).
\end{corollary}

In fact, the action of $ \mathbb{S}_b$ preserves degrees in
$k[Y_1,\dots ,Y_b]$, and each $\beta_j$ is a restriction of an
element in $\mathbb{S}_b$.

Let $F_b^{(e)}(Y_1)$ denote the class of $\Delta^e(F_b)(Z)$ in
$k[s_{b,1},\ldots,s_{b,b}][Y_1]$. Namely, set  $F_b^{(e)}(Y_1)=\Delta^e(F_b)(Y_1)$ )(\ref{par04}).
\begin{definition} \label{def013} Let $H_{F_b}$ be the $k$-subalgebra of $k[s_{b,1},\ldots,s_{b,b}]$
generated by the coefficients of the $b-1$ characteristic
polynomials, say $\psi_{F_b^{(e)}(Y_1)}(V),  for \hskip 0.1cm 1\leq e \leq
b-1.$
\end{definition}

\begin{lemma} \label{lemma014} $H_{F_b}$ is included in $k[Y_1-Y_2, Y_1-Y_3,\dots, Y_1-Y_b]$, and it is a graded subalgebra of
$k[s_{b,1},\ldots,s_{b,b}]$.
\end{lemma}
\proof

Recall that $F_b^{(e)}(Y_1)=(-1)^{b-e}t_{b-1,b-e}(Y_1-Y_j)$
(\ref{lemma2}), and note that $t_{b-1,b-e}(Y_1-Y_j)$ is
homogeneous in $k[Y_1-Y_2, Y_1-Y_3,\dots, Y_1-Y_b]$. The
coefficients of the characteristic polynomial of $F_b^{(e)}(Y_1)$
are symmetric polynomials on $\{F_b^{(e)}(Y_j) / 1\leq j \leq
b\}$, so they are also homogeneous elements in $k[Y_1-Y_2,
Y_1-Y_3,\dots, Y_1-Y_b]$. The second part of the claim follows
from the corollary.

\begin{proposition}\label{prop14} $\overline{R}_b$ (\ref{eqerrebarra}) is the integral closure of the graded ring
$H_{F_b}$ in $k[s_{b,1},\ldots,s_{b,b}]$.
\end{proposition}
\proof Observe that $(k[Y_1-Y_2, Y_1-Y_3,\dots,
Y_1-Y_b])^{ \mathbb{S}_{b}}=((k[Y_1-Y_2, Y_1-Y_3,\dots,
Y_1-Y_b])^{ \mathbb{S}_{b-1}})^{ \mathbb{S}_b}$,
 and recall that:

a) $(k[Y_1-Y_2, Y_1-Y_3,\dots, Y_1-Y_b])^{
\mathbb{S}_{b-1}}=k[F_b^{(1)}(Y_1),\dots ,
F_b^{(b-1)}(Y_1)](\subset k[s_{b,1},\ldots,s_{b,b}][Y_1]);$

b)$(k[Y_1-Y_2, Y_1-Y_3,\dots, Y_1-Y_b])^{
\mathbb{S}_{b}}=k[F_b^{(1)}(Y_1),\dots , F_b^{(b-1)}(Y_1)]\cap
k[s_{b,1},\ldots,s_{b,b}],$ as subrings of
$k[s_{b,1},\ldots,s_{b,b}][Y_1]$ (\ref{cor1}).

 Lemma \ref{lemma014} shows that the coefficients
of $\psi_{F_b^{(e)}(Y_1)}(V)$ are in the ring $\overline{R}_b$. In
fact, in $$k[Y_1-Y_2, Y_1-Y_3,\dots, Y_1-Y_b]\cap
k[s_{b,1},\ldots,s_{b,b}](= \overline{R}_b).$$

Finally use b) to show that $H_{F_b}\subset k[F_b^{(1)}(Y_1),\dots
, F_b^{(b-1)}(Y_1)]$  and notice that this is a finite ring
extension. In fact, the elements $F_b^{(e)}(Y_1)$ are integral
over $H_{F_b}$ since they satisfy the characteristic polynomial (i.e., they are roots of their own characteristic polynomials).

The claim follows now from:

i) $H_{F_b}\subset k[F_b^{(1)}(Y_1),\dots , F_b^{(b-1)}(Y_1)]\cap
k[s_{b,1},\ldots,s_{b,b}]$,

ii) $k[F_b^{(1)}(Y_1),\dots , F_b^{(b-1)}(Y_1)]$ is a finite
extension of $H_{F_b}$, and

iii) $k[F_b^{(1)}(Y_1),\dots , F_b^{(b-1)}(Y_1)]\cap
k[s_{b,1},\ldots,s_{b,b}]=\overline{R}_b$ is integrally closed in
$k[s_{b,1},\ldots,s_{b,b}]$.

As for iii): it was indicated in (\ref{rk09}) that $ \overline{R}_b$ is normal. Note that the quotient field $ \overline{R}_b$ is a subfield
of $k(Y_1-Y_2, Y_1-Y_3,\dots, Y_1-Y_b)$, whereas the quotient field 
of $k[s_{b,1},\ldots,s_{b,b}]$ is a subfield of 
$$k(Y_1, Y_2, Y_3,\dots, Y_b)=k(Y_1-Y_2, Y_1-Y_3,\dots, Y_1-Y_b)(Y_1),$$ which is purely transcendental over the latter.
\endproof

\centerline{\bf On multi-graded structures.} We have studied
invariants of a monic polynomial, motivated by elimination theory,
which were related to differential operators. However, elimination
applies also to more then one polynomial. An example of this is the
resultant of two polynomial. Here we generalize the previous
discussion to the case of several polynomials. This will lead us to
different graded rings, and also to extensions of graded rings. Of
special interest for further applications is the case in which
these ring extensions are integral (i.e., finite extensions).
\begin{parrafo}\label{rr1}{\em Fix positive integers $c_1,\dots,
c_r$, and set $b=\sum_ic_i$. For each index $i$ define
$$F_{c_i}=(Z-Y^{(i)}_1)\cdot (Z-Y^{(i)}_2)\cdots
(Z-Y^{(i)}_{c_i})\in k[Y^{(i)}_1,\dots ,Y^{(i)}_{c_i}][Z].$$ The
product of these polynomials is a polynomial of degree $b$ in $Z$,
say:
$$F_{c_1}(Z)\cdot F_{c_2}(Z)\cdots F_{c_r}(Z)=F_{b}(Z),$$
as polynomial in
$k[Y^{(1)}_1,\dots,Y^{(1)}_{c_1},Y^{(2)}_1,\dots,Y^{(2)}_{c_2},\dots
,Y^{(r)}_1,\dots,Y^{(r)}_{c_r}][Z]$, which we naturally identify
with $k[Y_1,\dots , Y_b][Z]$, as $c_1+\cdots+c_r=b$.

In other words, identify
$k[Y^{(1)}_1,\dots,Y^{(1)}_{c_1},Y^{(2)}_1,\dots,Y^{(2)}_{c_2},\dots
,Y^{(r)}_1,\dots,Y^{(r)}_{c_r}]$ with $k[Y_1,\dots , Y_b]$. The
permutation group $ \mathbb{S}_{c_i}$ acts on $k[Y^{(i)}_1,\dots
,Y^{(i)}_{c_i}]$, and there is an inclusion $ \mathbb{S}_{c_1}
\times
 \mathbb{S}_{c_2}\times \dots \times  \mathbb{S}_{c_r}$ in $ \mathbb{S}_b$. The main arguments
rely on two simple observations. First, that the permutation
groups $ \mathbb{S}_{c_i}$ and $ \mathbb{S}_b$ act linearly on
$k[Y^{(i)}_1,\dots ,Y^{(i)}_{c_i}]$ and  $k[Y^{}_1,\dots
,Y^{}_{b}]$ respectively, which asserts that the invariant rings
are graded. Second, that the inclusion of finite groups $
\mathbb{S}_{c_1} \times
 \mathbb{S}_{c_2}\times \dots \times  \mathbb{S}_{c_r}$ in $ \mathbb{S}_b$ provides a finite
extension of invariant rings.}

\end{parrafo}

\begin{parrafo}\label{qq1}
{\em The permutation group $ \mathbb{S}_{c_i}$ acts on $k[Y^{(i)}_1,\dots
,Y^{(i)}_{c_i}]$, and the inclusion $ \mathbb{S}_{c_1} \times  \mathbb{S}_{c_2}\times
\dots \times  \mathbb{S}_{c_r}$ in $ \mathbb{S}_b$ is such that: $$ (k[Y^{(1)}_1,\dots
,Y^{(1)}_{c_1}, \dots ,Y^{(r)}_1,\dots ,Y^{(r)}_{c_r}])^{ \mathbb{S}_{c_1}
\times  \mathbb{S}_{c_2}\times \dots \times  \mathbb{S}_{c_r} }=$$
$$=\otimes_ik[Y^{(i)}_1,\dots ,Y^{(i)}_{c_i}]^{ \mathbb{S}_{c_i}}.$$

The finite group $ \mathbb{S}_{c_1} \times  \mathbb{S}_{c_2}\times \dots \times
 \mathbb{S}_{c_r}$ also acts on the graded subalgebra
\begin{equation}\label{eqdeA}
A_{c_1,\dots, c_r}=k[Y^{(1)}_2-Y^{(1)}_1,\dots
,Y^{(1)}_{c_1}-Y^{(1)}_1, \dots ,Y^{(r)}_1-Y^{(1)}_1,\dots
,Y^{(r)}_{c_r}-Y^{(1)}_1]
\end{equation}
in a way that preserves the usual degrees. Therefore, the ring of
invariants, say
\begin{equation}\label{eqdeR}
\overline{R}_{c_1,\dots, c_r}=A_{c_1,\dots, c_r}^{ \mathbb{S}_{c_1} \times
 \mathbb{S}_{c_2}\times \dots \times  \mathbb{S}_{c_r}}
\end{equation}
 is a finitely
generated $k$-algebra, and a graded subring of $k[Y^{(1)}_1,\dots
,Y^{(1)}_{c_1}, \dots ,Y^{(r)}_1,\dots ,Y^{(r)}_{c_r}]$. Set
$$k[Y^{(i)}_1,\dots
,Y^{(i)}_{c_i}]^{ \mathbb{S}_{c_i}}=k[s^{(i)}_1,\dots,s^{(i)}_{c_i}].$$
 The graded rings $k[s^{(1)}_1,\dots,s^{(1)}_{c_1}, \dots ,s^{(r)}_1,\dots,s^{(r)}_{c_r}
 ]$ and $A_{c_1,\dots,
c_r}$ are both polynomial rings, and graded subrings of
$k[Y^{(1)}_1,\dots ,Y^{(1)}_{c_1}, \dots ,Y^{(r)}_1,\dots
,Y^{(r)}_{c_r}]$. As both are normal, so is
$$\overline{R}_{c_1,\dots, c_r}=A_{c_1,\dots, c_r} \cap
k[s^{(1)}_1,\dots,s^{(1)}_{c_1}, \dots
,s^{(r)}_1,\dots,s^{(r)}_{c_r}
 ]. $$

Since $\overline{R}_{c_1,\dots, c_r}$ is graded, it is generated
by weighted homogeneous polynomials, say
\begin{equation}\label{1ecdegl}
G_{c_1,\dots, c_r}^{l}=G_{c_1,\dots,
c_r}^{l}(s^{(1)}_1,\dots,s^{(1)}_{c_1}, \dots
,s^{(r)}_1,\dots,s^{(r)}_{c_r}), \hskip 0.5cm 1\leq l \leq
n_{c_1,\dots, c_r},
\end{equation}for some positive integer $n_{c_1,\dots, c_r}$.}
\end{parrafo}

\begin{theorem}
\label{vramloc}  Let $R$ be a $k$-algebra, and let
$Spec(R[Z]/\langle f(Z) \rangle
)\stackrel{\pi}{\rightarrow}Spec(R)$
 be the natural
finite morphism, where $$f(Z)=f_1(Z)\cdot f_2(Z)\cdots f_r(Z),
\hskip 0.5cm
f_i(Z)=Z^c+a^{(i)}_1Z^{c-1}+\ldots+a^{(i)}_{c_i-1}Z+a^{(i)}_{c_i}\in
R[Z],$$ and $c_1+c_2+\cdots +c_r=b$. Then:

i) $V(\langle G_{c_1,\dots,
c_r}^{l}(a^{(1)}_1,\dots,a^{(1)}_{c_1}, \dots
,a^{(r)}_1,\dots,a^{(r)}_{c_r})/ 1\leq l \leq n_{c_1,\dots, c_r}
\rangle) $ is the set of points in $ \mbox{Spec}(R)$ where the
finite morphism is purely ramified.

ii) If $R$ is regular, and $Q\in V(\langle f(Z) \rangle)$ is a
point of multiplicity $b$ of this hypersurface, then
$$\nu_S(G_{c_1,\dots, c_r}^{l}(a^{(1)}_1,\dots,a^{(1)}_{c_1},
\dots ,a^{(r)}_1,\dots,a^{(r)}_{c_r})) \geq deg \ G_{c_1,\dots,
c_r}^{l} $$ for $1\leq l \leq n_{c_1,\dots, c_r}^{}$, where
$S=R_P$, $P=\pi(Q)$.
\end{theorem}

\proof i) Set $f(Z)=Z^b+h_1Z^{b-1}+\ldots+h_{b-1}Z+h_b\in R[Z]$.
According to Theorem \ref{ramlocus}, the purely ramified locus is
the closed set in $Spec(R)$ defined by the ideal spanned by all
$G_i(h_1,\dots, h_b)$, where the $G_i$ are homogeneous generators
of the subring of $ \mathbb{S}_b$ invariants in $A_{c_1,\dots,
c_r}$ in (\ref{eqdeA}).

Here the elements $G_{c_1,\dots, c_r}^{l}$ span the subring of $
\mathbb{S}_{c_1} \times  \mathbb{S}_{c_2}\times \dots \times
\mathbb{S}_{c_r}$-invariants (\ref{1ecdegl}). The inclusion $
\mathbb{S}_{c_1} \times  \mathbb{S}_{c_2}\times \dots \times
\mathbb{S}_{c_r}\subset  \mathbb{S}_b$ defines a finite extension
of invariant rings, and both are generated by homogeneous
elements, say $k[{\{G_i\}}_{1\leq i \leq n_b}] \subset
k[{\{G_{c_1,\dots, c_r}^{l}\}}_{1\leq l \leq n_{c_1,\dots , c_r}}
].$

So the ideal in $k[{\{G_{c_1,\dots, c_r}^{l}\}}_{1\leq l \leq
n_{c_1,\dots , c_r}}]$ spanned by all elements $G_i$, is included,
and has the same radical, as that spanned by the $G_{c_1,\dots,
c_r}^{l}$'s (\ref{grocal}). As it was indicated in \ref{rkperma} this
property is preserved by arbitrary homomorphisms of $k$-algebras.
This proves (i).

ii) Follows by the same argument used in Theorem \ref{ramlocus}
ii).
\endproof
\begin{parrafo}\label{vpp2p2}
{\em We have defined  $F_b(Z)=(Z-Y_1)\cdot (Z-Y_2)\cdots
(Z-Y_b)\in k[Y_1,\dots ,Y_b][Z]$. In \ref{rkjun06} it was shown
that the coefficients of
$$Tay(F(Z))=F_b(T+Z)=(T+(Z-Y_1))\cdot (T+(Z-Y_2))\cdots
(T+(Z-Y_b))\in k[Y_1,\dots ,Y_b][Z,T],$$ in the variable $T$, are
the symmetric polynomials evaluated on the elements $Z-Y_j$,
$1\leq j\leq b$. Namely, $ F_b^{(e)}(Z)=(-1)^{b-e}s_{b,b-e}(Z-Y_1,
Z-Y_2,\dots , Z-Y_b), $ where $F_b^{(e)}(Z)$ denotes
$\Delta^{(e)}(F(Z))$(see \ref{eqjuni06}), and
$$k[Z-Y_1,\dots , Z-Y_b]^{ \mathbb{S}_b}=k[\{F_b^{(e)}(Z), e=0,1,\dots, b-1\}].$$

Fix, as in \ref{rr1}, positive integers $c_1, \dots ,c_r$ so that
$c_1+\cdots+c_r=b$, and set, for each index $i$,
$F_{c_i}(Z)=\prod_{1\leq j \leq c_i} (Z-Y^{(i)}_j)\in
k[Y^{(i)}_1,\dots,Y^{(i)}_{c_i}][Z]$. The product of these
polynomials is of degree $b$ in $Z$, say:
$$F_{c_1}(Z)\cdot F_{c_2}(Z)\cdots F_{c_r}(Z)=F_{b}(Z),$$
 in
$k[Y^{(1)}_1,\dots,Y^{(1)}_{c_1},Y^{(2)}_1,\dots,Y^{(2)}_{c_2},\dots
,Y^{(r)}_1,\dots,Y^{(r)}_{c_r}][Z]$, which we identified with
$k[Y_1,\dots , Y_b][Z]$, as $c_1+\cdots+c_r=b$. In fact, we identify
$k[Y^{(1)}_1,\dots,Y^{(1)}_{c_1},Y^{(2)}_1,\dots,Y^{(2)}_{c_2},\dots
,Y^{(r)}_1,\dots,Y^{(r)}_{c_r}]$ with $k[Y_1,\dots , Y_b]$, and
hence:
$$k[Z-Y_1,\dots , Z-Y_b]=k[\{Z-Y^{(i)}_j\} 1\leq i \leq r, 1\leq j\leq c_i].$$
Set, as before, $F_b^{(\alpha)}(Z)=\Delta^{(\alpha)}(F_b(Z))$ and
$F_{c_i}^{(\alpha)}(Z)=\Delta^{(\alpha)}(F_{c_i}(Z))$(
\ref{par04}).

}
\end{parrafo}

\begin{proposition}\label{propjun16} 1) For $0\leq \alpha \leq b-1$, and $0\leq \beta \leq c_i-1$, each $F_b^{(\alpha)}(Z)$ is homogeneous of
degree $b-\alpha$, and each $F_{c_i}^{(\beta)}(Z)$ is homogeneous
of degree $c_i-\beta$ in the graded ring
$k[Y^{(1)}_1,\dots,Y^{(1)}_{c_1},Y^{(2)}_1,\dots,Y^{(2)}_{c_2},\dots
,Y^{(r)}_1,\dots,Y^{(r)}_{c_r}][Z]=k[Y_1,\dots , Y_b][Z]$.

2) There is a finite inclusion of graded subrings of $k[Y_1,\dots
, Y_b][Z]$, defined by
$$k[\{F_b^{(e)}(Z), e=0,1,\dots, b-1\}]\subset k[\{F_{c_i}^{(\beta)}(Z),1\leq i\leq r, 0\leq \beta \leq c_i\} ].$$
3) The two algebras in 2) are also graded subalgebras in
$k[s^{(1)}_1,\dots,s^{(1)}_{c_1},\dots,s^{(r)}_1,\dots,
s^{(r)}_{c_r}][Z]$.
\end{proposition}

\proof

1) follows from \ref{eqjuni06}, and 2) from the fact that such
rings are the invariants by the finite groups $ \mathbb{S}_{c_1} \times
 \mathbb{S}_{c_2}\times \dots \times  \mathbb{S}_{c_r}\subset  \mathbb{S}_b$; acting on
$k[Z-Y_1,\dots , Z-Y_b]=k[\{Z-Y^{(i)}_j\} 1\leq i \leq r, 1\leq
j\leq c_i].$
\endproof
\begin{proposition} \label{propjun17} The graded algebra $\overline{R}_{c_1,\dots, c_r}$ in \ref{eqdeR} is a graded subring of
$k[\{F_{c_i}^{(\beta)}(Z),1\leq i\leq r, 0\leq \beta \leq c_i\}]$.
\end{proposition}
\proof

Note that $A_{c_1,\dots, c_r}$, in \ref{eqdeA} can be expressed as
a subring of $k[\{Z-Y^{(i)}_j; 1\leq i \leq r, 1\leq j\leq c_i\}]$
by setting $Y_j^{(i)}-Y_1^{(1)}=(Z-Y^{(1)}_1)-(Z-Y^{(i)}_j).$

\begin{corollary} \label{corjun17}Let $H_N\in \overline{R}_{c_1,c_2,\dots,c_r}$ be homogeneous
of degree $N$. Then there is a polynomial in variables $\{
W_{c_1}, W_{c_1}^{(1)}, W_{c_1}^{(2)},\dots
,W_{c_1}^{(c_1-1)},\dots , W_{c_r}^{}, W_{c_r}^{(1)},\dots
,W_{c_r}^{(c_r-1)} \},$ and coefficients in $k$, say
$$G\in k[  W_{c_1}, W_{c_1}^{(1)}, W_{c_1}^{(2)},\dots
,W_{c_1}^{(c_1-1)},\dots , W_{c_r}^{}, W_{c_r}^{(1)},\dots
,W_{c_r}^{(c_r-1)}],$$ such that \tiny $$G(F_{c_1}(Z),
F_{c_1}^{(1)}(Z), F_{c_1}^{(2)}(Z),\dots
,F_{c_1}^{(c_1-1)}(Z),\dots , F_{c_r}(Z), F_{c_r}^{(1)}(Z),
F_{c_r}^{(2)}(Z),\dots ,F_{c_r}^{(c_r-1)}(Z))=H_N.$$ \normalsize
Furthermore, we may assume that $G$ is weighted homogeneous of
degree $N$, if we assign degree $c_i-j$ to the variable
$W_{c_i}^{(j)}$.
\end{corollary}

\centerline{\bf Multi-graded structures and elimination.}

\begin{parrafo}\label{vpp2}
{\em We now extend Corollary \ref{cor1} and Proposition
\ref{prop14} to the context of several monic polynomials. Notice
that
$$F_{c_i}(Z)=\prod_{1\leq j \leq c_i} (Z-Y^{(i)}_j)\in
k[Y^{(i)}_1,\dots,Y^{(i)}_{c_i}][Z]$$ has coefficients in
$k[Y^{(i)}_1,\dots ,Y^{(i)}_{c_i}]^{
\mathbb{S}_{c_i}}=k[s^{(i)}_1,\dots,s^{(i)}_c]$; in fact
$$F_{c_i}(Z)=Z^{c_i}+(-1)s^{(i)}_1Z^{c_i-1}+\ldots+(-1)^{c_i-1}s^{(i)}_{c_i-1}Z+(-1)^{c_i}
s^{(i)}_{c_i}.$$

Set $i=1$ and let $ \mathbb{S}_{c_1-1}$ denote the subgroup of $ \mathbb{S}_{c_1}$,
consisting of those elements fixing $Y_1^{(1)}$. Observe that
$$k[s^{(1)}_1,\dots,s^{(1)}_{c_1}][Y^{(1)}_1]=k[Y^{(1)}_1,\dots
,Y^{(1)}_{c_1}]^{ \mathbb{S}_{c_1-1}}$$ can be identified with
$k[s^{(1)}_1,\dots,s^{(1)}_c][Z]/\langle F_{c_1}(Z)\rangle.$

On the other hand $$k[s^{(1)}_1,\dots,s^{(1)}_{c_1},\dots
s^{(r)}_1,\dots,s^{(r)}_{c_r}][Y^{(1)}_1]=$$ $$=
(k[Y^{(1)}_1,\dots ,Y^{(1)}_{c_1},\dots ,Y^{(r)}_1,\dots
,Y^{(r)}_{c_r}])^{ \mathbb{S}_{c_1-1}\times  \mathbb{S}_{c_2}
\cdots \times  \mathbb{S}_{c_r} }$$ can be identified with
$$k[s^{(1)}_1,\dots,s^{(1)}_{c_1},\dots
s^{(r)}_1,\dots,s^{(r)}_{c_r}][Z]/\langle F_{c_1}(Z)\rangle.$$ The
natural inclusion $A_{c_1,\dots, c_r}\subset k[Y^{(1)}_1,\dots
,Y^{(1)}_{c_1},\dots ,Y^{(r)}_1,\dots ,Y^{(r)}_{c_r}] $ (see
(\ref{eqdeA}))
 shows that
$$A_{c_1,\dots, c_r}^{ \mathbb{S}_{c_1-1}\times  \mathbb{S}_{c_2} \cdots \times
 \mathbb{S}_{c_r}}\subset  k[s^{(1)}_1,\dots,s^{(1)}_{c_1},\dots
s^{(r)}_1,\dots,s^{(r)}_{c_r}][Z]/\langle F_{c_1}(Z)\rangle.$$}
\end{parrafo}

Let $F^{(\alpha)}(Z)$ denote $\Delta^{(\alpha)}(F(Z))$ (as in
\ref{par04}).

\begin{lemma} \label{1lemdederi} The ring
$$A_{c_1,\dots, c_r}^{ \mathbb{S}_{c_1-1}\times  \mathbb{S}_{c_2} \cdots \times
 \mathbb{S}_{c_r}}=k[Y^{(1)}_2-Y^{(1)}_1,\dots ,Y^{(1)}_{c_1}-Y^{(1)}_1,
\dots ,Y^{(r)}_1-Y^{(1)}_1,\dots
,Y^{(r)}_{c_r}-Y^{(1)}_1]^{ \mathbb{S}_{c_1-1}\times  \mathbb{S}_{c_2} \cdots \times
 \mathbb{S}_{c_r}}$$
 is a graded $k$-algebra and a subring of
 $$k[s^{(1)}_1,\dots,s^{(1)}_{c_1},\dots
s^{(r)}_1,\dots,s^{(r)}_{c_r}][Y^{(1)}_1]=k[s^{(1)}_1,\dots,s^{(1)}_{c_1},\dots
s^{(r)}_1,\dots,s^{(r)}_{c_r}][Z]/\langle F_{c_1}(Z) \rangle, $$
generated by the class of the elements $$\{ F_{c_1}(Z),
F_{c_1}(Z)^{(1)}, F_{c_1}(Z)^{(2)},\dots
,F_{c_1}(Z)^{(c_1-1)},\dots ,   F_{c_r}(Z), F_{c_r}(Z)^{(1)},
F_{c_r}(Z)^{(2)},\dots ,F_{c_r}(Z)^{(c_r-1)} \}$$ in this quotient
ring. In other words,\tiny $$A_{c_1,\dots, c_r}^{ \mathbb{S}_{c_1-1}\times
 \mathbb{S}_{c_2} \cdots \times  \mathbb{S}_{c_r}}=$$ $$=k[F_{c_1}(Y^{(1)}_1),
F_{c_1}^{(1)}(Y^{(1)}_1), F_{c_1}^{(2)}(Y^{(1)}_1),\dots
,F_{c_1}^{(c_1-1)}(Y^{(1)}_1),\dots , F_{c_r}(Y^{(1)}_1),
F_{c_r}^{(1)}(Y^{(1)}_1), F_{c_r}^{(2)}(Y^{(1)}_1),\dots
,F_{c_r}^{(c_r-1)}(Y^{(1)}_1)].$$\normalsize
\end{lemma}

\begin{corollary}\label{cordecmul}

1) $\overline{R}_{c_1,c_2,\dots,c_r}\subset$ \tiny $$\subset
k[F_{c_1}(Y^{(1)}_1), F_{c_1}^{(1)}(Y^{(1)}_1),
F_{c_1}^{(2)}(Y^{(1)}_1),\dots ,F_{c_1}^{(c_1-1)}(Y^{(1)}_1),\dots
, F_{c_r}(Y^{(1)}_1), F_{c_r}^{(1)}(Y^{(1)}_1),
F_{c_r}^{(2)}(Y^{(1)}_1),\dots ,F_{c_r}^{(c_r-1)}(Y^{(1)}_1)]$$
\normalsize is a finite extension of graded rings.

2)Let $H_N\in \overline{R}_{c_1,c_2,\dots,c_r}$ be homogeneous of
degree $N$. There is a polynomial in variables $$\{ W_{c_1},
W_{c_1}^{(1)}, W_{c_1}^{(2)},\dots ,W_{c_1}^{(c_1-1)},\dots ,
W_{c_r}^{}, W_{c_r}^{(1)},\dots ,W_{c_r}^{(c_r-1)} \},$$ say $G\in
k[  W_{c_1}, W_{c_1}^{(1)}, W_{c_1}^{(2)},\dots
,W_{c_1}^{(c_1-1)},\dots , W_{c_r}^{}, W_{c_r}^{(1)},\dots
,W_{c_r}^{(c_r-1)}],$ such that \tiny $$G(F_{c_1}(Y^{(1)}_1),
F_{c_1}^{(1)}(Y^{(1)}_1), F_{c_1}^{(2)}(Y^{(1)}_1),\dots
,F_{c_1}^{(c_1-1)}(Y^{(1)}_1),\dots , F_{c_r}(Y^{(1)}_1),
F_{c_r}^{(1)}(Y^{(1)}_1), F_{c_r}^{(2)}(Y^{(1)}_1),\dots
,F_{c_r}^{(c_r-1)}(Y^{(1)}_1))=H_N.$$ \normalsize Furthermore, we
may assume that $G$ is weighted homogeneous of degree $N$, if we
assign weight (degree) $c_i-j$ to the variable $W_{c_i}^{(j)}$.
\end{corollary}
\begin{proof} 1) Both rings are defined as subrings of invariants of the finite groups
$$\mathbb{S}_{c_1-1}\times  \mathbb{S}_{c_2} \times \cdots \times  \mathbb{S}_{c_r}\subset
 \mathbb{S}_{c_1}\times  \mathbb{S}_{c_2} \times \cdots \times  \mathbb{S}_{c_r} $$ acting
linearly on $A_{c_1,\dots ,c_r}$ (see \ref{eqdeR} and
\ref{1lemdederi} ).

2) Same argument as in \ref{corjun17}. It follows from 1) and the
fact that each $ F_{c_i}^{(j)}(Y_1)$ is homogeneous of degree
$c_i-j$.
\end{proof}
\begin{remark} If, instead of $ \mathbb{S}_{c_1-1}\times  \mathbb{S}_{c_2} \cdots \times
 \mathbb{S}_{c_r}$  we consider $ \mathbb{S}_{c_1}\times  \mathbb{S}_{c_2-1}\times  \mathbb{S}_{c_3}
\cdots \times  \mathbb{S}_{c_r} $  , the same argument shows that
$A_{c_1,\dots, c_r}^{ \mathbb{S}_{c_1}\times
\mathbb{S}_{c_2-1}\times  \mathbb{S}_{c_3} \cdots \times
\mathbb{S}_{c_r} }$ is a graded $k$-algebra and a subring of
$$k[s^{(1)}_1,\dots,s^{(1)}_{c_1},\dots
s^{(r)}_1,\dots,s^{(r)}_{c_r}][Y^{(2)}_1]=k[s^{(1)}_1,\dots,s^{(1)}_{c_1},\dots
s^{(r)}_1,\dots,s^{(r)}_{c_r}][Z]/\langle F_{c_2}(Z)\rangle, $$
generated by the class of the elements $$\{ F_{c_1}(Z),
F_{c_1}(Z)^{(1)}, F_{c_1}(Z)^{(2)},\dots
,F_{c_1}(Z)^{(c_1-1)},\dots ,   F_{c_r}(Z), F_{c_r}(Z)^{(1)},
F_{c_r}(Z)^{(2)},\dots ,F_{c_r}(Z)^{(c_r-1)} \}$$ in this quotient
ring.
\end{remark}

\centerline{\bf Appendix 2: On normality of multi-graded invariant
ring.}

In this part we extend Lemma \ref{lemma014} and Proposition
\ref{prop14} to the case of several polynomials. Proofs are
straightforward generalizations.
\begin{parrafo}\label{mulpar010bis}
{\em The ring
$k[s^{(1)}_1,\dots,s^{(1)}_{c_1},\dots,s^{(r)}_1,\dots,
s^{(r)}_{c_r}][Y^{(1)}_1]$ is a free module of rank $b$ over

$k[s^{(1)}_1,\dots,s^{(1)}_{c_1},\dots,s^{(r)}_1,\dots,
s^{(r)}_{c_r}]$. Multiplication by an element $$\Theta \in
k[s^{(1)}_1,\dots,s^{(1)}_{c_1},\dots,s^{(r)}_1,\dots,
s^{(r)}_{c_r}][Y^{(1)}_1]$$ defines an endomorphism, say
$\phi_{\Theta}$, with characteristic polynomial, say
$$\psi_{\Theta}(V)=V^b+h_1V^{b-1}+\dots +h_b\in
k[s^{(1)}_1,\dots,s^{(1)}_{c_1},\dots,s^{(r)}_1,\dots,
s^{(r)}_{c_r}][V].$$}
\end{parrafo}

\begin{definition} \label{muldef013bis} Let $H_{F_{c_1},F_{c_2},\dots ,F_{c_r}}$ be the $k$-subalgebra of
$k[s^{(1)}_1,\dots,s^{(1)}_{c_1},\dots,s^{(r)}_1,\dots,
s^{(r)}_{c_r}]$ generated by the coefficients of the
characteristic polynomials $$\psi_{F_{c_i}^{(e)}(Y_1)}(V), \hskip
0.5cm 1\leq e \leq c_i-1, \hskip 0.25cm 1\leq i \leq r.$$
\end{definition}

\begin{lemma} \label{mullemma014bis} $H_{F_{c_1},F_{c_2},\dots ,F_{c_r}}$ is graded subalgebra both of
$A_{c_1,c_2,\dots,c_r}$ (\ref{eqdeA}), and of

$k[s^{(1)}_1,\dots,s^{(1)}_{c_1},\dots,s^{(r)}_1,\dots,
s^{(r)}_{c_r}]$.
\end{lemma}

\begin{proposition}\label{mulprop14bis} $\overline{R}_{c_1,c_2,\dots,c_r}$ is the integral closure of the graded ring
$H_{F_{c_1},F_{c_2},\dots ,F_{c_r}}$  in

$k[s^{(1)}_1,\dots,s^{(1)}_{c_1},\dots,s^{(r)}_1,\dots,
s^{(r)}_{c_r}]$.
\end{proposition}
\centerline{\bf Keeping track of the grade (the variable W), and
elimination algebras.}
\begin{parrafo}\label{lavw}{\em Let $G$ be an $\mathbb{N}$-graded ring. Add a variable $W$, say $G[W]$, and grade this new algebra with
$W$ of degree one.  We define now what we call a {\em graded
inclusion}, say $G\to G[W]$, as follows: If $G=\sum_{k\geq 0}
[G]_k $, define $\sum_{k\geq 0} [G]_k W^k$ as subalgebra in
$G[W]$. Note that if $G=\sum_{k\geq 0} [G]_k \subset
G'=\sum_{k\geq 0} [G']_k$ is a finite extension of
$\mathbb{N}$-graded algebras, then $\sum_{k\geq 0} [G]_k
W^k\subset G'=\sum_{k\geq 0} [G']_kW^k$ is also finite. If $G$ is
a $k$-algebra generated by elements, say $\{H_1,\dots, H_s\}$,
where each $H_i$ is homogeneous of degree $d_i$; then the graded
inclusion of $G$ is the $k$-subalgebra in $G[W]$ generated, over
$k$, by $\{H_1W^{d_1},\dots, H_sW^{d_1}\}$.

 }
\end{parrafo}

\begin{parrafo}\label{1univuno} {\em Set $T=k[s^{(1)}_1,\dots,s^{(1)}_{c_1},\dots,s^{(r)}_1,\dots,
s^{(r)}_{c_r}]$. We attach to the subalgebras
$\overline{R}_{c_1,c_2,\dots,c_r}$ and $H_{F_{c_1},F_{c_2},\dots
,F_{c_r}}$  two graded $k$-subalgebras included in a polynomial
ring $T[W]$, as in \ref{lavw}.

Now set $U^{(1)}_{c_1,c_2,\dots,c_r}$ (set
$U^{(2)}_{c_1,c_2,\dots,c_r}$) as the $T$-subalgebra of $T[W]$
generated by such $k$-subalgebra. In other words, the algebras
generated over $T$ by all elements of the form
$$H(s^{(1)}_1,\dots,s^{(1)}_{c_1},\dots,s^{(r)}_1,\dots,
s^{(r)}_{c_r})\cdot W^{d_H},$$ where $H$ is an homogeneous element
of degree $d_H$ in $\overline{R}_{c_1,c_2,\dots,c_r}$ (in
$H_{F_{c_1},F_{c_2},\dots ,F_{c_r}}$ ).

We can also express  $U^{(1)}_{c_1,c_2,\dots,c_r}=\sum I^{(1)}_k
W^k (\subset T[W])$, and $U^{(2)}_{c_1,c_2,\dots,c_r}=\sum
I^{(2)}_k W^k (\subset T[W])$, where $I^{(1)}_k$ and $I^{(2)}_k$
are ideals in $T$.

 Lemma \ref{mullemma014bis} ensures that $I^{(2)}_k\subset
I^{(1)}_k$ for each positive index $k$. Namely, that
$$U^{(2)}_{c_1,c_2,\dots,c_r}\subset
U^{(1)}_{c_1,c_2,\dots,c_r}.$$ Furthermore, Prop
\ref{mulprop14bis} asserts that this extension is finite.

Note that both algebras $U^{(1)}_{c_1,c_2,\dots,c_r}$ and
$U^{(2)}_{c_1,c_2,\dots,c_r}$ are finitely generated over $T$.
This follows from the fact that the graded algebras
$\overline{R}_{c_1,c_2,\dots,c_r}$ and $H_{F_{c_1},F_{c_2},\dots
,F_{c_r}}$ are finitely generated.

For example $U^{(1)}_{c_1,c_2,\dots,c_r}$ is the subring of $T[W]$
generated over $T$ by: $$\mathcal{F}=\{ G_{c_1,\dots,
c_r}^{l}(s^{(1)}_1,\dots,s^{(1)}_{c_1}, \dots
,s^{(r)}_1,\dots,s^{(r)}_{c_r})\cdot W^{ deg \ G_{c_1,\dots,
c_r}^{l}} , \hskip 0.5cm 1\leq l \leq n_{c_1,\dots, c_r}\} \hskip
0.5cm (\mbox{ see } \ref{1ecdegl}).$$}
\begin{remark} $\overline{R}_{c_1,\dots,c_r}$ and $H_{F_{c_1},F_{c_2},\dots ,F_{c_r}}$  are graded subalgebra both of
$A_{c_1,c_2,\dots,c_r}$ (\ref{eqdeA}), and of $T$, which in turn
is a graded ring in
$k[Y^{(1)}_1,\dots,Y^{(1)}_{c_1},\dots,Y^{(r)}_1,\dots,
Y^{(r)}_{c_r}].$ The algebras $U^{(2)}_{c_1,c_2,\dots,c_r}\subset
U^{(1)}_{c_1,c_2,\dots,c_r}$ are defined in $T[W]$, where the
variable $W$ keeps track of the degree.

\end{remark}

\centerline{\bf Specialization: From universal to concrete.}

\begin{definition}\label{1def29} Fix a $k$-algebra $S$ and monic
polynomials $$f_i(Z)= Z^{c_i}-a^{(i)}_1+\cdots +
(-1)^{c_i}a^{(i)}_{c_i}\in S[Z],\hskip 0.5cm  1\leq i \leq r.$$
Set $F_{c_1},F_{c_2},\dots,F_{c_r}$ and
$\overline{R}_{{c_1},{c_2},\dots,{c_r}}$ as before. So
$F_{c_i}(Z)\in k[s^{(i)}_1,\dots,s^{(i)}_{c_i}][Z]$, and each
$f_i(Z)$ is obtained from $F_{c_i}(Z)$ by the change of base rings
$\pi: k[s^{(i)}_1,\dots,s^{(i)}_{c_i}]\to S$ defined by setting
$\pi(s^{(i)}_j)=a^{(i)}_j$.

In this way the polynomials $f_i(Z)\in S[Z]$, $1\leq i \leq r$
define a morphism of $k$-algebras:
\begin{equation}\label{114422}T=k[s^{(1)}_1,\dots,s^{(1)}_{c_1},\dots,s^{(r)}_1,\dots,
s^{(r)}_{c_r}]\to S,
\end{equation} which extends to a morphism $T[W] \to
S[W]$. We define the {\em elimination algebra}, say
\begin{equation}\label{1eqdrbar}
\overline{\mathcal{R}}_{f_{c_1},f_{c_2},\dots,f_{c_r}} (\subset
S[W])
\end{equation}
as the subalgebra of $S[W]$ generated by the image of
$U^{(1)}_{c_1,c_2,\dots,c_r}=\sum I^{(1)}_k W^k (\subset T[W])$.

In the same way we define
\begin{equation}\label{1eqdrbar2}
\mathcal{H}_{f_{c_1},f_{c_2},\dots,f_{c_r}} (\subset S[W])
\end{equation}
as the subalgebra of $S[W]$ generated by the image of
$U^{(2)}_{c_1,c_2,\dots,c_r}=\sum I^{(2)}_k W^k (\subset T[W])$.

The previous observations show that
\begin{equation}\label{1eqdrbarmier}
\mathcal{H}_{f_{c_1},f_{c_2},\dots,f_{c_r}}\subset
\overline{\mathcal{R}}_{f_{c_1},f_{c_2},\dots,f_{c_r}}
\end{equation}
 and that
this ring extension is finite. Note also that both are finitely generated, for instance, the
elimination algebra
$\overline{\mathcal{R}}_{f_{c_1},f_{c_2},\dots,f_{c_r}}$ is
spanned over $S$ by the finite set $$\mathcal{F}=\{ G_{c_1,\dots,
c_r}^{l}(s^{(1)}_1,\dots,s^{(1)}_{c_1}, \dots
,s^{(r)}_1,\dots,s^{(r)}_{c_r})\cdot W^{ deg \ G_{c_1,\dots,
c_r}^{l}} , \hskip 0.5cm 1\leq l \leq n_{c_1,\dots, c_r}\} \hskip
0.5cm (\ref{1ecdegl}).$$

\end{definition}
\end{parrafo}
\begin{parrafo}\label{uni660} {\em
Set $T= k[s^{(1)}_1,\dots,s^{(1)}_{c_1},\dots,s^{(r)}_1,\dots,
s^{(r)}_{c_r}]$ as above. Prop \ref{propjun16} asserts that
there is a finite extension, say:
$k[\{F_b^{(e)}(Z), e=0,1,\dots, b-1\}]\subset
k[\{F_{c_i}^{(\beta)}(Z),1\leq i\leq r, 0\leq \beta \leq c_i \}]$
of graded $k$-subalgebras of (the graded algebra)
$T[Z]=k[s^{(1)}_1,\dots,s^{(1)}_{c_1},\dots,s^{(r)}_1,\dots,
s^{(r)}_{c_r}][Z]$. This defines a finite extension of $k$-subalgebras in $T[Z,W]$ (\ref{lavw}). Finally, as each such
$k$-algebra extends to an algebra over $T[Z]$, there is a finite extension of $T[Z]$-algebras:
\begin{equation}\label{eq660-1}
T[Z][\{F_b^{(e)}(Z)W^{b-e}, 0\leq e \leq b-1\}]\subset
T[Z][\{F_{c_i}^{(\beta)}(Z)W^{c_i-\beta},1\leq i\leq r, 0\leq
\beta \leq c_i \}]
\end{equation}
Consider, as in \ref{1def29}, a ring $S$ and monic polynomials
$f_{c_i}(Z)= Z^{c_i}-a^{(i)}_1+\cdots + (-1)^{c_i}a^{(i)}_{c_i}\in
S[Z],\hskip 0.5cm  1\leq i \leq r.$ 
Set $b=c_1+\dots +c_r$, and
$f_b(Z)=f_{c_1}(Z)\cdot f_{c_2}(Z)\cdots f_{c_r}(Z)$, which is
monic of degree $b$. The morphism in (\ref{114422}), and (\ref{eq660-1}),
define a finite inclusion
\begin{equation}\label{eq660}
S[Z][\{f_b^{(e)}(Z)W^{b-e}, 0\leq e \leq b-1\}]\subset
S[Z][\{f_{c_i}^{(\beta)}(Z)W^{c_i-\beta},1\leq i\leq r, 0\leq
\beta \leq c_i \}].
\end{equation}

Similarly, the inclusion in Prop \ref{propjun17} yields a (non-finite)
inclusion of the elimination algebra
\begin{equation}\label{eq6V60}
\overline{\mathcal{R}}_{f_{c_1},f_{c_2},\dots,f_{c_r}}\subset
S[Z][\{f_{c_i}^{(\beta)}(Z)W^{c_i-\beta},1\leq i\leq r, 0\leq
\beta \leq c_i \}].
\end{equation}
Where the left hand side is 
a subalgebra of $S[W]$, and the variable $Z$ has been eliminated.
 }
\end{parrafo}

%%%%%%%%%%%%%%%%%%%%%%%%%%%%%%%%%%%%%%%%%%%%%%%%%%%%%%%%%%%%%%%%%

%%%%%%%%%%%%%%%%%%%%%%%%%%%%%%%%%%%%%%%%%%%%%%%%%%%%%%%%%%%%%%%%%
\section{ Rees algebras and relative differential Rees algebras.}

In this section a class of $\mathbb{N}$-graded algebras is
introduced, which extend the notion of Rees rings of ideals. The
geometric application of these algebras will lead us to consider
them only up to integral closure (not to distinguish algebras with the same integral closure). Special attention will be devoted to
finite extensions among these algebras. In \ref{IIpar2} we begin
the study algebras with a natural action of differential
operators. 

The results in this Section appear in Sections 1 of \cite{VV1}. They are included here because of their relevance in all the further discussion.

\begin{parrafo}\label{IIpar1}  {\em Fix a noetherian ring $B$, and a sequence of ideals $\{I_k\}$, $k\geq 0$,
which fulfill the conditions:

{\bf i)} $I_0=B$, and

{\bf ii)} $I_k\cdot I_s \subset I_{k+s}$.

This defines a graded subring $\bigoplus_{k\geq 1} I_k\cdot W^k$
of the polynomial ring $ B[W]$.  We say that $\bigoplus I_k\cdot
W^k$ is a {\em Rees algebra} if this
subring is a (noetherian) finitely generated $B$-algebra.}

\end{parrafo}

\begin{remark}\label{IIrk1}
1) Examples of Rees algebras are the Rees rings of an ideal
$I\subset B$, where $I_k=I^k$ for each $k\geq 1$. These are the Rees algebras generated, as $B$-algebras,  in degree one.

2) Whenever $\bigoplus I_k\cdot W^k \subset (\subset B[W])$ is a
Rees algebra, define  $\bigoplus
I'_k\cdot W^k $ by setting $$ I'_k=\sum_{r\geq k} I_r.$$
If $\bigoplus I_k\cdot W^k  $ is generated by $\mathcal{F}=\{
g_{n_i}W^{n_i}, 1\leq i\leq m, n_i>0 \}$. Namely, if: $$ \bigoplus
I_k\cdot W^k =B[\{g_{n_i}W^{n_i}\}_{g_{n_i}W^{n_i}\in
\mathcal{F}}],$$ then $\bigoplus I'_k\cdot W^k $ is generated by
the finite set $\{ g_{n_i}W^{n'_i}, 1\leq i\leq m, 1 \leq n'_i\leq
n_i \}$.

Note that $I'_k \supset I'_{k+1}$, and that $\bigoplus I_k\cdot
W^k \subset \bigoplus I'_k\cdot W^k $ is a finite extension. It suffices to check that given an element $g\in I_{k}$,
then $g\cdot W^{k-1}$ is integral over $\bigoplus I_k\cdot W^k $.
Notice that $$g\in I_{k}\Rightarrow g^{k-1} \in
I_{k(k-1)}\Rightarrow g^{k} \in I_{k(k-1)},$$ so $g\cdot W^{k-1}$
fulfills the equation $Z^{k}-(g^{k}\cdot W^{k(k-1)})=0$.

Up to integral closure we may assume that a Rees algebra has
the additional condition:

{\bf iii)}  $I_k\supset I_{k+1}.$

\end{remark}

\begin{parrafo}\label{IIpar3}
{\em In what follows we define a Rees algebra, say
$\bigoplus_{n\geq 0} I_nW^n$ in $B[W]$, by fixing a set of
generators, say $\mathcal{F}=\{g_{n_i}W^{n_i}/ n_i>0, 1\leq i \leq
m\}$. So if $f \in I_n$, then $f=F_n(g_{n_1},\dots , g_{n_m}),$
where $F_n(Y_1,\dots , Y_m)$ is weighted homogeneous in $m$
variables, and each $Y_j$ has degree $n_j$.}
\end{parrafo}

\begin{parrafo}\label{IIpar2}
{\em Let $B=S[Z]$ be a polynomial ring, and let $Tay: B \to B[U]$
be the $S$-homomorphism defined by setting $Tay(Z)=Z+U$. For
$f(Z)\in B$ set $Tay(f(Z))=\sum_{\alpha \geq 0}
\Delta^{\alpha}(f(Z))U^{\alpha}.$

The operators $\Delta^{\alpha}$ are $S$-differential operators
($S$ linear). Furthermore, for every positive integer $N$,
$\{\Delta^{\alpha}, 0\leq \alpha \leq N\}$ is a basis of the
$B$-module of $S$-differential operators on $B$, of order $\leq
N$.}

\end{parrafo}

\begin{definition}\label{IIdef1}
Set $B=S[Z]$ as before, a polynomial ring over a noetherian ring
$S$. A Rees algebra $\bigoplus I_k\cdot W^k\subset B[W]$ is a
differential Rees algebra, relative to $S$, when:

i) $I_k\supset I_{k+1}$ for $k\geq 0$.

ii) For all $n>0$ and $f\in I_n$, and every index $0\leq j \leq n$
and every $S$-differential operator of order $\leq j$, say $D_j$:
$D_j(f)\in I_{n-j}.$
\end{definition}

\begin{remark}\label{IIrk2} Let $Diff^N_S(B)$ denote the module of
$S$-differential operators of order at most $N$. Then
$Diff^N_S(B)\subset Diff^{N+1}_S(B)\subset \dots$. For this
reason it is natural to require condition (i) in our previous
definition. Note also that \ref{IIpar2} asserts that (ii) can be
reformulated as:

\bigskip

ii') For all $n>0$ and $f\in I_n$, and for every index $0\leq
\alpha \leq n$: $\Delta^{\alpha}(f)\in I_{n-\alpha}.$

In fact, (i)+(ii) is equivalent to (i)+(ii'):
\end{remark}

\begin{theorem}\label{IIth} Fix $B=S[Z]$ as before, and a finite set $\mathcal{F}=\{
g_{n_i}W^{n_i},  n_i>0 , 1\leq i\leq m \}$, with the following
properties:

a) For $1\leq i\leq m$, and every $n'_i$, $0< n'_i\leq n_i$:
$g_{n_i}W^{n'_i}\in \mathcal{F}.$

b) For $1\leq i\leq m$, and for every index $0\leq \alpha < n_i$:
$\Delta^{\alpha}(g_n)W^{n_i-\alpha}\in \mathcal{F}.$

Then the $B$ subalgebra of $B[W]$, generated by $\mathcal{F}$ over
the ring $B$, a differential Rees algebra relative to $S$.
\end{theorem}
\proof

Condition (i) in Def \ref{IIdef1} is by \ref{IIrk1}, 2). Let
$I_NW^N$ be the homogeneous component of degree $N$ of the $B$
subalgebra generated by $\mathcal{F}$. We prove that for all $h\in
I_N$, and $0\leq \alpha \leq N$, $\Delta^{\alpha}(h)\in
I_{N-\alpha}$. The ideal $I_N \subset B$ is generated by all elements of the form
\begin{equation}\label{IIeq1}
H_N=g_{n_{i_1}}\cdot g_{n_{i_2}}\cdots g_{n_{i_p}} \hskip 1cm
n_{i_1}+n_{i_2}+\cdots n_{i_p}=N,
\end{equation}
with the $g_{n_{i_i}}W^{n_{i_i}}\in \mathcal{F}$ not necessarily
different.

Since the operators $\Delta^{\alpha}$ are linear, it suffices to
prove that $\Delta^{\alpha}(a\cdot H_N)\in I_{N-\alpha}$, for
$a\in B$, $H_N$ as in \ref{IIeq1}, and $0\leq \alpha \leq N$. We
proceed in two steps, by proving:

1)  $\Delta^{\alpha}(H_N)\in I_{N-\alpha}$.

2)  $\Delta^{\alpha}(a\cdot H_N)\in I_{N-\alpha}$.

We first prove 1). Set $Tay: B=S[Z] \to B[U]$, as in \ref{IIpar2}.
Consider, for every $g_{n_{i_l}}W^{n_{i_l}}\in \mathcal{F}$,
$$Tay(g_{n_{i_l}})=\sum_{\beta \geq 0}
\Delta^{\beta}(g_{n_{i_l}})U^{\beta}\in B[U].$$

Hypothesis (b) states that for $0\leq \beta < n_{i_l}$,
$\Delta^{\beta}(g_{n_{i_l}})W^{n_{i_l}-\beta}\in \mathcal{F}$. On
the one hand $$Tay(H_N)=\sum_{\alpha \geq 0}
\Delta^{\alpha}(H_N)U^{\alpha},$$ and, on the other hand
$$Tay(H_N)=Tay(g_{n_{i_1}})\cdot Tay(g_{n_{i_2}})\cdots
Tay(g_{n_{i_p}})$$ in $B[U]$. So for a fixed $\alpha$ ($0\leq
\alpha \leq N$), $\Delta^{\alpha}(H_N)$ is a sum of elements of
the form:
$$\Delta^{\beta_1}(g_{n_{i_1}})\cdot
\Delta^{\beta_2}(g_{n_{i_2}})\cdots
\Delta^{\beta_p}(g_{n_{i_p}}),\hskip 0.5cm \sum_{1\leq s\leq p}
\beta_s=\alpha.$$

Therefore it suffices to show that each of these summands is in
$I_{N-\alpha}$. Note here that $$\sum_{1\leq s\leq p}
(n_{i_s}-\beta_s)=N-\alpha,$$ and that some of the integers
$n_{i_s}-\beta_s$ might be zero or negative. Set $$G=\{ r: 1\leq r
\leq p, \mbox{ and } n_{i_r}-\beta_r > 0\}.$$ So
$N-\alpha=\sum_{1\leq s\leq p} (n_{i_s}-\beta_s)\leq \sum_{r\in G}
(n_{i_r}-\beta_r)=M.$

Hypothesis (b) ensures that $\Delta^{\beta_r}(g_{n_{i_r}})\in
I_{n_{i_r}-\beta_r}$ for every index $r\in G$, in
particular:$$\Delta^{\beta_1}(g_{n_{i_1}})\cdot
\Delta^{\beta_2}(g_{n_{i_2}})\cdots
\Delta^{\beta_p}(g_{n_{i_p}})\in I_M.$$ Finally, since $M \geq
N-\alpha,$ $I_M\subset I_{N-\alpha}$, and this proves Case 1).

For Case 2), fix $0\leq \alpha \leq N$. We claim that
$\Delta^{\alpha}(a\cdot H_N)\in I_{N-\alpha}$, for $a\in B$ and
$H_N$ as in \ref{IIeq1}. At the ring $B[U]$, $$Tay(a\cdot
H_N)=\sum_{\alpha \geq 0} \Delta^{\alpha}(a \cdot
H_N)U^{\alpha},$$ and, on the other hand $$Tay(a \cdot
H_N)=Tay(a)\cdot Tay(H_N).$$ This shows that
$\Delta^{\alpha}(a\cdot H_N)$ is a sum of terms of the form $
\Delta^{\alpha_1}(a)\cdot \Delta^{\alpha_2}(H_N)$,$ \alpha_i\geq
0$, and $\alpha_1+\alpha_2=\alpha$. In particular $\alpha_2 \leq
\alpha$; and by Case 1), $\Delta^{\alpha_2}(H_N)\in
I_{N-\alpha_2}$. On the other hand $N-\alpha_2\geq N-\alpha$, so
$\Delta^{\alpha_2}(H_N)\in I_{N-\alpha}$, and hence
$\Delta^{\alpha}(a\cdot H_N)\in I_{N-\alpha}$.
\endproof

\begin{corollary}\label{2corth} A Rees algebra in $B[W]$, generated over
$B$ by $$\mathcal{F}=\{ g_{n_i}W^{n_i},  n_i>0 , 1\leq i\leq m
\},$$ extends to a smallest differential Rees algebra relative to
$S$, which is generated by
$$\mathcal{F'}=\{\Delta^{\alpha}(g_n)W^{n'_i-\alpha};\mbox{ for
all }g_{n_i}W^{n_i}\in \mathcal{F}\mbox{ and } 0\leq \alpha
<n'_i\leq n_i\}.$$
\begin{parrafo}\label{parj}{\em Of particular interest in our
development will be the case $B=S[Z]$ where $S$ is a local regular
ring. In particular both $S$ and $S[Z]$ will be unique
factorization domains. We will consider graded subalgebras in
$B[W]$, and always up to integral closure within this ring.

Assume that $\bigoplus I_k\cdot W^k\subset B[W]$ is a differential
Rees algebra relative to $S$. If, for some positive integer $k$
there is a polynomial , say $f(Z)\in I_k$, which is monic of
degree, say $a<k$, then $\Delta^{a-k}(f(Z))=1$ In this case so $W^{a-k}\in
\bigoplus I_k\cdot W^k$, and the integral closure of this algebra
is all $B[W]$.

Assume now that for some positive integer $b$, there is a monic
polynomial of degree $b$, say $f_b(Z)\in I_b$. Let $f_b^{(e)}$ denote $\Delta^{e}(f_b)$, note that $S[Z][\{f_b^{(e)}(Z)W^{b-e}, 0,\leq e \leq b-1\}]\subset
\bigoplus I_k\cdot W^k$.

If, in addition, there is a factorization of $f_b(Z)$, of the
form: $f_b(Z)=f_{c_1}(Z)\cdot f_{c_2}(Z)\cdots f_{c_r}(Z)$, where
each factor $f_{c_i}(Z)$ is a monic polynomial of degree $c_i$,
then \ref{eq660} asserts that
$$S[Z][\{f_b^{(e)}(Z)W^{b-e}, 0,\leq e \leq b-1\}] \subset S[Z][\{f_{c_i}^{(\beta)}(Z)W^{c_i-\beta}\},1\leq i\leq r, 0\leq
\beta \leq c_i ]$$ is finite. In particular,
each element $f_{c_i}(Z)W^{c_i}$ is integral over $\bigoplus
I_k\cdot W^k$.

 }
\end{parrafo}
\end{corollary}

%%%%%%%%%%%%%%%%%%%%%%%%%%%%%%%%%%%%%%%%%%%%%%%%%%%%%%%%%%%%%%%%%%%%%%%%%%%%%%%%%%%
%%%%%%%%%%%%%%%%%%%%%%%%%%%%%%%%%%%%%%%%%%%%%%%%%%%%%%%%%%%%%%%%%%%%%%%%%%%%%%%%%%%%%

%%%%%%%%%%%%%%%%%%%%%%%%%%%%%%%%%%%%%%%%%%%%%%%%%%%%%%%%%%%%%%%%%%%%%%%%%%%%%%%%%%%%%%%%%%

%%% ----------------------------------------------------------------------

\section{ Differential Rees algebras on smooth schemes}

Since the outstanding Theorem of Hironaka was published in \cite{Hir64} ,
an important effort was done to simplify the proof of this Theorem. One of the major  steps in this direction was achieved by Jean Giraud in 
\cite{Gir} and \cite{Giraud1975}. This simplification grows from the new technics, which were introduced there, involving differential operators.

In Section 3 of \cite{VV1} it is proved that every Rees algebra can be naturally extended, in a unique manner, to a new Rees algebra enriched by the action of differential operators. These are called 
Differential Rees algebras, which are to be discussed in this section. 
For the ease of the exposition we also recall some properties of these extensions, introduced in \cite{VV1}, as they are to be used in the coming sections.

\begin{parrafo}\label{par3.0}

{\rm Let $V$ be a smooth scheme over a field. A sequence of coherent ideals on $V$, say
$\{I_n\}_{n\in \mathbb{N}}$, such that $I_0=\mathcal{O}_V$, and
$I_k\cdot I_s\subset I_{k+s}$, defines a graded sheaf of algebras
$\bigoplus_{n\geq 0} I_n\cdot W^n \subset  \mathcal{O}_{V}[W]$.

We say that this algebra is a Rees algebra over $V$, if there is an open
covering of $V$ by affine open sets $\{U_i\}$, so that
$\bigoplus_n I_n(U_i) W^n \subset \calo_V(U_i)[W]$ is a finitely
generated $\mathcal{O}_V(U_i)$-algebra. In what follows $V$ will
denote a smooth scheme of a perfect field $k$, and $Diff^r_k(V)$, or
simply $Diff^r_k$, denotes the locally free sheaf of $k$-linear
differential operators of order at most $r$.}

\end{parrafo}

\begin{definition} \label{3def1} {\bf 1)} A Rees algebra $\mathcal{G} $ defined by $\{I_n\}_{n\in
\mathbb{N}}$ is a {\em differential Rees algebra}, relative to the
field $k$, or simply a {\em Diff-algebra}, if:

i) $I_n\supset I_{n+1}$.

ii) There is open covering of $V$ by affine open sets $\{U_i\}$,
and for every $D\in Diff^{(r)}(U_i)$, and $h\in I_n(U_i)$,
 $D(h)\in I_{n-r}(U_i)$, provided that $n\geq r$.
 
 {\bf 2)} We sometimes fix a smooth morphism $\beta: V\to V'$, of smooth schemes over $k$. In this case $\mathcal{G} $ is said to be a 
 relative differential algebra, or a $\beta$-differential algebra, when
 ii) holds only for those $D\in Diff_{\beta}^{(r)}(U_i)$, namely for $\beta$-linear differential operators.

\end{definition}
The abbreviated notation of  "{\em Diff-algebra}" omits reference
to the underlying field $k$, but this will be clear from the
context. Due to the local nature of the definition, we reformulate
this notion in terms of $k$-algebras.

\begin{definition}
In what follows $R$ will denote a smooth algebra over a field, or
a localization of such algebra at a closed point (a regular local
ring). A Rees algebra is defined by a sequences of ideals
$\{I_k\}_{k\in \mathbb{N}}$ such that:

1) $I_0=R$, and $I_k\cdot I_s \subset I_{k+s}$.

2) $\bigoplus I_k W^k$ is a finitely generated $R$-algebra.

We say that the Rees algebra is a
differential Rees algebra relative to $k$, or simply a
Diff-algebra, if, in addition to the previous conditions:

3) $I_n\supset I_{n+1}$. 

4) given $D\in Diff_k^{(r)}(R)$, $D(I_n)\subset I_{n-r}$.

\end{definition}

\begin{theorem}\label{3th1}
Fix a smooth scheme $V$ over a perfect field $k$. Assume that $\mathcal G=\bigoplus I_k\cdot W^k$ is a Rees algebra
over $V$. Then there is a smallest extension of it
to a differential Rees algebra relative to the field $k$.
\end{theorem}
The Theorem says that a Rees algebra $\mathcal G$, over a smooth scheme $V$
extends to a smallest Diff-algebra (i.e., included in every other
Diff-algebra containing it). The latter is called the Diff-algebra spanned 
by $\mathcal G$. We refer here to Th 3.4 in \cite{VV1}
for the proof. Let us indicate that it follows easily from the argument for the
one-variable case (Th \ref{IIth}). Fix
a closed point $x\in V$, and a regular system of parameters
$\{x_1,\dots,x_n\}$ for $\calo_{V,x}$, then smoothness of $V$
locally at $x$ asserts that there is a ring homomorphism, say:
$$Tay: \calo_{V,x}\to \calo_{V,x}[[T_1,\dots ,T_n]], \hskip 1cm
Tay(f)=\sum_{\alpha\in (\ent)^n}\Delta^{\alpha}(f)T^{\alpha}$$
where $Tay(x_i)=x_i+T_i$. Furthermore, as the underling field $k$ is perfect, $\{ \Delta^{\alpha}, \alpha
\in (\nat)^n , 0 \leq |\alpha| \leq c\}$ generate the
$\calo_{V,x}$-module $Diff^c_k(\calo_{V,x})$.

The proof of the previous Theorem shows that, at a suitable affine
neighborhood of $x$ in $V$, where say  $\bigoplus I_k\cdot W^k$ is
generated by $\mathcal{F}=\{ g_{n_i}W^{n_i}, n_i>0 , 1\leq i\leq m
\},$ then
\begin{equation}\label{eq3422} \small
\mathcal{F'}=\{\Delta^{\alpha}(g_n)W^{n'_i-\alpha}/
g_{n_i}W^{n_i}\in \mathcal{F}, \alpha=(\alpha_1, \alpha_2,\dots ,
\alpha_n) \in (\mathbb{N})^n, \mbox{ and } 0\leq |\alpha| < n'_i
\leq n_i\}
\end{equation}
generates the smallest extension of $\bigoplus I_k\cdot W^k$ to a
Diff-algebra (relative to $k$). In particular, a Rees ring of an
ideal (all $n_i=1$ in $\mathcal{F}$) is a Diff-algebra.

\begin{corollary}\label{coragre} Given inclusions of
Rees algebras, say $$\mathcal{G}=\bigoplus I_n W^n \subset
\mathcal{G'}=\bigoplus I'_n W^n \subset \mathcal{G}''=\bigoplus
I''_n W^n, $$ where $\mathcal{G}''$ is the Diff-algebra spanned by
$\mathcal{G}$, then $\mathcal{G}'' $ is also the Diff-algebra
spanned by $\mathcal{G}'$.

\end{corollary}

\begin{center}

{\bf Differential Rees algebras and singular locus.}

\end{center}

\begin{parrafo}

{\rm The notion Diff-algebras, over a smooth $k$-scheme $V$, is
closely related to the notion of {\em order} at the local regular
rings of $V$ when $k$ is a perfect field. Recall that the order of a non-zero ideal $I$ at a
local regular ring $(R,M)$ is the highest integer $b$ for which
$I\subset M^b$.

If $I\subset \calo_V$ is a sheaf of ideals, $V(Diff^{b-1}_{k}(I))$
is the closed set of points of $V$ where the ideal has order at
least $b$. We analyze this fact locally at a closed point $x$.

Let $\{x_1,\dots ,x_n\}$ be a regular system of parameters for
$\calo_{V,x}$, and consider the differential operators
$\Delta^{\alpha}$, defined on $\calo_{V,x}$ in terms of these
parameters, as in the Theorem \ref{3th1}. So at $x$,
$$(Diff^{b-1}_{k}(I))_x=\langle \Delta^{\alpha}(f)/ f\in I, 0\leq
| \alpha| \leq b-1\rangle.$$ One can now check at $\calo_{V,x}$,
or at the ring of formal power series $\hat{\calo}_{V,x}$, that
$Diff^{b-1}_{k}(I)$ is a proper ideal if and only if $I$ has order
at least $b$ at the local ring.

The operators $\Delta^{\alpha}$ are defined globally at a suitable
neighborhood $U$ of $x$. So if $\bigoplus I_n\cdot W^n \subset
\calo_{V}[W]$ is a Diff-algebra, and $x\in V$ is a closed point, the
Diff-algebra $\bigoplus (I_n)_x\cdot W^n \subset \calo_{V,x}[W]$ will
be properly included in $\mathcal{O}_{V,x}[W]$, if and only, for
each index $k\in \mathbb{N}$, the ideal $(I_k)_x$ has order at
least $k$ at the local regular ring $\mathcal{O}_{V,x}$.}

\end{parrafo}

\begin{definition}\label{3defsing} The {\em singular locus} of a Rees algebra $\mathcal{G}=\bigoplus I_n\cdot W^n
\subset \calo_{V}[W]$, will be  $$Sing(\mathcal{G})=\cap_{r\geq 0}
V(Diff^{r-1}_{k}(I_r)) (\subset V).$$ It is the set of points
$x\in V$ for which all $(I_r)_x$ have order at least $r$ (at
$\mathcal{O}_{V,x}$).
\end{definition}

\begin{remark}\label{3rkdelsing}

 Assume that $f\in (I_r)_x$ has order $r$ at $\mathcal{O}_{V,x}$.
Then, locally at $x$,  $Sing(\mathcal{G})$ is included in the set
of points of multiplicity $r$ (or say, $r$-fold points) of the
hypersurface $V(\langle f \rangle)$.

In fact $Diff^{r-1}_{k}(f)\subset Diff^{r-1}_{k}(I_r)$, and the
closed set defined by the first ideal is that of points of
multiplicity $r$.
\end{remark}
\begin{proposition}\label{3propsing}
\begin{enumerate}

\item[(1)]  If $\mathcal{G}$  is a Rees algebra generated over $\calo_V$ by
$\mathcal{F}=\{ g_{n_i}W^{n_i}, n_i>0 , 1\leq i\leq m \}$, then
$Sing(\mathcal{G})=\cap V( Diff^{n_i}(\langle g_i \rangle ) ).$

\item[(2)]  If $\mathcal{G}=\bigoplus I_n\cdot W^n$ and $ \mathcal{G}'=\bigoplus I'_n\cdot W^n$ are Rees algebras with
the same integral closure (e.g. if $\mathcal{G}\subset
\mathcal{G}'$ is a finite extension), then
$Sing(\mathcal{G})=Sing(\mathcal{G}').$

\item[(3)]  Let $\mathcal{G}''=\bigoplus I''_n\cdot W^n$ be the extension
of $\mathcal{G}$ to a Diff-algebra, as defined in Theorem
\ref{3th1}, then $Sing(\mathcal{G})=Sing(\mathcal{G}'')$.

\item[(4)]  For every Diff-algebra $\mathcal{G}''=\bigoplus I''_n\cdot
W^n$, $Sing(\mathcal{G}'')=V(I''_1)$.

\item[(5)]  Let  $\mathcal{G}''=\bigoplus I''_n\cdot W^n$ be a
Diff-algebra. For every index $r$, $Sing(\mathcal{G}'')=V(I''_r)$.
\end{enumerate}
\end{proposition}
\proof

1) We have formulated 1) with a global condition on $V$, however
this is always the case locally. In fact, there is a covering of
$V$ by affine open sets, so that the restriction of $\mathcal{G}$
is generated by finitely many elements. Let $U$ be such open set,
so $\mathcal{G}(U)=\bigoplus I_k(U)\cdot W^k$ is generated by
$\mathcal{F}=\{ g_{n_i}W^{n_i}, n_i>0 , 1\leq i\leq m \}$,
$g_{n_i}\in \mathcal{O}(U)$. The claim is that
 $y\in Sing(\mathcal{G})\cap U$ if and only if
the order of $g_{n_i}$ at $\calo_{V,y}$ is at least $n_i$, for
$1\leq i\leq m$.

The condition is clearly necessary. Conversely, if
$\mathcal{G}=\bigoplus I_n =\calo_U[\{g_iW^{n_i}\}_{g_iW^{n_i}\in
\mathcal{F}}]$, and each $g_{n_i}$ has order at least $n_i$ at
$\calo_{V,y}$, then $I_n$ (generated by weighted homogeneous
expressions on the $g_i$'s) has order at least $n$ at
$\calo_{V,y}$.

2) Rees algebras are, locally, finitely generated over $\calo_{V}$. This ensures that there is a (infinite) semigroup $\mathbb{M}$ in $\mathbb{N}$,
so that for every $n\in \mathbb{M}$ both $I_n$ and $I'_n$ have the
same integral closure. We may also define both singular loci as
the points where these ideals have order at least $n$. The
equality in 2) holds because the order of an ideal in a local
regular ring, is the same as the order of the integral closure
(\cite{ZS} Appendix 3).

3) We argue as in 1), here we may also assume that there is $x\in
U$, a regular system of parameters $\{x_1,\dots ,x_n\}$ at $x$,
and differential operators $\Delta^{\alpha}$ as in the Theorem
\ref{3th1}, defined globally at $U$. The Diff-algebra
$\mathcal{G}''$ in Theorem \ref{3th1}, is defined by
$$\mathcal{F'}=\{\Delta^{\alpha}(g_n)W^{n_i-\alpha}/
g_{n_i}W^{n_i}\in \mathcal{F}, \alpha=(\alpha_1, \alpha_2,\dots ,
\alpha_n) \in (\mathbb{N})^n, \mbox{ and } 0\leq |\alpha| <
n_i\}.$$

Note finally that if the order of $g_{n_i}$ at a local ring is
$\geq n_i$, then the order of $\Delta^{\alpha}(g_n)$ is $\geq
n_i-|\alpha|$.

4) The inclusion $Sing(\mathcal{G}'')\subset V(I''_1)$ holds, by
definition, for every Rees algebra. On the other hand, the
hypothesis ensures that $Diff^{r-1}(I''_r)\subset I''_1$, so
$Sing(\mathcal{G}'')\supset V(I''_1)$.

5) Follows from 4).
\endproof

%%%%%%%%%%%%%%%%%%%%%%%%%%%%%%%%%%%%%%%%%%%%%%%%%%%%%%%%%%%%%%%%%%%%%%%%%%%%%%%%%%%%%%%%%%%%%%%%%%%%%%%%%%%%%%%%

%%%%%%%%%%%%%%%%%%%%%%%%%%%%%%%%%%%%%%%%%%%%%%%%%%%%%%%%%%%%%%%%%%%%%%%%%%%%%%%%%%%%%%%%%%%%%%%%%%%%%%%%%%%%%%%%%%%

%%%%%%%%%%%%%%%%%%%%%%%%%%%%%%%%%%%%%%%%%%%%%%%%%%%%%%%%%%%%%%%%%%%%%%%%%%%%%%%%%%%%%%%%%%%%%%%%%%%%%%%%%%%%%%%%%%%%%%

%%%%%%%%%%%%%%%%%%%%%%%%%%%%%%%%%%%%%%%%%%%%%%%%%%%%%%%%%%%%%%%%%

\section { Simple differential algebras and
projections}
 
This work, as a whole, is motivated by the problem of resolution of singularities over arbitrary fields. Generally speaking, the open problem of resolution of singularities over a perfect field reduces to that of constructing {\em resolutions} of Rees algebras. A second fundamental reduction of the open problem of resolution is that one can take, as starting point, a Rees algebra 
which is, in addition, a differential Rees algebra. 

We do not define here resolution of Rees algebras, neither do we address these two fundamental reductions in this work, as the subject has been thoroughly studied. But let us indicate that this already justifies our particular attention on {\em differential} Rees algebras.
%due to its relevance in this open problem.

 So, summarizing, the problem of embedded resolution over perfect fields reduces 
to that of resolution of differential Rees algebras. This is what we know, but what we do not know, at least at the moment, is to define resolutions of this kind of Rees algebras. 

The hope is to prove the latter by using some form of induction, a strategy which works in characteristic zero. The form of induction we are searching for leads us firstly 
to the discussion of the $\tau$-invariant in \ref{partau}. In fact, the
$\tau$-invariant is a positive integer, and the larger this integer is, the simpler it is to construct a resolution.

This Section aims to the presentation of an entirely new form of induction.
 It is sustained, essentially, on elimination theory, and makes use of the algebras that were introduced in Definition \ref{defmathcalR}. The properties of these are gathered in Theorem \ref{4th10}, which is the main result in this Section. 

%Elimination algebras are defined here in terms of universal graded rings, as discussed in \ref{ontges}. Elimination is defined here with a natural compatibility  with finite extensions of Rees algebras.

\vskip 1cm

Hironaka's Theorem of resolution of singularities in characteristic zero 
uses a form of induction which is based on a reformulation of the problem but now in a smooth hypersurface. So his induction is on the dimension on the ambient space. The form of induction we present here is different from that of Hironaka, but leads to the same result over fields of characteristic zero. Here the restriction to a smooth hypersurface is replaced by the elimination of one variable.
Both approaches  make use of Hironaka's notion of $\tau$-invariants.
This is a powerful invariant attached to a point, as it indicates the total number of variables that can be eliminated, at least in a neighborhood of the point (see also \ref{par511}).

Fix a Rees algebra over a smooth scheme $V$, say
$\mathcal{G}=\bigoplus \cali_k\cdot W^k (\subset \calo_V[W])$, and
a closed point $x\in Sing(\mathcal{G})$. Let
$\mathcal{G}_x=\bigoplus I_k\cdot W^k (\subset R[W])$ be the
localization at $R=\calo_{V,x}$.
 Observe that the Diff-algebra spanned by $\mathcal{G}$ induces,
at $x\in V$, the Diff-algebra spanned by $\mathcal{G}_x$.
\begin{definition}\label{defsmpt}
 A Rees algebra $\mathcal{G}=\bigoplus \cali_k\cdot W^k (\subset \calo_V[W])$ is said to be simple at a point $x\in Sing(\mathcal{G})$,  if, for some $n$, the order
of $I_n$ is $n$ at the local ring $ \mathcal{O}_{V, x}$.
\end{definition}

\begin{parrafo}\label{partau}{\em Here $R=\mathcal{O}_{V,x}$ is a local regular ring. The
graded algebra of the maximal ideal, say $gr_M(R)$ is a polynomial
ring. Recall that $Spec( gr_M(R))$ is the tangent space of $V$ at $x$. Attach
to $\mathcal{G}$ an homogeneous ideal in $gr_M(R)$, called the
initial (or tangent) ideal of $\mathcal{G}_x$, spanned by
$In_k(I_k)$, for all index k. It defines a closed set in $Spec(
gr_M(R))$  called the {\em tangent cone} of $\mathcal{G}_x$.

The tangent ideal of a Rees algebra at a closed point $x\in
Sing(\mathcal{G})$ is zero unless $x$ is a simple point.
In this case $gr_M(R)=k'[Z_1,\dots, Z_n]$ (polynomial ring in
$n$ variables), where $k'$ denotes the residue field at $x$. In the case in which $\mathcal{G}$ is a Diff-algebra
it is easy to check that the tangent ideal is closed when applying
homogeneous differential operators of the form $\Delta^{\alpha}$,
for every multi-index $\alpha=(\alpha_1,\dots, \alpha_n)$, where these operators are defined
by taking Taylor expansions in terms of the variables $Z_i$. In
other words, if $H$ is an homogeneous element of degree $N$ in the
tangent ideal, then $\Delta^{\alpha}(H)$ is homogeneous of degree
$N-|\alpha|$, and belongs to the tangent ideal.

Homogeneous ideals with this property were studied in
(\cite{Hironaka70}). If $k'$ is a field of
characteristic zero, then the ideals with this property are
exactly  those generated by linear forms, which we may take to be
$ Z_1,..,Z_{\tau}$. If $k'$ is a field of characteristic $p$, the
ideals with this property are generated by elements of the form
\begin{equation}\label{eqptau}
l_1,\dots. l_{s_0}, l_{s_0+1}, \dots ,l_{s_1}, \dots,
l_{s_{r-1}}\dots l_{s_{r}},
\end{equation} where each $l_j$ is a linear
combination of powers $Z_j^{p^t}$, for $ s_{t}\leq j \leq
s_{t+1}$; and no $l_j$ is in the ideal spanned by the previous
elements. When $k'$ is a perfect field, one can take $l_j=Z_j^{p^t}$, for $ s_{t}\leq j \leq
s_{t+1}$

It is said that the initial ideal (defining the tangent cone at the point) is spanned by a flag of
Frobenius-linear ideals in powers of the characteristic. The tangent cone is a subspace of the tangent space in characteristic zero, and a additive subgroup in arbitrary characteristic.

There are
two important invariants defined by the tangent cone at the point:

1) the integer $s_{r}$, usually called the invariant $\tau$ of the
singularity (at the point $x$), and

2) the smallest integer $e_0$ so that  $p^{e_0}$ is the smallest
power which arises in the description of the elements $l_i$ in the
previous flag. When $\mathcal{G}$ is a Diff-algebra the order (at $R$) of $\cali_n$ is $n$ if and only if
$n$ is a multiple of $p^{e_0}$.

We refer to the work of T. Oda ( \cite{Oda1973}, \cite{Oda1983},
and \cite{Oda1987}) where these notions are studied. 
Related to this notion, but different and original, is that introduced by Kawanoue in \cite{kaw}, consisting on graded sub-algebras of tangent algebra $gr_M(R)$ (as opposed the previous discussion, based on homogeneous ideals in $gr_M(R)$).
}

\end{parrafo}

\begin{parrafo}\label{p66au}{ Algebras generated by monic polynomials.} 

{\em Take a Rees algebra and a closed point $x$ as above, and let $q$ be an integer such that, for all natural number $m$,
$I_{mq}$ has order $mq$ at the local regular ring $R=\mathcal{O}_{V,x}$. Take now the completion, so assume that
$R$ is complete. 

Once we fix a regular system of parameters for $R$, say $\{ z,
x_1,\dots x_{d-1}\}$, then one can set $gr_M(R)=k'[ Z,
X_1,\dots X_{d-1}]$ where $Z$ is the initial form of $z$, and each 
$X_i$ is the initial form of $x_i$. Moreover, after enlarging the base field one can take the regular system of parameters 
so that the line $V(<X_1,\dots X_{d-1}>)$ is
transversal to the tangent cone in $Spec(gr_M(R))$. The Weierstrass Preparation Theorem  ensures that, for some index $n=mq$,
there is an element $f\in I_n$ of order $n$, that multiplied by a
suitable unit of $R$, can be expressed as a monic polynomial of
degree $n$ in $S[z]$. Here $S$ is a regular local ring with coordinates
$\{ x_1,\dots x_{d-1}\}$, and say $$ F(z)=z^n+a_1 z^{n-1}+\dots
+a_n.$$

Observe that in this case $$I_{(n)}:= I_n \cap S[z]$$ is
an ideal spanned by monic polynomials of degree $n$. To check
this, set $$ A=S[z]/\langle F \rangle (=R/\langle f \rangle),$$
and note that each $g\in I_n$ has a class, say:
$b_1z^{n-1}+b_2z^{n-2}+\dots +b_n$. On the other hand
$$G(z)=(z^n+a_1 z^{n-1}+\dots +a_n)+ (b_1z^{n-1}+b_2z^{n-2}+\dots
+b_n)$$ is a monic polynomial in $I_n \cap S[z]$, and  all
monic polynomials arising in this manner span $I_{(n)}$. In this
case $I_n=I_{(n)}R$, thus $I_n$ is generated by monic polynomials of
degree $n$ in $S[z]$.

Note here that after an enlargement of the base field, $\bigoplus I_k\cdot W^k (\subset R[W])$ can be
generated by elements in $S[z][W]$. In fact, let $\mathcal{F}=\{
g_{n_i}W^{n_i}, n_i>0 , 1\leq i\leq m \} \subset R[W]$ be a set of
generators of this graded $R$ subalgebra of $R[W]$. We can always 
take the system of coordinates so that, up to multiplication by a unit in $R$, each $g_{n_i}$ is a monic polynomial in $z$, of some degree, say $m_i$.
%%VVVVVVVVVVVVVVV
%We may always
%choose the integer $n$, above, to be bigger that all $n_i$. As $I_n\subset
%I_{n_i}$, then $F(Z)\in I_{n_i}$ and each $g_{n_i}W^{n_i}$ can be
%replaced by $G_{n_i}W^{n_i}$, where $G_{n_i}\in S[Z]$ denotes the
%class of $g_{n_i}$ in $A$.

A remarkable fact about simple Rees algebras at a 
local regular ring $R=\mathcal{O}_{V,x}$ of a point $x\in Sing(\mathcal G)$, is that, up to integral closure, they can be
generated by {\em monic} polynomials in $S[z]$, say
$$\mathcal{F}=\{ G_{n_i}W^{n_i}, n_i>0 , 1\leq i\leq m \} \subset
S[Z][W],$$ where now each $G_{n_i}$ is monic of degree $n_i$ in $S[z]$.  In fact, if we choose $n$ to be divisible by all $n_i$, it is
clear that $\bigoplus I_k\cdot W^k$ is a finite extension of
$G'=R[I_nW^n](\subset R[W^n])$; and, as it was indicated above, $I_n$ can
be generated by monic polynomials of degree $n$ in $S[Z]$.

Since $\bigoplus I_k\cdot W^k$ is a noetherian subalgebra of
$R[W]$, we may assume that so is $$\bigoplus I_{(k)}\cdot W^k
(\subset S[Z][W]).$$

This observation, together with the discussion in  \ref{1def29} ,
where elimination algebras where defined, will lead us to a form
of elimination for Diff-algebras, locally at simple points.

}
\end{parrafo}

%%%%%%%%%%%%%%%%%%%%%%%%%%%%%%%%%%%%%%%%%%%%%%%%%%%%%%%%%%%%%%%%%%%%%%%%%%%%%%%
%%%%%%%%%%%%%%%%%%%%%%%%%%%%%%%%%%%%%%%%%%%%%%%%%%%%%%%%%%%%%%%%%%%%%%%%%%%%%%%%%%

\begin{center} {\bf On Rees algebras and integral extensions.}

\end{center}

We discuss here the concept of finite extensions of Rees algebras,
as it will be used later in the definition of elimination (or
projection of Diff-algebras). Fix a noetherian ring $B$ and ideals
defining a Rees algebra
 $\bigoplus_{k\geq 0} I_k\cdot W^k \subset (\subset B[W])$  as in \ref{IIpar1}. Set
$$\bigoplus_{k\geq 0} I_k\cdot W^k =B[\{I_nW^n, n\geq 0\}].$$

An inclusion of Rees algebras $B[\{I_nW^n, n\geq 0\}]\subset
C[\{J_nW^n, n\geq 0\}]$, is defined by a ring extension $B\subset
C$, and an inclusion of ideals $I_n\subset J_n$ for each $n$. They
arise in various ways:

Given a Rees algebra $B[\{I_nW^n, n\geq 0\}]$ and a positive
integer $m$ define
\begin{equation}\label{4eq1} V^{(m)}(B[\{I_nW^n, n\geq 0\}])=\bigoplus_{n\geq0} I_{mn}W^{mn}
(\subset B[\{I_nW^n, n\geq 0\}]).
\end{equation}
An inclusion, which is an integral extension, is that of:
\begin{equation}\label{4eq1a} V^{(m)}(B[\{I_nW^n, n\geq 0\}])\subset B[\{I_nW^n, n\geq
0\}].
\end{equation}

Let $A\subset B$ be a ring extension, and $B[\{I_nW^n, n\geq 0\}]
$ a Rees algebra. An inclusion of  Rees algebras arises by
setting
\begin{equation}\label{4eq2}
B[\{I_nW^n, n\geq 0\}]\cap A[W]\subset B[\{I_nW^n, n\geq 0\}],
\end{equation}  where the left hand side is a graded subring of $A[W]$.
Finally, given rings $B\subset B'$, then:
\begin{equation}\label{4eq3} B[\{I_nW^n, n\geq 0\}]\subset B'[\{I'_nW^n, n\geq 0\}],
\end{equation}
 where $I'_n=I_nB'$, also defines a graded extension.

\begin{parrafo} {\em Let $B[\{I_nW^n, n\geq 0\}] $ be a Rees algebra, and assume that $A\subset B$ is a finite extension of rings. In
this case one could expect that $B[\{I_nW^n, n\geq 0\}] $ be a
finite extension of the intersection algebra (\ref{4eq2}). Example
\ref{ejdqno} shows that this is not so in general. However this
will be the case for our notion of elimination for Diff-algebras
(see Theorem \ref{4th10}).

}
\end{parrafo}

\begin{remark}\label{rk4i5}1) The extension (\ref{4eq1}) is integral, and so
is (\ref{4eq3}) when $B\subset B'$ is integral.

2) The extension $B[\{I_nW^n, n\geq 0\}]\subset C[\{J_nW^n, n\geq
0\}]$ is integral if and only if
$$V^{(m)}(B[\{I_nW^n, n\geq 0\}])\subset V^{(m)}(C[\{J_nW^n, n\geq 0\}])$$ is integral for
some $m$.

\end{remark} \proof 1) is clear; 2) follows from the finiteness in
(\ref{4eq1a}).

\begin{example}\label{ejdqno}
Set $A=k[x_1,\dots,x_n]_M\subset B=A[Z]/\langle f(Z) \rangle$,
where $M=\langle x_1,\dots ,x_n\rangle$, and $f(Z)$ is a monic
polynomial of degree $e$, and $f(Z)\in \langle M,Z \rangle$. Let
$M$ and $N$ denote the maximal ideals of $A$ and $B$, and assume
that $In(f(Z))\in gr_{\langle M,Z
\rangle}(A[Z])=k[X_1,\dots,X_n,Z]$ is such that $\{ In(f),
X_1,\dots, X_n\}$ is a regular sequence.

In this case $gr_M(S)\to gr_N(B) (= gr_{\langle M,Z
\rangle}/\langle In(f) \rangle)$ is flat. Note that:

i) $\bigoplus N^k\cdot W^k \cap A[W] = \bigoplus M^k\cdot W^k .$

ii) The ring extension (\ref{4eq2}) is not finite in this case,
unless $f(Z)\in \langle M,Z \rangle^e$.

To prove i), use the fact that flatness ensures an inclusion
$gr_M(S)\subset gr_N(B)$.

\end{example}

\begin{center}

{\bf On differential Rees algebras and projections.}

\end{center}
\begin{parrafo}\label{4rksindelp}

{\rm Fix a Diff-algebra $\mathcal{G}$ over a smooth scheme $V$, and
a closed point $x\in Sing(\mathcal{G})$. Define
$\mathcal{G}_x=\bigoplus I_k\cdot W^k (\subset
\mathcal{O}_{V,x}[W])$ by localization, which, in a natural sense,
is also a Diff-algebra, since differential operators also act on
$\mathcal{O}_{V,x}$. Assume that $x\in Sing(\mathcal{G})$ is a simple point,
and hence that $\mathcal{G}_x$ is a simple Diff-algebra.

We make use of the following handy observations:

\bigskip
{\em Observation 1):} Fixed a point $x$ in a smooth scheme $V$ of dimension $d$. It is simple to define smooth schemes $V'$, of dimension $d-1$, together with a smooth morphisms, say
 $\pi: V\to V'$. 
 
 The construction of such morphisms require restrictions to  \'etale neighborhoods of $x$.

\bigskip
{\em Observation 2):} If $\mathcal{G}$ is a Diff-algebra over $V$, and if $\pi: V\to V'$ is a smooth morphism over $k$, then $\mathcal{G}$ is also a differential algebra relative to this morphism (see Def \ref{3def1}).

\bigskip
As for Observation 1) recall that $(V,x)$ is an \'etale neighborhood
of, say $({\mathbb A}^d,x')$, where $d$ is the dimension of $V$ at $x$.
\'Etale maps are of course smooth, so the claim reduces to the affine case. The underlying field $k$ is assumed to be perfect. So after suitable separable extension one can assume that both $x\in V$ and $x'\in {\mathbb A}^d$ are rational points over $k$. In this case the assertion is simple, as one can we can even take linear projections on smaller dimensional affine schemes.

 Set $S\subset \mathcal{O}_{V',\pi(x)}$, so $S\subset \mathcal{O}_{V,x}$ is the
inclusion of regular rings defined as
in \ref{p66au}, by a transversal projection $\pi$. Set $R=S[Z]_{<M_S,Z>}(\subset \mathcal{O}_{V,x})$.
Choose elements $f_{c_i}(Z)\in I_{({c_i})}$, monic of degree $c_i$
in $S[Z]$, $1\leq i \leq r$, and set
$\mathcal{F}=\{f_{c_i}(Z)W^{c_i}; 1\leq i \leq r\}\subset
S[Z][W]$. So $S[Z]\subset R$, and 
\begin{equation}\label{471eq}
S[Z][\{f_{c_i}(Z)W^{c_i}; 1\leq i \leq
r\}]\subset \mathcal{G}_x.
\end{equation}
Moreover, the discussion in  \ref{p66au} ensures that such inclusion can be taken to be a finite extension.
%XXXXXXXXXXXXXXBORRAR
%We make use of the following handy observations:

%\bigskip
%{\em Observation 1):} If $\mathcal{G}$ is a Diff-algebra over $Z$, and if $Z\to Z'$ is a smooth morphism over $k$, then $\mathcal{G}$ is also a Diff-algebra relative to this morphism. 

%\bigskip
%{\em Observation 2):} Fixed a point $x$ in a smooth scheme $Z$, it is simple to define smooth schemes $Z'$ together with smooth morphisms $Z\to Z'$ at an \'etale neighborhood of $x$.

%\bigskip
%As for Observation 2) recall that $(Z,x)$ is an \'etale neighborhood
%of, say $({\mathbb A}^d,x')$, where $d$ is the dimension of $Z$ at $x$.
%\'Etale maps are of course smooth, so the claim reduces to the affine case.Usually the underlying field $k$ is assumed to be perfect. So after suitable separable extension one can assume that both $x\in Z$ and $x'\in {\mathbb A}^d$ are rational points over $k$. In this case the assertion is simple, as one can we can even take linear projections on smaller dimensional affine schemes.
%XXXXXXXXXXXXXXXXXXXXXXXXXXXXX

Observation 2) ensures that, as $\mathcal{G}_x$
is a differential Rees algebra relative to the structure field
$k$, it is also closed by differentials relative to $S$; and hence
$$S[Z][\{\Delta^{\alpha}(f_{c_i})W^{n_i-\alpha}/ f_{c_i}W^{n_i}\in
\mathcal{F}, \mbox{ and } 0\leq \alpha < c_i\}]\subset
\mathcal{G}_x \  (\ref{2corth}).$$

 Recall that the elimination algebra
$\overline{\mathcal{R}}_{f_{c_1},f_{c_2},\dots,f_{c_r}}$
(\ref{1eqdrbar}) is defined as the graded $S$ subalgebra of $
S[W]$, generated by $\mathcal{F}=\{G_{c_1,\dots,
c_r}^{l}(a^{(1)}_1,\dots,a^{(1)}_{c_1}, \dots
,a^{(r)}_1,\dots,a^{(r)}_{c_r}))W^{deg \ G_{c_1,\dots, c_r}^{l}}
\}$ (notation as in Th \ref{vramloc}), where $G_{c_1,\dots,
c_r}^{l}$ runs among the generators of the algebra
$\overline{R}_{c_1,\dots, c_r}$ (\ref{eqdeR}), and $deg \
G_{c_1,\dots, c_r}^{l}$ is the degree of the weighted homogeneous
polynomial $ G_{c_1,\dots, c_r}^{l}$.

Here $S$ is a regular local ring, and
$Sing(\overline{\mathcal{R}}_{f_{c_1},f_{c_2},\dots,f_{c_r}})$ is
the closed set in $Spec(S)$ defined by the points where each
element $G_{c_1,\dots, c_r}^{l}(a^{(1)}_1,\dots,a^{(1)}_{c_1},
\dots ,a^{(r)}_1,\dots,a^{(r)}_{c_r})) $ has order at least $deg \
G_{c_1,\dots, c_r}^{l}$ (Prop \ref{3propsing}, 1)).}

\end{parrafo}

\begin{lemma}\label{3lemuno}
Let $\mathcal{G}_x=\bigoplus I_k\cdot W^k (\subset
\mathcal{O}_{V,x}[W])$ (localization at $x\in Sing(\mathcal{G})$),
be a simple Rees algebra at the point, and set $R=S[Z]_{<M_S,Z>}(\subset
\mathcal{O}_{V,x})$ for a suitable transversal projection. Assume that $\mathcal{G}_x$ is a differential algebra relative to $S$ (e.g., that $\mathcal{G}_x$ is a Diff-algebra)(see \ref{IIdef1}).

 Choose elements $f_{c_i}\in I_{c_i}$
which are monic polynomials of degree $c_i$ in $S[Z]$ and order
$c_i$ at $R=S[Z]_{<M_S,Z>}$, for $1\leq i \leq r$. Set
$f(Z)=f_{c_1}\cdot f_{c_2}\cdots f_{c_r}\in S[Z]$, $b=c_1+\cdots +
c_r$, and $ \pi:Spec(S[Z]/\langle f(Z) \rangle)\to Spec(S)$. Then:

1) Locally at $x$ the closed set $Sing(\mathcal{G})$ is included
in the set of points of multiplicity $b$ of the hypersurface
$V(\langle f(Z)\rangle )$.

2)  $\pi(Sing(\mathcal{G}_x))\subset
Sing(\overline{\mathcal{R}}_{f_{c_1},f_{c_2},\dots,f_{c_r}})$
\end{lemma}
\begin{proof} Note that $f(Z)\in I_b$, is an element of order $b$
in $R=S[Z]_{<M_S,Z>}$ and hence in $\mathcal{O}_{V,x}$, so 1)
holds by Remark \ref{3rkdelsing}. In order to prove 2) it suffices
to show that the set of points of multiplicity $b$ of $V(\langle
f(Z)\rangle )$ map into
$Sing(\overline{\mathcal{R}}_{f_{c_1},f_{c_2},\dots,f_{c_r}})$.
This follows from part ii) in Theorem \ref{vramloc}, and Prop
\ref{3propsing}, 1).
\end{proof}
\begin{lemma}\label{3lemdos} Assume here that $\mathcal{G}=\bigoplus I_k\cdot W^k$ is a Diff-algebra, and set, as before a simple point $x\in Sing(\mathcal{G})$.
Then elements $f_{c_1},\dots, f_{c_r}$ can be chosen, as in the previous
Lemma, so that $$\pi(Sing(\mathcal{G}_x))=
Sing(\overline{\mathcal{R}}_{f_{c_1},f_{c_2},\dots,f_{c_r}}).$$
\end{lemma}
\begin{proof} It suffices to prove that, for suitable
${f_{c_1},f_{c_2},\dots,f_{c_r}}$ as above,
$Sing(\overline{\mathcal{R}}_{f_{c_1},f_{c_2},\dots,f_{c_r}})\subset
\pi(Sing(\mathcal{G}))$. For technical reasons, to be used later, we
first choose and fix an element $f_{c_1}$ of order $c_1$ in
$I_{c_1}$. Since $\mathcal{G}_x$ is finitely generated, there is
an integer $n_0$ so that $I_{n_0}$ has order $n_0$ at
$\mathcal{O}_{V,x}$, and such that $\mathcal{G}_x$ is an integral
extension of $\mathcal{O}_{V,x}[I_{n_0}W^{n_0}] (\subset
\mathcal{O}_{V,x}[W])$. In this case
$V(I_{n_0})=Sing(\mathcal{G}_x )$ (Prop \ref{3propsing}, 5)), and
there are elements $f_{c_2},\dots,f_{c_r}\in I_{n_0}$, all of
order $n_0$, that generate the ideal $I_{n_0}$ in a neighborhood
of the point.

Recall that
$Sing(\overline{\mathcal{R}}_{f_{c_1},f_{c_2},\dots,f_{c_r}})$
consists of all primes $P$ in $S$, such that, at $S_P$:
$$\nu_{S_P}(G_{c_1,\dots, c_r}^{l}(a^{(1)}_1,\dots,a^{(1)}_{c_1},
\dots ,a^{(r)}_1,\dots,a^{(r)}_{c_r})))\geq deg \ G_{c_1,\dots,
c_r}^{l}\hskip 1cm (\mbox{see }\ref{4rksindelp}).$$ So if $P\in
Sing(\overline{\mathcal{R}}_{f_{c_1},f_{c_2},\dots,f_{c_r}})$,
Theorem \ref{vramloc}, i) asserts that
$$Spec(S[Z]/<f(Z)>)\stackrel{\pi}{\longrightarrow}Spec(S)$$ is
purely ramified over $P$. In particular $\pi^{-1}(P)=Q$ ( is a unique
point in $Spec(S[Z]/\langle f(Z)\rangle )$. Here
$f(Z)=f_{c_1}\cdot f_{c_2}\cdots f_{c_r}$, and $Spec(S[Z]/\langle
f_{c_i}(Z)\rangle)$ is  closed in $Spec(S[Z]/\langle f(Z)\rangle
)$, and maps surjectively into $Spec(S)$. Identify $Q$ with a
prime in $S[Z]$, say $Q$ again, so $$Q \in V(\langle
f_{c_1}\rangle)\cap  V(\langle f_{c_2}\rangle)\cap \dots \cap
V(\langle f_{c_r}\rangle).$$ In particular, $$\langle
f_{c_2},\dots,f_{c_r}\rangle= I_{n_0}\subset Q,$$ so $Q\in
V(I_{n_0})=Sing(\mathcal{G})$, and hence $P \in
\pi(Sing(\mathcal{G}))$ as was to be shown.
\end{proof}

\begin{definition}\label{defmathcalR}

We now define the {\em elimination algebra} of $\mathcal{G}$
relative to a transversal projection $\pi$, under the assumption that 
$\mathcal{G}$ is a $\pi$-relative differential algebra (\ref{3def1}).

The elimination algebra, say: $\mathcal{R}_{\mathcal{G}}$, is defined as
the smallest subalgebra of $ S[W]$, containing all (elimination)
algebras $\overline{\mathcal{R}}_{f_{c_1},f_{c_2},\dots,f_{c_r}},$
for all choices of $r$, and of elements $f_{c_i}(Z)\in
I_{({c_i})}$ monic of degree $c_i$ in $S[Z]$. Note that, as
graded subalgebra, it can be expressed in terms of ideals $J_k$ in
$S$, namely:

$$\mathcal{R}_{\mathcal{G}}=\bigoplus J_k\cdot W^k (\subset
S[W])$$ for suitable ideals $J_k$ in $S$.

\end{definition}

There is a natural inclusion of graded algebras
$\overline{\mathcal{R}}_{f_{c_1},f_{c_2},\dots,f_{c_r}}\supset
\overline{\mathcal{R}}_{f_{c_2},\dots,f_{c_r}}.$ So if we fix
$f_{c_1}(Z)\in I_{({c_1})}$, monic of degree $c_1$ in $S[Z]$, we
may also define $\mathcal{R}_{\mathcal{G}}$ as the smallest
subalgebra containing all those of the form
$\overline{\mathcal{R}}_{f_{c_1},f_{c_2},\dots,f_{c_r}}$
(including the fixed element $f_{c_1}$).

On the other hand we can define $B=S[Z]/\langle f_{c_1}(Z)
\rangle$, and consider the algebra induced by restriction of
$\mathcal{G}$, say: $$\overline{\mathcal{G}}=\bigoplus
\overline{I}_k\cdot W^k (\subset B[W]),$$ where
$\overline{I}_k=I_kB$.

\begin{theorem}\label{4th10} Fix an algebra $\mathcal{G}$ over a smooth scheme $V$, a closed point $x\in Sing(\mathcal{G})$, and
$\pi: V\to V'$ transversal at $x$.
Set $R= \calo_{V,x}$ and $S=\calo_{V',\pi(x)}$ as before, namely $R=S[Z]_{<M_S,Z>}$. 
Assume that $\mathcal{G}_x$ is a differential algebra relative to $S$ (e.g., that $\mathcal{G}_x$ is a Diff-algebra)(\ref{IIdef1}).
%and a simple algebra $\mathcal{G}=\bigoplus I_k\cdot W^k (\subset R[W])$.

Fix $f_{c_1}(Z)\in I_{{c_1}}$, monic of degree $c_1$ in $S[Z]$.
Set $B=S[Z]/\langle f_{c_1}(Z)\rangle $,
$\overline{\mathcal{G}}\subset B[W]$ as above, and $\pi:
Spec(S[Z]/\langle f_{c_1}(Z)\rangle )\to Spec(S) \rangle$ (restriction of $\pi$). Then, locally at $x$:
\begin{enumerate}
\item[(i)] $Sing(\mathcal{G}) \subset V(Diff^{c_i-1}(\langle  f_{c_1}(Z)
\rangle ))$, and
$\pi(Sing(\mathcal{G}))\subset Sing(\mathcal{R}_{\mathcal{G}})$.
Moreover, $$\pi(Sing(\mathcal{G}))=Sing(\mathcal{R}_{\mathcal{G}})$$
if $\mathcal{G}$ is a Diff-algebra.

\item[(ii)] The elimination algebra $\mathcal{R}_{\mathcal{G}}$ is
included in $\overline{\mathcal{G}}\cap S[W]$ (as subalgebras of
$S[W]$).

\item[(iii)] $\overline{\mathcal{G}}$ is integral over
$\mathcal{R}_{\mathcal{G}}$ (in particular
$\overline{\mathcal{G}}\cap S[W]$ is integral over
$\mathcal{R}_{\mathcal{G}}$).

\item[(iv)] The algebra $\overline{\mathcal{G}}\cap S[W]$ is, up to
integral closure, independent of the choice of $f_{c_1}(Z)\in
I_{{c_1}}$.

\item[(v)] If $\mathcal{G}\subset \mathcal{G}'$ is a finite extension,
then $\mathcal{R}_{\mathcal{G}}\subset \mathcal{R}_{\mathcal{G}'}$
is also finite.
\end{enumerate}
\end{theorem}
\proof

i) The first inclusion is \ref{3rkdelsing}. The equality follows
from Lemmas \ref{3lemuno} and \ref{3lemdos}.
\medskip

 ii)It suffices to show that each algebra
$\overline{\mathcal{R}}_{f_{c_1},f_{c_2},\dots,f_{c_r}}$
($f_{c_1}$ as above) is included in $\overline{\mathcal{G}}\cap
S[W]$ as graded algebra. This will follow, on the one hand from
\ref{cordecmul}; and, on the other hand, by the fact that
$\mathcal{G}$ is closed by the action of differential operators in
the variable $Z$. In fact, recall that
$\overline{\mathcal{R}}_{f_{c_1},f_{c_2},\dots,f_{c_r}}$ was
defined in terms of the universal polynomials
$F_{c_1},F_{c_2},\dots,F_{c_r}$, and the Rees algebra
$\overline{R}_{{c_1},{c_2},\dots,{c_r}}$ (see (\ref{1eqdrbar}),
and \ref{1univuno}). Fix an homogeneous element of degree $m$, say
$G_m$, in $\overline{R}_{{c_1},{c_2},\dots,{c_r}}$. $G_m$ can be
expressed as a polynomial in $ \{ F^{(j)}_{c_i}, 0\leq j \leq c_i,
1\leq i \leq r \},$ ($F^{(j)}_{c_i}$ defined in terms of
differential operators) with coefficients in the field $k$ .
Furthermore, such expression of $G_m=G_m(F^{(j)}_{c_i})$ is
weighted homogeneous, provided $F^{(j)}_{c_i}$ is given weight
$c_i-j$ (see Corollary \ref{cordecmul}, 2)).

The elements $f_{c_1},f_{c_2},\dots,f_{c_r}$ are defined from
$F_{c_1},F_{c_2},\dots,F_{c_r}$ by specialization, and each $f_{c_i}$
is homogeneous of degree $c_i$ in $\mathcal{G}$ (\ref{1eqdrbar}). Therefore the
elements $f^{(j)}_{c_i}$ (defined in terms of differential
operators) are homogeneous of degree $c_i-j$ in the Diff-algebra
$\mathcal{G}$; and $G_m(f^{(j)}_{c_i})$ (image of
$G_m(F^{(j)}_{c_i})$) is homogeneous of degree $m$ in
$\mathcal{G}$.
 This proves that
$\overline{\mathcal{R}}_{f_{c_1},f_{c_2},\dots,f_{c_r}} \subset
\overline{\mathcal{G}}\cap S[W]$ as graded algebras, since
$\overline{\mathcal{R}}_{f_{c_1},f_{c_2},\dots,f_{c_r}}$ is the
graded algebra generated by all $G_m({f}^{(j)}_{c_i})$.
\medskip

iii) Choose a positive integer $n_0$ with two conditions. First,
that $V^{(n_0)}(\mathcal{G})$ be a usual Rees ring defined by
the ideal $I_{n_0}$ (i.e.,
$V^{(n_0)}(\mathcal{G})=\oplus_k(I_{n_0})^kW^{kn_0}).$ And second,
that the order of $I_{n_0}$ at the local ring $R$ is $n_0$. Such
choice of $n_0$ and $I_{n_0}$ is possible since $\mathcal{G}$ is
finitely generated and simple.

The ideal $I_{n_0}$ can be generated by elements of order $n_0$ in
the local regular ring $R$; and, replacing $R$ by its completion,
we may assume that it is generated by monic polynomials, say
$f_2(Z), \dots f_r(Z)$ in the variable $Z$ ($I_{n_0}=<f_2(Z),
\dots, f_r(Z)>$).

Recall that $B=S[Z]/\langle f_{c_1} \rangle $, and set
$\overline{I}_n=I_nB$, and $\overline{f}_i\in I_nB$ as the class
of ${f}_i$. In order to prove that $\overline{\mathcal{G}}$ is
finite over the subalgebra $\mathcal{R}_{\mathcal{G}}$, it
suffices to prove that the elements $\overline{f}_i$ are
integral over $\mathcal{R}_{\mathcal{G}}$ (see Remark
\ref{rk4i5}). Observ that:

a) the elements $F_{c_i}(Y_1)^{}$, $ 1\leq i \leq r$ are integral
over $\overline{R}_{c_1,c_2,\dots,c_r}$ (see Corollary
\ref{cordecmul}, 1)).

b)  $F_{c_i}(Y_1)^{}$, $ 1\leq i \leq r$ are elements in
$k[s^{(1)}_1,\dots,s^{(1)}_{c_1},\dots
s^{(r)}_1,\dots,s^{(r)}_{c_r}][Z]/<F_{c_1}(Z)>$ (see Lemma
\ref{1lemdederi}).

c) $B$ and ${\overline{\mathcal{R}}}_{f_{c_1},f_{c_2},\dots,f_{c_r}}$ are
defined from $k[s^{(1)}_1,\dots,s^{(1)}_{c_1},\dots
s^{(r)}_1,\dots,s^{(r)}_{c_r}][Z]/<F_{c_1}(Z)>$ and from
$\overline{R}_{c_1,c_2,\dots,c_r}$ by base change:
$k[s^{(1)}_1,\dots,s^{(1)}_{c_1},\dots
s^{(r)}_1,\dots,s^{(r)}_{c_r}]\to S ( \ref{1def29} ).$

This shows that the elements $\overline{f}_i$ are integral over
${\overline{\mathcal{R}}}_{f_{c_1},f_{c_2},\dots,f_{c_r}}$, and
hence over $\mathcal{R}_{\mathcal{G}}$.

\medskip

(iv) Follows from (iii), since $\mathcal{R}_{\mathcal{G}}$ is
defined independently of the choice of $f_{c_1}(Z)\in I_{{c_1}}$.

\medskip

(v) The image at $B=S[Z]/\langle f_{c_1}(Z) \rangle$ defines
$\overline{\mathcal{G}} \subset \overline{\mathcal{G}}' $, which
is also a finite extension, so the claim follows from (iii).
\endproof

The elimination algebra $\mathcal{R}_{\mathcal{G}}$ has been
defined as a direct limit of algebras
$\overline{\mathcal{R}}_{f_{c_1},f_{c_2},\dots,f_{c_r}}$
(\ref{1eqdrbar}). We can also define
 $\mathcal{H}_{\mathcal{G}}$ as a direct limit of
 $\mathcal{H}_{f_{c_1},f_{c_2},\dots,f_{c_r}} (\subset S[W])$ (\ref{1eqdrbar2}).

\begin{corollary}\label{4cor10} The singular locus of the
Diff-algebra $\mathcal{G}_x=\bigoplus I_k\cdot W^k$ maps
bijectively to $Sing(\mathcal{R}_{\mathcal{G}})$, which
coincides with $Sing(\mathcal{H}_{\mathcal{G}})$. In fact
$\mathcal{R}_{\mathcal{G}}$ is a finite extension of
$\mathcal{H}_{\mathcal{G}}$.
\end{corollary}
\begin{proof}
The claim follows from the finiteness of the extension in
(\ref{1eqdrbarmier}).
\end{proof}
\begin{theorem}\label{4th10-1}Fix a 
Rees algebra $\mathcal{G}$, which is simple at a closed point $x\in Sing(\mathcal{G})$, and $\pi: V\to V'$ transversal at $x$. Set $R= \calo_{V,x}$ and $S=\calo_{V',\pi(x)}$ as before, namely $R=S[Z]_{<M_S,Z>}$; and assume that 
$\mathcal{G}_x=\bigoplus I_k\cdot W^k (\subset R[W])$ is a differential Rees algebra relative to $S$. The
elimination algebra $\mathcal{R}_{\mathcal{G}} (\subset S[W])$ is
a graded subalgebra of $\mathcal{G}_x(\subset R[W])$, via the
inclusion $S[W]\subset R[W]$.
\end{theorem}

\proof Recall, as in the proof of the previous corollary, the
definition of  $\mathcal{R}_{\mathcal{G}}$ in terms of algebras
$\overline{\mathcal{R}}_{f_{c_1},f_{c_2},\dots,f_{c_r}}$
(\ref{defmathcalR}). Recall also that
$\overline{\mathcal{R}}_{f_{c_1},f_{c_2},\dots,f_{c_r}}$ is the
pull-back of $\overline{{R}}_{f_{c_1},f_{c_2},\dots,f_{c_r}}$ (see
\ref{1univuno}). The claim follows from Prop \ref{propjun17}, and
Corollary \ref{corjun17}, which show that
$\overline{\mathcal{R}}_{f_{c_1},f_{c_2},\dots,f_{c_r}}$ is
generated by elements which are weighted homogeneous on elements
on $\{\Delta^e(f_{c_i})W^{c_i-e}$, $0\leq e \leq c_i-1\}$; and
hence by homogeneous elements in $\mathcal{G}$.
\begin{corollary} 
1) If $\{f_{c_i}(Z)W^{c_i}; 1\leq i \leq r\}\subset
S[Z][W]$ are chosen so that $$S[Z][\{f_{c_i}(Z)W^{c_i}; 1\leq i \leq
r\}]\subset \mathcal{G}_x$$ is a finite extension, then
$\mathcal{R}_{\mathcal{G}}$
and
$\overline{\mathcal{R}}_{f_{c_1},f_{c_2},\dots,f_{c_r}}$ have the same integral closure (see (\ref{471eq}).

2) If $\mathcal{G}$ is, in addition, a Diff-algebra, then the Rees algebra $\mathcal{R}_{\mathcal{G}}$
is also a Diff-algebra.
\end{corollary}
\proof Claim 1) follows from the two Theorems. As for 2) note that the inclusion $S\subset S[Z]\subset
\calo_{V,x}$ provides an inclusion of higher order differentials
on $S$, as differentials on $ \calo_{V,x}$. \endproof

%%%%%%%%%%%%%%%%%%%%%%%%%%%%%%%%%%%%%%%%%%%%%%%%%%%%%%%%%%%%%%%%%%%%%%%%%%%%%%%%%%%%%%%%%%%%%%%%%%%%%%%%%%%%%%%%

%%%%%%%%%%%%%%%%%%%%%%%%%%%%%%%%%%%%%%%%%%%%%%%%%%%%%%%%%%%%%%%%%%%%%%%%%%%%%%%%%%%%%%%%%%%%%%%%%%%%%%%%%%%%%%%%%%%

%%%%%%%%%%%%%%%%%%%%%%%%%%%%%%%%%%%%%%%%%%%%%%%%%%%%%%%%%%%%%%%%%%%%%%%%%%%%%%%%%%%%%%%%%%%%%%%%%%%%%%%%%%%%%%%%%%%%%%

%%%%%%%%%%%%%%%%%%%%%%%%%%%%%%%%%%%%%%%%%%%%%%%%%%%%%%%%%%%%%%%%%

\section { On differential invariants and projections.}

 The local ring $\calo_{V,x}$ of the smooth scheme
$V$ at $x$, is regular. There is a well
defined notion of order for ideals in a local regular ring. So a
sheaf of ideals, say $I\subset \calo_{V}$, defines a function on
$V$ with values in the integers, by considering, at each point $x\in V$ the
order of the ideal at $ \calo_{V,x}$. We first extend this notion of order
to the case of
 Rees algebras $\mathcal{G}$ over $V$. In Theorem \ref{th55}, which is the main result in this Section,  we study the behavior of the 
order with our notion of elimination of one variable. This elimination is defined in
terms of a projection. There can be different projections and the Theorem studies the order function of the elimination algebras for different projections.

In Proposition \ref{4th11} we show that Hironaka´s
$\tau$-invariant has the expected behavior in positive
characteristic when considering elimination of one variable (as
the known behavior of the $\tau$-invariant in characteristic zero).

\begin{parrafo}\label{agre1}
{\em The notion of Rees algebras $\mathcal{G}=\bigoplus_{k\geq 1}
I_k\cdot W^k$ parallels that of idealistic exponents in
\cite{Hironaka77}, and the notion of singular locus
$Sing(\mathcal{G})$, is the natural analog of that defined for
idealistic exponents. In what follows assume that $I_k$ is
non-zero for some index $k>0$.

We recall the definition of a function, which is the natural
reformulation of that defined by Hironaka for idealistic exponents. We
follow here the presentation of section 6 in \cite{VV3}. Fix $x\in
\Sing(\mathcal{G})$. Given $f_nW^n\in I_nW^n$, $f_n\neq 0$, set
$ord_x (f_n)=\frac{\nu_x(f_n)}{n}\in \mathbb{Q};$ called the order
of $f_n$ (weighted by $n$), where $\nu_x$ denotes the order at the
local regular ring $\calo_{V,x}$. As $x\in \Sing(\mathcal{G})$ it
follows that $ord_x (f_n)\geq 1.$ Define
$$ord_x (\mathcal{G})= inf_{}\{ord_x (f_n); f_nW^n\in I_nW^n\};$$
or equivalently, $ord_x (\mathcal{G})= inf_{n\geq
0}\{\frac{\nu_x(I_n)}{n}\}$. In general $ord_x (\mathcal{G})\geq
1$ for all $x\in \Sing(\mathcal{G})$; and  a Rees algebra
$\mathcal{G}$ is simple at $x$ if and only if $ord_x
(\mathcal{G})= 1$ (\ref{defsmpt}). }
\end{parrafo}
\begin{proposition}\label{3propsingZ}
\begin{enumerate}
\item[(1)] If $\mathcal{G}$  is  generated over $\calo_V$ by
$\mathcal{F}=\{ g_{n_i}W^{n_i}, n_i>0 , 1\leq i\leq m \}$, then
$ord_x (\mathcal{G})= inf_{}\{ord_x (g_{n_i}); 1\leq i\leq m\}.$
And if $N$ is a common multiple of all $n_i , 1\leq i\leq m$, then
$ord_x (\mathcal{G})=\frac{\nu_x(I_N)}{N}$.

\item[(2)] If $\mathcal{G}$ and $ \mathcal{G}'$ are Rees algebras with the
same integral closure (e.g., if $\mathcal{G}\subset \mathcal{G}'$
is a finite extension), then, for all $x\in
Sing(\mathcal{G})(=Sing(\mathcal{G}'))$: $ord_x
(\mathcal{G})=ord_x (\mathcal{G}').$

\item[(3)] Let $\mathcal{G}''=\bigoplus I''_n\cdot W^n$ be the extension
of $\mathcal{G}$ to a differential Rees algebra relative to $k$,
as defined in Theorem \ref{3th1}, then for all $x\in
Sing(\mathcal{G})(=Sing(\mathcal{G}''))$: $ord_x
(\mathcal{G})=ord_x (\mathcal{G}'').$
\end{enumerate}
\end{proposition}
All assertions are easy to check . The assumption in (1) holds at
an affine open set of $V$.
\begin{parrafo}\label{par410}{\em Fix a Diff-algebra $\mathcal{G}=\bigoplus_{k\geq
1} I_k\cdot W^k$ over a smooth scheme $V$ and a closed point $x\in
Sing(\mathcal{G})$. Assume that $x$ is a simple point, or
equivalently, that $ord_x (\mathcal{G})=1$. In such a case there
must be an index $n$ so that $\nu_x(I_n)=n$ (order of $I_n$ at
$\calo_{V,x}$). In other words, there is an homogeneous element of
degree $n$, say $f_n\cdot W^n\in \mathcal{G}=\bigoplus_{k\geq 1}
I_k\cdot W^k$, so that $f_n$ has order $n$ at $\calo_{V,x}$. We
claim now that if $\mathcal{G}$ is integrally closed (i.e., equal
to its integral closure in $\calo_{V}[W]$), then $f_n$ can be
chosen to be analytically irreducible at the local regular ring
$\calo_{V,x}$. This would show, in particular,  that
$f_{c_1}(Z)\in S[Z]$ can be chosen to be irreducible in Theorem
\ref{4th10}. This fact, already interesting in itself, will be used in our forthcoming Theorem \ref{th55}.

Note here that if $\mathcal{G}'$ denotes the integral closure of
$\mathcal{G}$ in $\calo_{V}[Z]$, then $\mathcal{G} \subset
\mathcal{G}'$ is a finite extension. In Theorem  6.12 of
\cite{VV1}, it is proved that if $\mathcal{G} \subset
\mathcal{G}'$ is a finite extension of Rees algebras, there is an
inclusion of the Diff-algebras spanned by each of them, and this extension is also
finite. In particular, if $\mathcal{G}$ is already a Diff-algebra,
its integral closure is also a Diff-algebra.

Recall that for any smooth morphism $\pi_1: V\to V^{(1)}$,
$\mathcal{G}$ is also a differential Rees algebra relative to  $\pi_1$
(see Observation (2) in \ref{4rksindelp}). 
This relative differential structure is as much as we use in order to the define the elimination algebra
 $\mathcal{R}_{\mathcal{G}}^{}$ (\ref{defmathcalR}).
 
 Let us first discuss the claim, that
 $f_{c_1}(Z)\in S[Z]$ can be chosen to be irreducible, 
 under the assumption that $ \mathcal{G}$ is an $S$-relative differential algebra. Set $S={\calo}_{V^{(1)},x}$, and let $m_S$ be the maximal ideal.  Suppose that a polynomial $S[Z]$ is defined so that 
$S[Z]\subset {\calo}_{V,x}$, and $S[Z]_{(m_S,Z)}\subset {\calo}_{V,x}$ is an \'etale extension of local rings.
 
 Assume that there is an element $f_nW^n\in \mathcal{G}$, with $f_n$ of order $n$ at $\calo_{V,x}$, and that up to multiplication by a unit 
 $f_n=F_n(Z)\in S[Z]$ is a monic polynomial of degree $n$.
 
 Set $F_n(Z)=G_{r_1}\cdot G_{r_2}\cdots G_{r_s}$, a product of
irreducible polynomials in $S[Z]$, where each $G_{r_i}(Z)$ is monic
of degree $r_i<n$. In this case, it follows from
\ref{uni660} that each $G_{r_i}(Z)\cdot W^{r_i}$ is in the
integral closure of $\mathcal{G}$. As we assume that $\mathcal{G}$ is integrally closed, $G_{r_i}(Z)\cdot W^{r_i}\in  \mathcal{G}$ and 
$G_{r_i}(Z)$ has order $r_i (<n)$ at ${\calo}_{V,x}$. If $r_i=1$ then 
$G_{r_i}(Z)$ is analytically irreducible in this local ring.

We prove our claim, in full generality, by using the previous argument. 
At this point we shall make use of the Weierstrass Preparation Theorem, on the one hand, and also use the fact that, locally at $x\in V$, there are plenty 
of smooth morphisms as $\pi_1: V\to V^{(1)}$. Moreover, we show here that constructing such morphisms is very simple, at least over perfect fields.

\bigskip

% VVVVVVVVVVVVVVVV
%The previous discussion also applies in this relative setting: the integral closure is also a differential Rees algebra relative to  $\pi_1$. 
%That is as much as we need in for the existence of the irreducible 
%${f_1}(Z)\in S[Z]$.

Assume that $f_n$ has order $n$ at $\calo_{V,x}$, and let
$f_n=g_{r_1}\cdot g_{r_2}\cdots g_{r_s}$ be a factorization as a
product of irreducible elements at the henselization of
${\calo}_{V,x}$, say $R$. Then, and after replacing $V$ by a suitable 
\'etale neighborhood of $x$, we may assume that such factorization holds at   ${\calo}_{V,x}$, and that each $g_{r_i}$ is analytically
irreducible. Notice that if $r_i$ denotes the multiplicity of each
$g_{r_i}$ at the local regular ring, then $r_1+r_2+\cdots
+r_s=n$.

The Weierstrass Preparation Theorem holds at henselian local
rings. We claim now that after taking again a suitable \'etale neighborhood, one can construct a subring $S\subset {\calo}_{V,x}$, and a polynomial
ring $S[Z]\subset {\calo}_{V,x}$, so that up to multiplication by  units,
$f_n=F_n(Z)$ and each $g_{r_i}=G_{r_i}(Z)$, where
$F_n(Z)=G_{r_1}\cdot G_{r_2}\cdots G_{r_s}$ is a product of
irreducible polynomials in $S[Z]$, and each $G_{r_i}(Z)$ is monic
of degree $r_i$. 
To check this claim note first that one can choose a regular system of parameters, say $\{x_1, \dots, x_d\}$ for ${\calo}_{V,x}$, so that 
$f_n$ also has order $n$ when restricted to ${\calo}_{V,x}/\langle x_1, \dots ,x_{d-1}\rangle$. Since $V$ is smooth over $k$, ${\calo}_{V,x}$ contains the polynomial ring $k[x_1, \dots, x_d]$.

As $k$ is a perfect field, and $x$ is a closed point, after suitable finite extension of the base field, ${(V,x)}$ is an 
\'etale neighborhood of the affine space ${\mathbb A}^d$ at the origin
(see \cite{AK}). 
So there is an inclusion of local rings, say 
$$k[x_1, \dots, x_{d-1}]_{\langle x_1, \dots , x_{d-1} \rangle}\subset {\calo}_{V,x},$$
and a smooth morphism, say $V\to {\mathbb A}^{d-1}$, defined at a suitable neighborhood of $x$. Let $X\subset V$ be the hypersurface 
$V(\langle f_n \rangle )$, and set $X\to {\mathbb A}^{d-1}$ by restriction. Notice that $x$ is an isolated point in the fiber over the origin of ${\mathbb A}^{d-1}$. Under these conditions an expression $F_n(Z)=G_{r_1}\cdot G_{r_2}\cdots G_{r_s}$, as above, can be defined in $S[Z]$, where $S$ denotes the henselization of $k[x_1, \dots, x_{d-1}]_{\langle x_1, \dots , x_{d-1} \rangle}$ (see \cite{R}). 

As each $G_i(Z)$ involves finitely many
coefficients in $S$, we may also assume that all statements hold
al a suitable \'etale neighborhood of ${\mathbb A}^{d-1}$ at the origin.
Take, finally, an \'etale neighborhood
of the closed point $x\in V$ containing the previous local ring.
%
%As every
%differential Rees algebra over the field $k$, is also a relative
%differential Rees algebra (over $S$); 

Therefore, at a suitable \'etale neighborhood of $x$ we may assume that $G_{r_i}(Z)\cdot W^{r_i}$ is in the
integral closure of $\mathcal{G}$ 
(see (\ref{eq660})). This proves the claim.  Namely that if we pass from $\mathcal{G}$ to its integral
closure, we may assume that the element $f_{c_1}(Z)$, in Theorem
\ref{4th10}, 
 is irreducible in
$S[Z]$ and, moreover, analytically irreducible in $\calo_{V,x}$.
% at a suitable \'etale neighborhood of the point. 
}
\end{parrafo}
\begin{parrafo}\label{par411}{\em Given a Diff-algebra $\mathcal{G}$, a simple point $x\in Sing(\mathcal{G})$, and a general projection,
then an elimination algebra $\mathcal{R}_{\mathcal{G}}$ has been
defined over a regular scheme $Spec(S)$. The previous discussion
shows how to construct transversal smooth morphisms over
a field. 
%In other words, that the elimination algebras are 
%defined within the class of smooth schemes over a field. 
In fact,
let $d$ denote the dimension of the smooth scheme $V$, then, after restriction to a suitable
\'etale neighborhood of $x$, there is a smooth scheme $V^{(1)}$ of
dimension $d-1$, and a smooth morphism $\pi_1: V\to V^{(1)}$, so
that $\mathcal{R}_{\mathcal{G}}^{(1)}$ (\ref{defmathcalR}) can be
defined at $V^{(1)}$. 
%So $S$ can be taken to be
%$\calo_{V^{(1)},\pi(x)}$ in Theorem \ref{4th10}.

Furthermore, we may take $f_{c_1}(Z)$ to be a global section at a  neighborhood of $x$, and define
$\mathcal{R}_{\mathcal{G}}^{(1)}\subset \calo_{V^{(1)}}[W]$. Moreover, if $\mathcal{G}$ is a Diff-algebra, 
then
$\pi_1(Sing(\mathcal{G}))=Sing(\mathcal{R}_{\mathcal{G}}^{(1)})$.
The setting of Theorem \ref{4th10} holds at every closed point
$y\in Sing(\mathcal{G})$ by taking $S=\calo_{V^{(1)},\pi_1(y)}$.
As $Sing(\mathcal{G})$ is included in the set of points of multiplicity
$c_1$ of the hypersurface defined by $f_{c_1}(Z)$, it follows that
the components of $Sing(\mathcal{G})$ of codimension one in $V$
are open and closed in $Sing(\mathcal{G})$, and furthermore, they
are also smooth. In fact this is the case for the highest
multiplicity locus of a hypersurface in a smooth scheme (see
\cite{EncVil97:Tirol}, Ex. 13.11, where a characteristic free
proof is suggested). 

We will assume here that $Sing(\mathcal{G})$ has
codimension at least 2 locally at the simple point. This ensures that
$\mathcal{R}_{\mathcal{G}}^{(1)}$ is a direct sum of non-zero
ideals in $ S=\calo_{V^{(1)},\pi_1(y)}$, and in this case
$ord_{\pi_1(x)}(\mathcal{R}_{\mathcal{G}}^{(1)})$ is defined
(\ref{agre1}).

We took, as stating point, that  $\mathcal{G}$ was a Diff-algebra.
This of course ensures that for any smooth morphism $\pi_1: V\to V^{(1)}$, $\mathcal{G}$ as also a differential Rees algebra relative to $\pi_1$
(see Observations 2) in \ref{4rksindelp}). However this latter condition is all we need to define $\mathcal{R}_{\mathcal{G}}^{(1)}$.

\begin{theorem} \label{th55} Fix a Diff-algebra $\mathcal{G}$, a simple closed point $x\in Sing(\mathcal{G})$, and assume that the local codimension of this closed set ( in $V$) is at least 2.
Let $\pi_1: V\to V^{(1)}$ and $\pi_2: V\to V^{(2)}$ be, as above,
two morphisms of smooth schemes ($dim(V^{(1)})=dim(V^{(2)})=d-1$);
defining elimination algebras, say
$\mathcal{R}_{\mathcal{G}}^{(1)} \subset \calo_{V^{(1)}}[W]$ and
$\mathcal{R}_{\mathcal{G}}^{(2)} \subset \calo_{V^{(2)}}[W]$. Then, 
$$
ord_{\pi_1(x)}(\mathcal{R}_{\mathcal{G}}^{(1)})=ord_{\pi_2(x)}(\mathcal{R}_{\mathcal{G}}^{(2)}).$$
\end{theorem}

This shows that $ord_{\pi(x)}(\mathcal{R}_{\mathcal{G}}^{})$ is an
invariant of $x\in Sing(\mathcal{G})$  (i.e., independent of the
projection $\pi$). We will introduce some notation, and discuss
some preliminary results before we address the proof in \ref{55e}.
}
\end{parrafo}

\begin{parrafo} \label{55a} {\em  Let $k[[x_1,\dots,x_d]]$ be the ring of formal power series over a field $k$; and let $f_{c}$ an irreducible
element of multiplicity $c$. Let $B$ denote the quotient
$k[[x_1,\dots,x_d]] / \langle f_c \rangle$, which is a domain with quotient field, say $L$.
The Weierstrass Preparation Theorem asserts that for a sufficiently general choice of coordinates we may assume that, up to
multiplication by a unit, $f_c$ is a monic polynomial of degree $c$ in the variable $x_d$, and that the class of
$ x_1,\dots,x_{d-1}$  in $B$  span a reduction of the maximal ideal.

Here $B$ is a finite extension, and a free module of rank $c$, over the subring of formal power series, say $S= k[[x_1,\dots,x_{d-1}]] $.
Let $K$ denote the total quotient field of $S$, and note that $L$ is a finite extension of degree $c$ over $K$.

For each discrete valuation ring, say $V$, in $L$, we consider the
restriction, say $$ V_K=V\cap K.$$ Identify the group of values
of  $V$ with the integers $\mathbb{Z}$, , and define the { \em
ramification index} of $V$ over $V_K$ to be the index in
$\mathbb{Z}$ of the subgroup of values of $V_K$.

Let $Spec(B) \leftarrow F$ denote the normalized blowup of $B$ at the maximal ideal. Let $H_1, H_2,
\dots, H_l$ denote the irreducible exceptional hypersurfaces of $F$, and let $V_1, V_2,
\dots, V_l$ denote the discrete valuation rings in $L$ , where each $V_i$  is  the local ring  of the normal scheme $F$  at the generic point of  the hypersurface $H_i$.

Let $Spec(S) \leftarrow Y$  denote the blow up of $S$ at the
maximal ideal (i.e., the quadratic transformation). So here $Y$ is
regular, and has only one exceptional hypersurface, say $h$, let $V_S$
be the local ring of $Y$ at the generic point of $h$. $V_S$ is the
valuation at $K$ which extends the order at the local regular ring
$S$, for elements in $S$.

\begin{lemma}\label{55b} For each discrete valuation ring $V_i$ , $i=1,...,l$  as above, the ramification index of $V_i$ over $(V_i)_K$ is the
order at $V_i$ of the maximal ideal of $B$.
\end{lemma}
\proof Let $Spec(B) \leftarrow  \overline{Y}$ denote the fiber product of $Spec(S) \leftarrow Y$ with the finite morphism
 $Spec(B) \to Spec(S)$. Note that $Spec(B) \leftarrow  \overline{Y}$ is also the blow up of $B$ at the ideal spanned by the
 elements $x_1,\dots,x_{d-1}$, say $\langle x_1,\dots,x_{d-1} \rangle B$.
So $\overline{Y}$ is a finite extension of $Y$, and the total
quotient field of $\overline{Y}$ is $L$. Since $\langle
x_1,\dots,x_{d-1} \rangle B$ is a reduction of the maximal ideal
of $B$, it follows that $F$ is the normalization of
$\overline{Y}$, and that $(V_i)_K=V_S$ for all $i=1,\dots ,l$.
Furthermore, as the ideal spanned by $x_1,\dots,x_{d-1}$ has order
one at $V_S$, it also follows that the ramification index of $V_i$
over $V_S$ is the order of  the ideal spanned by
$x_1,\dots,x_{d-1}$  at $V_i$. But this is the order of the
maximal ideal of $B$ at $V_i$. In fact, the maximal ideal and
$\langle x_1,\dots,x_{d-1} \rangle B$ have the same integral
closure.

\begin{remark}\label{55c} Fix notation as above, and let  $e_i$ denote the ramification index of $V_i$ over $(V_i)_K=V_S$, for $i=1, \dots ,l$.  Let $J$ be an ideal in the local regular ring $S$.
The order of $J$ at $S$ is the valuation of $J$ at $V_S$, say $b\in \mathbb{Z}$.  It follows from Lemma \ref{55b}  that the order of the extended ideal $JB$ at the valuation $V_i$ is  the integer $b\cdot e_i$.

\end{remark}

\begin{corollary} \label{55d} The ramification index of each $V_i$ over $(V_i)_K$ is independent of the choice of
$x_1,\dots,x_{d-1}$  (i.e., of  $S= k[[x_1,\dots,x_{d-1}]] \subset B$), as far as $\langle x_1,\dots,x_{d-1} \rangle B$ is a reduction of the maximal ideal of $B$.
\end{corollary}

}
\end{parrafo}

\begin{parrafo}\label{55e} Proof of Theorem \ref{th55}: {\em  Here $\mathcal{G}=\bigoplus
{I}_k\cdot W^k (\subset \calo_V[W])$ is a Diff-algebra and $x\in
Sing(\mathcal{G})$ is a simple point. Set $\pi_1: V\to V^{(1)}$
and $\mathcal{R}_{\mathcal{G}}^{(1)}=\bigoplus J^{(1)}_k\cdot W^k
(\subset \calo_{V^{(1)}}[W])$ for suitable ideals $J^{(1)}_k$ in
$\calo_{V^{(1)}}$. Let $\mathcal{G}_x=\bigoplus ({I}_k)_x\cdot W^k
(\subset \calo_{V,x}[W])$ be the localization at $x\in V$. As the
point $x$ is simple there must be an index $c$ and an element
$f_c\in ({I}_c)_x$ of order $c$ at $\calo_{V,x}$. At a suitable
\'etale neighborhood of $x$ and $\pi_1(x)$, $\pi_1$ induces a
finite morphism from the subscheme defined by $\langle f_c
\rangle$ to $V^{(1)}$. In particular, a finite morphism
$$\pi_1:
Spec(\calo_{V,x}/\langle f_{c}\rangle )\to
Spec(\calo_{V^{(1)},x}). $$

The Weierstrass Preparation Theorem asserts that (at a suitable
\'etale neighborhood), setting $S=\calo_{V^{(1)},x}$, there is an
inclusion of regular local rings, say $R=S[Z]_{<M_S,Z>}\subset
\calo_{V,x}$, and a monic polynomial of degree $c$, say $f_c(Z)\in
S[Z]$, so that $$B=\calo_{V,x}/\langle f_{c}\rangle=S[Z]/\langle
f_{c}(Z)\rangle.$$
 Let $\overline{\mathcal{G}}=\bigoplus \overline{I}_k\cdot W^k (\subset
B[W])$ be the algebra induced by restriction of $\mathcal{G}$,
where $\overline{I}_k=I_kB$.

Theorem \ref{4th10}, (iii), states that up to integral closure,
the localization of $\mathcal{R}_{\mathcal{G}}^{(1)}$ at
$\pi_1(x)$ is $\overline{\mathcal{G}}\cap S[W]$. In particular
$ord_{\pi_1(x)}(\mathcal{R}_{\mathcal{G}}^{(1)})=ord_{\pi_1(x)}(\overline{\mathcal{G}}\cap
S[W])$ (\ref{3propsingZ}, (2)).

If $\mathcal{G}\subset \mathcal{G}'$ is a finite extension of
Diff-algebras in $\calo_{V}[W]$, then the restrictions to $B[W]$,
say $\overline{\mathcal{G}}\subset \overline{\mathcal{G}}'$, is
also a finite extension . Therefore $\overline{\mathcal{G}}\cap
S[W]\subset \overline{\mathcal{G}}'\cap S[W]$ is a finite
extension, and hence $ord_{\pi_1(x)}(\overline{\mathcal{G}}\cap
S[W])=ord_{\pi_1(x)}(\overline{\mathcal{G}}'\cap S[W])$
(\ref{3propsingZ}, (2)).

As the integral closure of a Diff-algebra is a Diff-algebra
(\ref{par410}), we may assume, for the proof of this Theorem, that
$\mathcal{G}$ is integrally closed. In particular, as it was indicated in
(\ref{par410}), we may assume here that $f_c\in I_c$ is chosen to
be analytically irreducible. In other words, that the completion
of local ring $B$ is irreducible.

Here  $\pi_1: V\to V^{(1)}$ and $\pi_2: V\to V^{(2)}$ are defined in a neighborhood of $x\in V$. Set $R=\mathcal{O}_{V,x}$. At the tangent space of $x\in V$, say $Spec(gr_M(R))$, the tangent cone of the Diff-algebra $\mathcal{G}$ is a linear subspace, say $T$, in $Spec(gr_M(R))$ (see \ref{partau}). As $\pi_1$ and $\pi_2$ are smooth, their differential maps define one dimensional subspaces,
say $ker (d(\pi_1))$ and $ker (d(\pi_2))$, in $Spec(gr_M(R))$.
As both smooth morphisms are transversal to $\mathcal{G}$ at $x$:
$T$, $ker (d(\pi_1))$ and $ker (d(\pi_2))$, are three subspaces in general position. Using this fact we conclude that the analytically irreducible element $f_c$ can be taken to be transversal to both 
smooth morphisms.

Set $S_1=\calo_{V^{(1)},\pi_1}$ and $S_2=\calo_{V^{(2)},\pi_2}$.
Both play the role of $S$ in the previous discussion. Set now, by
localization at $\pi_i(x)$ :
$$(\mathcal{R}_{\mathcal{G}}^{(i)})_{\pi_i(x)}=\bigoplus J^{(i)}_k\cdot W^k
(\subset S_i[W]) \ (i=1,2).$$

Each $\mathcal{R}_{\mathcal{G}}^{(i)}$  is a  Rees algebras over
the smooth scheme $V^{(i)}$. The claim is that
$$ord_{\pi_1(x)}(\mathcal{R}_{\mathcal{G}}^{(1)})=
ord_{\pi_2(x)}(\mathcal{R}_{\mathcal{G}}^{(2)}).$$

We can choose an integer $N$, so that both
$ord_{\pi_i(x)}(\mathcal{R}_{\mathcal{G}}^{(i)})=\frac{\nu_{\pi_i(x)}(J^{(i)}_N)}{N}$,
for $i=1,2$ (see Prop \ref{3propsingZ},(1)). In particular it
suffices to show that
$\nu_{\pi_1(x)}(J^{(1)}_N)=\nu_{\pi_2(x)}(J^{(2)}_N)$, where
$\nu_{\pi_1(x)}$ denotes the order at the local regular ring $S_i$
(i=1,2).

Within the local ring $B$ there are two regular local
rings $S_1$ and $S_2$. Theorem \ref{4th10},(iii) asserts that the
inclusions
$$(\mathcal{R}_{\mathcal{G}}^{(i)})_{\pi_i(x)}=\bigoplus
J^{(i)}_k\cdot W^k \subset \overline{\mathcal{G}}=\bigoplus
\overline{I}_k\cdot W^k (\subset B[W])$$
are both finite, for $i=1,2$. We can choose the integer $N$ so
that, in addition to the previous conditions,
$V^{(N)}(\overline{\mathcal{G}})$,
$V^{(N)}(\mathcal{R}_{\mathcal{G}}^{(1)})_{\pi_1(x)})$, and
$V^{(N)}(\mathcal{R}_{\mathcal{G}}^{(2)})_{\pi_2(x)})$ are all
Rees rings of ideals. As
$V^{(N)}(\mathcal{R}_{\mathcal{G}}^{(1)})_{\pi_1(x)})\subset
V^{(N)}(\overline{\mathcal{G}})$ is a finite extension of Rees
rings (see Remark \ref {rk4i5}), it follows that $J^{(1)}_N\cdot B$,  and
$\overline{I}_N$ have the same integral closure. A similar
argument proves that the three ideals: $J^{(1)}_N\cdot B$,
$J^{(2)}_N\cdot B$, and $\overline{I}_N$ have the same integral
closure in $B$.

Fix notation as in \ref{55a}, where $Spec(B)\leftarrow F$ denotes
the normalized blow up of $B$ at the maximal ideal, and $V_1,\dots
,V_l$ are valuation rings corresponding to the irreducible
exceptional hypersurfaces in $F$. Set $S=S_1$ and let
$Spec(S)\leftarrow Y$ and $V_S$ also as in \ref{55a}.

Fix a valuation ring among $V_1,\dots ,V_l$, say $V_1$. As $V_1$
dominates the local domain $B$, and $J^{(1)}_N\cdot B$, and
$\overline{I}_N$ have the same integral closure in $B$, it follows
that both ideals have the same valuation at $V_1$.

Let $e_1$ denote the ramification index of $V_1$ over $V_S$.
Remark \ref{55c} says that the valuation of the ideal
$J^{(1)}_N\cdot B$ at the valuation ring $V_1$ is
$\nu_{S_1}(J^{(1)}_N)\cdot e_1$.  Finally Corollary \ref{55d}
ensures that $$\nu_{S_1}(J^{(1)}_N)\cdot
e_1=\nu_{S_2}(J^{(2)}_N)\cdot e_1,$$ as both coincide with the
order of $\overline{I}_N$ at $V_1$. In particular
$\nu_{S_1}(J^{(1)}_N)=\nu_{S_2}(J^{(2)}_N)$.
\endproof

 }
\end{parrafo}

\begin{center}

{\bf On Hironaka's $\tau$-invariant and projections.}

\end{center}
\begin{parrafo}\label{par511}{\em We now discuss a property of elimination of
one variable which parallels well-known properties that hold for inductive arguments used in desingularization over fields of
characteristic zero. To clarify this point recall that the
elimination of a variable, when passing from a Diff-algebra
$\mathcal{G}$ to the Diff-algebra $\mathcal{R}_\mathcal{G}$, is
defined locally at a point $x\in Sing(\mathcal{G})$ when this
point is simple. For simple points we have defined an
homogeneous tangent ideal and a positive integer $\tau=s_r$
(\ref{eqptau}). In the case of characteristic zero $s_0=s_r$, and
the algebra $\mathcal{R}_\mathcal{G}$ will also be simple, unless
$\tau=s_0=1$. In fact if $\tau >1$ at $x\in Sing(\mathcal{G})$,
the invariant $\tau'$ in the elimination algebra is $\tau-1$. The
following result shows that this also holds in the context of
positive characteristic. }
\end{parrafo}
\begin{proposition}\label{4th11}
If in the previous setting $\tau (\mathcal{G}_x)>1,$ then the
elimination algebra $\mathcal{R}_{\mathcal{G}}$ is a simple
Diff-algebra.
\end{proposition}
\proof In what follows we assume that $R$ is the completion of the
local ring at a closed point. Set $\mathcal{G}_x=\oplus
I_nW^n\subset R[W]$. After change of base field, which does not
affect our arguments, we may assume that the closed point is
rational over a perfect field $k$. Let $\{z,x_1,\dots,x_{d-1}\}$
be a regular system of parameters for $M$ (the maximal ideal of
$R$), and set $gr_M(R)=k[Z,X_1,\dots ,X_{d-1}]$ the graded ring
where the variables are the initial forms of the parameters. Over
perfect fields we may assume that the ideal of the tangent cone is
generated by $p$-th powers of linear forms. Assume, for
simplicity, that there are two elements, say $Z^{p^{e'}},
X_1^{p^{e''}} \in k[Z,X_1,\dots ,X_{d-1}]$, in the tangent ideal.
For a suitable $p$-th power, say $n$, there is an element $F_n \in
I_{n}$ ($n=p^e$), and an element $g_n \in I_n$, such that:

i) $In_{p^e}(F_{p^e})=Z^{p^{e}},$ and

ii) $In_{p^e}(g_{p^e})=X_1^{p^{e}}.$

The Weierstrass Preparation Theorem allows us to assume that
$F_{p^e}=F_{p^e}(z)\in S[z],$ is a monic polynomial in the
variable $z$, where $(S,N)$ is the formal power series ring, and
the maximal ideal $N$ is generated by the regular system of
parameters $\{x_1,\dots,x_{d-1}\}$. In fact, multiplication by a
unit modifies the initial form by multiplication by a non-zero
constant in the field. Set
$$F_{p^e}(z)=z^{p^e}+a_1z^{p^{e}-1}+\dots +a_{p^e}\in S[z],$$ and
note that $\nu_S(a_i)>i$, for $1\leq i \leq p^e$. The surjection $S[z]\to
B=S[z]/<F_{p^e}(z)>$ defines a graded ring
$\overline{\mathcal{G}}_x=\oplus \overline{I}_nW^n\subset B[W]$,
where $\overline{I}_n=I_nB$. In particular, the class of
$g_{p^e}$, say $\overline{g}_{p^e}$, is an element of degree
${p^e}$. Set
$$\overline{g}_{p^e}= b_1 \overline{z}^{p^{e}-1}+
b_2\overline{z}^{p^{e}-2}\dots +b_{p^e}\in S[z]/<F_{p^e}(z)>.$$ We
claim that $ \nu_S(b_{i})>i$ for $1\leq i\leq p^e-1$, and
$In_N(b_{p^e})=X_1^{p^{e}}.$ This can be checked at the formal
power series ring $R$. In fact the class is obtained by replacing
the powers $z^{N}$, for all $N>p^e$, in the formal expression of
$g_{p^e}$, by smaller powers. This is done by means of the
relation:$$z^{p^e}=-a_1z^{p^{e}-1}-\dots -a_{p^e}.$$

But $-a_1z^{p^{e}-1}-\dots -a_{p^e}\in M^{p^{e}+1}$, so this
operation does not affect the initial form of $g_{p^e}$.

The ring $S[z]/<F_{p^e}(z)>$ is a free $S$ module, and
multiplication by $\overline{g}_{p^e}$ defines an $S$-linear  endomorphism. So a characteristic
polynomial is assigned to it. The norm, say $$|\overline{g}_{p^e}|=\prod_{1\leq i
\leq p^e} (b_1 \overline{z}_i^{p^{e}-1}+
b_2\overline{z}_i^{p^{e}-2}\dots +b_{p^e}),$$ is defined formally
as this product,
 where each
$\overline{z}_i$ is (formally) the image of $\overline{z}$ by the
different embeddings. We view each factor $(b_1
\overline{z}_i^{p^{e}-1}+ b_2\overline{z}_i^{p^{e}-2}\dots
+b_{p^e})$ with the same weight as $\overline{g}_{p^e}$, namely
$p^e$. The product is an $S$-linear combination on the symmetric
functions on $\overline{z}_i$, which in turn are weighted
functions on the coefficients $a_i$. Since $\nu_S(a_i)> i$, the
order of these elements is higher than the expected order, so
$$In_S(|b_1 \overline{z}^{p^{e}-1}+ b_2\overline{z}^{p^{e}-2}\dots
+b_{p^e}|)=X_1^{(p^e)^2},$$ and the weight of the norm
$|\overline{g_{p^e}}|\in S$ is precisely $(p^e)^2$ (i.e.,
$|\overline{g_{p^e}}|W^{(p^e)^2}\in \mathcal{H}_{\mathcal{G}}$).

This shows that $\mathcal{H}_{\mathcal{G}}$ is simple, and hence
that $\mathcal{R}_{\mathcal{G}}$ is simple (\ref{4cor10}).
\endproof
\begin{remark}\label{rkast} In most resolution problems ideals in smooth
schemes are considered up to integral closure. The natural analog,
in the context of Rees algebras, is to consider them up to
integral closure (see Seccion 5 in \cite{VV3}). For this reason it
is important to introduce invariants of Rees algebras, or of
Diff-algebras, which coincide for algebras with the same integral closure.

Proposition \ref{4th11}, together with Theorem \ref{4th10}, and
Prop \ref{3propsingZ}, 2), already show that the $\tau$-invariant
of a Diff-algebra at a singular point, is the same as that of its
integral closure. In particular the $\tau$-invariant is well
defined up to integral closure.

The following example shows that this is not the case for the
invariant $e_0$, attached to a singular point of a Diff-algebra in
\ref{partau}.

\begin{example} Set $G_1=X^2+Y^5$, $G_2=X$,$G_3=X$, and $F=X^4+X^2Y^5$
in $R=k[X,Y]$, where $k$ is a field of characteristic two. All
polynomials are monic in $X$, and $F=G_1\cdot G_2\cdot G_3$. Let
$\mathcal{G}=\bigoplus I_k\cdot W^k (\subset R[W])$ denote
Diff-algebra spanned by the Rees algebra, generated over $R$, by
the element $F\cdot W^4$, as in Theorem \ref{3th1}. It follows
easily from (\ref{eq3422}) that the $e_0$ invariant at the origin
is $2$, namely that $I_4$ has order $4 (=p^2)$ at the origin; but
the order of $I_i < i$, for $i=1,2,3$, at such point. Let
$\mathcal{G}'=\bigoplus I'_k\cdot W^k (\subset R[W])$ denote the
integral closure. The discussion in \ref{par410} shows that
$X=G_1\in I'_1$, so the $e_0$ invariant of $\mathcal{G}'$ at the
origin is zero.

\end{example}

\end{remark}
%%%%%%%%%%%%%%%%%%%%%%%%%%%%%%%%%%%%%%%%%%%%%%%%%

\section{ Monoidal transformations and differential Rees algebras. Examples.}

\begin{definition}\label{3def072}

 Let  $\mathcal{G}=\bigoplus I_k\cdot W^k\subset \calo_{V}$ be a Rees algebra over a smooth scheme $V$. A
monoidal transformation with smooth center $Y\subset V$, say
$V\longleftarrow V_1$, is said to be {\em permissible} for
$\mathcal{G}$ if $Y\subset Sing(\mathcal{G})$. The exceptional
locus is a smooth hypersurface, say $H\subset V_1$. Since $I_n$ has
order at least $n$ along $Y$,  $I_n \mathcal{O}_{V'}\subset
I(H)^{n}$, for $n\geq 0$. In particular , there is a factorization, say  $I_n
\mathcal{O}_{V_1}=I^{(1)}_n\cdot I(H)^{n},$
 for a unique sheaf of ideals $I^{(1)}_n$, which we call the {\em weighted transform} of $I_n$. The {\em weighted
transform} of $\mathcal{G}$ will be the Rees algebra:
$$(\mathcal{G})_1=\bigoplus I^{(1)}_n\cdot W^k (\subset \mathcal{O}_{V_1}[W]).$$
\end{definition}

Assume that, after restriction to an affine open set of $V$,
$\mathcal{G}$ is generated by $\mathcal{F}=\{ g_{n_i}W^{n_i},
n_i>0 , 1\leq i\leq m \}$. The total transform of $g_{n_i}$ is an
element of $I(H)^{n_i}$ (i.e., vanishes along the hypersurface $H$
with order at least $n_i$). The total transform of $\mathcal{G}$
is the algebra generated by the same $\mathcal{F}$, but now as sub-algebra in $\mathcal{O}_{V'}[W]$.

There is an open covering of $V_1$, so that locally
$g_{n_i}\mathcal{O}_{V'}=<g'_{n_i}>\cdot I(H)^{n_i}$ , for a
principal ideal spanned by some $g'_{n_i}$. This defines
$g'_{n_i}$ locally, and up to a unit. Every such $g'_{n_i}$ will
be called a {\em weighted transform} of $g_{n_i}$.

\begin{proposition}\label{prop6.6.2} Let $\mathcal{G}=\bigoplus I_k\cdot W^k (\subset \calo_V[W])$ be
a Rees algebra over $V$, generated by elements
$\{g_{N_1}W^{N_1},\dots ,g_{N_s}W^{N_s}\}$; and let $V\leftarrow
V_1$ be a permissible monoidal transformation. Then the weighted
transform of $\mathcal{G}$ is generated by
$\{g'_{N_1}W^{N_1},\dots ,g'_{N_s}W^{N_1}\}$, where each
$g'_{N_i}$ is a weighted transform of $g_{N_i}$.
\end{proposition}
(see \cite{VV3} Prop. 1.3).
\begin{remark} If $\mathcal{G}=\bigoplus I_k\cdot W^k \subset \mathcal{G}'=\bigoplus I_k\cdot W^k$ is a finite extension of graded
algebras in $\calo_V[W]$, then
$Sing(\mathcal{G})=Sing(\mathcal{G}')$ (\ref{3propsing},2)), and
there is an inclusion of the weighted transforms, say
$(\mathcal{G})_1 \subset (\mathcal{G}')_1$ which is a finite
extension of algebras over $V_1$.

\end{remark}

\begin{remark} It is convenient to pass from an arbitrary Rees algebra to a
Diff-algebra, namely to the differential extension. In fact for the sake of resolution replacing a Rees algebra by Diff-algebra it spans can be taken for free. On the other hand, our previous results show that Diff-algebras have very powerful properties. 
The differential extension does not affect the singular locus, so a permissible transformation for one is permissible for both. However
the weighted transform of a Diff-algebra is not necessarily a
Diff-algebra.
 But still elimination algebras are well defined. This fundamental property, which is based on the {\em stability of transversality}, will be discuss in \ref{par667}.
% This raises some questions of compatibility of these
%extensions and with monoidal transformations.

Suppose that a hypersurface $X$, embedded in a smooth scheme $V$,
has points of multiplicity $b$ and no point of multiplicity $b+1$.
In order to study the b-fold points we consider the Rees algebra,
say $\mathcal{G}=\calo_V[I(X)W^b]$.  Note that $Sing(\mathcal{G})$
is the closed set of points where the hypersurface has
multiplicity $b$. Let $V\leftarrow V_1 $ be a monoidal
transformation at $Y\subset Sing(\mathcal{G})$. Let
$\mathcal{G}_1$ be the weighted transform of $\mathcal{G}$, and
set $X_1$ as the strict transform of $X$. Then $\mathcal{G}_1$ is
the Rees algebra $\calo_{V_1}[I(X_1)W^b]$, so
$Sing(\mathcal{G}_1)$ is the set of $b$-fold points of $X_1$.

Consider again the situation at $V$. We pass from an ordinary Rees
algebra $\mathcal{G}$ to the Diff-algebra that it spans, say
$\mathcal{G}'$. 
%So locally at $x\in Sing(\mathcal{G})=
%Sing(\mathcal{G}')$, we have defined a transversal projection
%$\mathcal{R}_{\mathcal{G}'}$.
The monoidal transform $V\leftarrow V_1 $ defines a transform of
$\mathcal{G}'$, say $\mathcal{G}'_1$, and clearly
$\mathcal{G}_1\subset \mathcal{G}'_1.$

Our interest is on the $b$-fold points of $X_1$, namely on
$Sing(\mathcal{G}_1)$ . We can consider the Diff-algebra generated
by $\mathcal{G}_1$.  It is therefore clear, that in a step by step
argument, we would like to relate the Diff-algebra spanned by
$\mathcal{G}_1$, with that spanned by $\mathcal{G}_1'$. In Theorem
\ref{dcodgcmo} we address this form of compatibility of
transformations with extensions by differential operators.
\end{remark}

\begin{parrafo} {\em Fix
$\mathcal{G}=\bigoplus I_n W^n \subset \mathcal{G'}=\bigoplus I'_n
W^n \subset \mathcal{G}''=\bigoplus I''_n W^n, $ over $V$ as in
Corollary \ref{coragre}, so that $\mathcal{G}''$ is the
Diff-algebra spanned by $\mathcal{G}$. Then
$Sing(\mathcal{G})=Sing(\mathcal{G'})=Sing(\mathcal{G''})$ (Prop
\ref{3propsing}). In particular, a monoidal transformation
$V\longleftarrow V_1$, at a smooth center $Y \subset
Sing(\mathcal{G})$, defines weighted transforms, say
$\mathcal{G}^{(1)}$, $\mathcal{G'}^{(1)}$, and
$\mathcal{G''}^{(1)}$.

The following result is the so-called Giraud's Lemma, formulated
here in the context of Rees algebras.

}

\end{parrafo}

\begin{theorem}\label{dcodgcmo} Let $\mathcal{G}\subset  \mathcal{G'} \subset \mathcal{G''}  (\subset \mathcal{O}_{V}[W])$ be
an inclusion of Rees algebras as above, so that $\mathcal{G}''$ is
the Diff-algebra spanned by  $\mathcal{G}$. Fix a monoidal
transformation $V\longleftarrow V_1$ with center $Y\subset
Sing(\mathcal{G})$. Then

1) $\mathcal{G}^{(1)}\subset  \mathcal{G'}^{(1)} \subset
\mathcal{G''}^{(1)} (\subset \mathcal{O}_{V_1}[W])$, and

2) all three Rees algebras in 1) span the same Diff-algebras over
$V_1$.

In particular, the condition on the inclusion $\mathcal{G} \subset
\mathcal{G'}$ in Corollary \ref{coragre}, is preserved by weighted
transformations of Rees algebras.

\end{theorem}

\begin{center} {\bf Elimination and monoidal transformations.}
\end{center}

\begin{parrafo}\label{par667} {\rm Consider the localization, say $\mathcal{G}_x=\bigoplus I_k\cdot W^k
(\subset \calo_{V,x}[W])$, of a Diff-algebra $\mathcal{G}$ at a
simple closed point $x\in Sing(\mathcal{G})$. Since the
Diff-algebra is simple, there is an index $c_1$ and an element
$f_{c_1}\in I_{c_1}$ of order $c_1$ at $\calo_{V,x}$. At a
suitable \'etale neighborhood, a smooth transversal morphism 
$\beta: V \to V'$ can be defined on a smooth scheme $V'$, $\dim V'=\dim V-1$.

At a suitable \'etale neighborhood
%$\hat{\calo}_{V,x}$, 
we can also assume that $f_{c_1}$ is a monic
polynomial of degree $c_1$ in $S[Z]$, and of order $c_1$ in
$R=S[Z]_{<M_S,Z>}$, where $S=\calo_{V',x}$, $R\subset \calo_{V,x}$,
and the inclusion is \'etale.
%$\hat{R}=\hat{\calo}_{V,x}$. 
This defines by restriction a finite morphism, say $\beta$ again:
$$\beta: Spec(S[Z]/\langle
f_{c_1}(Z)\rangle )\to Spec(S), $$
together with an elimination algebra $\mathcal{R}_{\mathcal{G}}\subset S[W].$

Here we view $X=Spec(S[Z]/\langle f_{c_1}(Z)\rangle )$ as a
hypersurface in $V$. Theorem \ref{4th10}, (i), asserts that,
locally at $x$, $Sing(\mathcal{G})$ is included in the $c_1$-fold
points of this hypersurface (i.e., in $V(Diff^{c_i-1}(\langle
f_{c_1}(Z) \rangle ))$. Moreover, as $\mathcal{G}$ is a Diff-algebra,
$\beta(V(Diff^{c_i-1}))=Sing(\mathcal{R}_{\mathcal{G}})$.

Under these conditions the multiplicity formula asserts $\beta$ is
one to one on the closed set $V(Diff^{c_i-1}(\langle  f_{c_1}(Z)
\rangle ))$. Furthermore, if $Y$ is a closed and smooth in
$V(Diff^{c_i-1}(\langle f_{c_1}(Z) \rangle ))$, then $Y$ is
isomorphic to $\beta(Y)(\subset Spec(S))$ (see \cite{ZS}, Corollary
1, page 299). So both $Y$ in $V$, and $\beta(Y)$ in $Spec(S)$, are
regular centers. 
Let now $V\leftarrow V_1$, and $Spec(S)\leftarrow U$, denote the
monoidal transformations at $Y$ and $\pi(Y)$ respectively. As $\pi(Y)\subset Sing(\mathcal{R}_{\mathcal{G}})$,  $Spec(S)\leftarrow U$ defines a weak transform, say:
$$ (\mathcal{R}_{\mathcal{G}})_1\subset \calo_U[W]$$
Let
$X'$ denote the strict transform of $X$. The hypersurface $X'$ has
at most points of multiplicity $c_1$. Let $F(\subset X')$ denote
the closed set of points of multiplicity $c_1$. It is well known that locally at a neighborhood of $F$, there is a smooth morphism 
$\beta_1: V_1\to U$ that defines a commutative square. Moreover, it defines, by restriction, finite morphism say $\beta_1:X' \to U$,
compatible with $\beta$.
As the center $Y$ was chosen in $Sing(\mathcal{G})$, a
weighted transform, say
$$\mathcal{G}_1=\bigoplus I^{(1)}_n\cdot W^k (\subset \mathcal{O}_{V'}[W])$$
is defined. As $I(X)\subset I_{c_1}$, it follows that
$I(X')\subset I^{(1)}_{c_1}$ (where $I(X)$ and $I(X')$ are the
sheaves of ideals defining $X$ and $X'$ respectively). In
particular $Sing(\mathcal{G}_1)\subset F$ (points of multiplicity
$c_1$).
So locally at $F$ there is a finite morphism
$\beta_1: X'=Spec(S'[Z]/\langle
f'_{c_1}(Z)\rangle )\to U , $ where $f'_{c_1}$ is a weighted (or
strict) transform of $f_{c_1}$. 
%And $\Sing(\mathcal{G}_1)\subset F$
%(the $c_1$ fold points of $X'$.

One can check that $\beta_1$ is transversal to $X'$ and to  $\mathcal{G}_1$ at any point $y\in \Sing(\mathcal{G}_1)\subset F$ mapping to $x$. 
%that there is a smooth transversal morphism $\pi_1:V'\to U$,  in a neighborhood of $\Sing(\mathcal{G}_1)$, defining 
In can be proved that $\beta_1$ defines an elimination algebra, say $\mathcal{R}_{\mathcal{G}_1}$, and moreover:
$$\mathcal{R}_{\mathcal{G}_1}=(\mathcal{R}_{\mathcal{G}})_1.$$
To be precise, we mean here that both algebras have the same integral closure. This provides a natural commutativity of elimination with monoidal transformations, sustained on the so called {\em stability of transversality } by monoidal transformations. This property is not addressed in this work. The point is that  $\mathcal{G}_1$ is no longer a Diff-algebra, but it is a $\beta_1$-differential algebra, which is all we need to define $\mathcal{R}_{\mathcal{G}_1}$.

In particular we can only guarantee that
$\beta_1(\Sing(\mathcal{G}_1))\subset \Sing(\mathcal{R}_{\mathcal{G}_1})$. 
%

%If $\mathcal{G}'_1$ denotes the
%Diff-algebra spanned by $\mathcal{G}_1$, then Theorem
%\ref{4th10} ensure that locally at $\pi_1(y)$ there is an elimination
%algebra, say $\mathcal{R}_{\mathcal{G}'_1}\subset \calo_{U,\pi_1(y)}[W],$ and 
%$$\pi_1(\Sing(\mathcal{G}_1))=\pi_1(\Sing(\mathcal{G}'_1))\subset \Sing(\mathcal{R}_{\mathcal{G}'_1}).$$

}
\end{parrafo}

%\begin{remark} Suppose that locally at the simple point $x\in Sing(\mathcal{G})$, as
%above, we fix two projections, say $\pi_1: V\to V^{(1)}$ and
%$\pi_2: V\to V^{(2)}$ (two morphisms of smooth schemes
%($dim(V^{(1)})=dim(V^{(2)})=d-1$). The monoidal transformation
%$V\leftarrow V'$ at a center $Y\subset Sing(\mathcal{G})$ defines
%a transform, say $\mathcal{G}_1$, and this in turn spans a
%Diff-algebra $\mathcal{G}'_1$.

%The previous discussion shows that one can define, say
%$V^{(1)}\leftarrow U^{(1)}$, the monoidal transformation with
%center $\pi_1(Y)$; and $V^{(2)}\leftarrow U^{(2)}$, with center
%$\pi_2(Y)$. It also indicates that if we fix a point $y\in
%Sing(\mathcal{G}'_1)$, mapping to $x$, we can define two
%Diff-algebras, say

%$$\mathcal{R}^{(1)}_{\mathcal{G}'_1}\subset \calo_{U^{(1)},\pi_1'(y)}[W] \\ \mbox{ and }
%\mathcal{R}^{(2)}_{\mathcal{G}'_2}\subset \calo_{U^{(2)},\pi_2'(y)}[W].$$

%The same argument used in the proof of Theorem \ref{th55} (see
%\ref{55e}), shows that
%$$
%ord_{\pi'_1(y)}(\mathcal{R}^{(1)}_{\mathcal{G}'_1})=ord_{\pi'_2(y)}(\mathcal{R}^{(2)}_{\mathcal{G}'_2}).$$

%Furthermore, a similar result holds after applying several
%permissible monoidal transformations over the original simple
%Diff-algebra $\mathcal{G}$.

%\end{remark}

\begin{parrafo}\label{668} {\em

%Note that in the setting of \ref{par667}, since $Y \subset
%Sing(\mathcal{G})$, $\pi(Y)\subset
%Sing(\mathcal{R}_{\mathcal{G}})$, and there is also a weighted
%transform, say
%$$(\mathcal{R}_{\mathcal{G}})_1 \subset \mathcal{O}_{U'}[W].$$
%The question now is to relate the Rees algebra
%$(\mathcal{R}_{\mathcal{G}})_1$ with
%$\mathcal{R}_{\mathcal{G}'_1}$, locally at the point $\pi(y)$. One
%can check, from the definition of elimination algebra, that there
%is an inclusion:
%\begin{equation}\label{efin-1}
%(\mathcal{R}_{\mathcal{G}})_1\subset \mathcal{R}_{\mathcal{G}'_1}.
%\end{equation}

Here $(\mathcal{R}_{\mathcal{G}})_1$ is the transform of
$\mathcal{R}_{\mathcal{G}}$ by one monoidal transformation. If we
could ensure that $Sing(\mathcal{R}_{\mathcal{G}})_1=
\beta_1(Sing(\mathcal{G}_1))$, we could identify the singular locus
of $\mathcal{G}_1$ (i.e., of $\mathcal{G}'_1$) with the singular
locus of the transform of $\mathcal{R}_{\mathcal{G}}$. If
furthermore, this link between $\mathcal{G}$ and
$\mathcal{R}_{\mathcal{G}}$ is preserved by every sequence of
monoidal transformations, then we have achieved a way of
representing the singular locus of $\mathcal{G}$ which is stable
by monoidal transformations (see property (C), of stability of
elimination, in the Introduction). So this would lead to resolution of singularities.

This will be the case when the invariant $e_0$ attached to the
simple point $x$ in \ref{partau} is zero. In this case $f_{c_0}$
can be taken to be a monic polynomial of degree 1, defining
therefore a smooth hypersurface of maximal contact, and then
$\mathcal{R}_{\mathcal{G}}$ is naturally identified with the {\em restriction} of
$\mathcal{G}$ at such smooth hypersurface. In such case
% When this holds, then 
%the terms in
%(\ref{efin-1}) span the same Diff-algebra. In particular:
\begin{equation}\label{efin1}
Sing((\mathcal{R}_{\mathcal{G}})_1) (=
Sing(\mathcal{R}_{\mathcal{G}_1}))=\beta_1(Sing(\mathcal{G}_1)).
\end{equation}
And furthermore, a similar result holds for every sequence of, say
$r$ transformations, namely  
$(\mathcal{R}_{\mathcal{G}})_r=
\mathcal{R}_{\mathcal{G}_r}$ 
%and both terms span the same
%Diff-algebra.
 In particular:
\begin{equation}\label{efin2}
Sing((\mathcal{R}_{\mathcal{G}})_r)=
Sing(\mathcal{R}_{\mathcal{G}_r})=\beta_r(Sing(\mathcal{G}_r)).
\end{equation}
This is not the case, at least in general, when the invariant
$e_0$ is non-zero. 
%%In fact, in such case it can even happen that 
%%$(\mathcal{R}_{\mathcal{G}})_1$ and $\mathcal{R}_{\mathcal{G}'_1}$
%do not span the same Diff-algebra. The following examples illustrate this pattern
%of behavior.

 }

\end{parrafo}
\begin{example} The following example illustrates the general fact that
$(\mathcal{R}_{\mathcal{G}})_1$ and $\mathcal{R}_{\mathcal{G}_1}$
coincide. 
%Moreover, in characteristic
%zero, $e_0$ is always zero, so (\ref{efin1}) also
%holds.

Set $R=\mathbb{Q}[Y,Z]$. Consider the Rees algebra in $R[W]$
generated (over $R$) by the element $(Z^2+Y^5)W^2$. Let
$\mathcal{G}=\bigoplus I_k\cdot W^k (\subset R[W])$ be the
Diff-algebra spanned by this Rees algebra. According to Theorem
\ref{3th1} and formula (\ref{eq3422}), $\mathcal{G}$ is generated
by $\{2ZW, 5Y^4W,(Z^2+Y^5)W^2\}$, or say $\{ZW, Y^4W\}$. As $Z\in
I_1$, we can choose $c_1=1$ and $f_{c_1}(Z)=Z$ in
Theorem \ref{4th10}. One can check now that the projection in
$S=k[Y]$, namely $\mathcal{R}_{\mathcal{G}}$ is generated by
$\{Y^4W\}$.

Consider the quadratic transformation at the relevant chart
$\mathbb{Q}[Y,Z_1]$, where $Y\cdot Z_1=Z$. Here $\mathcal{G}_1$ is
generated by $\{Z_1 W,Y^3W\}$ , which is also a Diff-algebra. The elimination algebra, say $\mathcal{R}_{\mathcal{G}_1}$ is the $S$
subalgebra in $S[W]$ generated by $\{Y^3W\}$. Finally check
that $\mathcal{R}_{\mathcal{G}_1}=(\mathcal{R}_{\mathcal{G}})_1$.

\end{example}

There is yet another natural observation. If (\ref{efin2}) modified by taking $\mathcal{G}'_r$ as the Diff-algebra spanned by $\mathcal{G}_r$, then Theorem \ref{4th10}, (i),  says that 
$Sing(\mathcal{R}_{\mathcal{G}'_r})=\beta_r(Sing(\mathcal{G}'_r))$.
Resolution of singularities would be achieved if we could show that 
$$Sing(\mathcal{R}_{\mathcal{G}_r})=Sing(\mathcal{R}_{\mathcal{G}'_r}).$$
Unfortunately this does not hold in positive characteristic.

\begin{example}\label{exf2} The following is an example of a pathology that can only occur in positive characteristic. Namely that 
$(\mathcal{R}_{\mathcal{G}})_1$ and $\mathcal{R}_{\mathcal{G}'_1}$
do not span  Diff-algebras with the same integral closure, where 
$\mathcal{G}'_1$ denotes the Diff-algebra spanned by $\mathcal{G}_1$.

Fix a field $k$ of characteristic two, and set $R=k[Y,Z]$. Let
$\mathcal{G}=\bigoplus I_k\cdot W^k (\subset R[W])$ be, as before,
the Diff-algebra generated by $(Z^2+Y^5)W^2$. Here $\mathcal{G}$
is generated by $\{Y^4W,(Z^2+Y^5)W^2\}$. The elimination in
$S=k[Y]$, namely $\mathcal{R}_{\mathcal{G}}$, is generated by
$\{Y^8W^2\}$, so up to integral closure it is generated
$\{Y^4W\}$. To check this fact consider $$\beta: Spec(S[Z]/\langle
f_{c_1}(Z)\rangle )\to Spec(S) ,$$ where $f_{c_1}=Z^2+Y^5$. This
is a purely inseparable extension, and 
$\mathcal{R}_{\mathcal{G}}$ is (up to integral closure) generated
by the coefficients of the characteristic polynomial of
multiplication by $Y^4$ in the free $S$-module $B=S[Z]/\langle
f_{c_1}(Z)\rangle$ (see Corollary \ref{4cor10}). In this case
$\mathcal{H}_{\mathcal{G}}$ is generated by $Y^8W^2$, so up to
integral closure $\mathcal{R}_{\mathcal{G}}$ is generated by
$Y^4W$.

Consider the quadratic transformation at the relevant chart
$k[Y,Z_1]$, where $Y\cdot Z_1=Z$. Here $\mathcal{G}_1$ is
generated by $\{Y^3W,(Z_1^2+Y^3)W^2\}$.

Let $\mathcal{G}_1'$ be the Diff-algebra spanned by
$\mathcal{G}_1$. Then
$\mathcal{R}_{\mathcal{G}'_1}$ is, up to integral closure,
generated by $\{Y^2W\}$, which is already a Diff-algebra. So 
$\mathcal{R}_{\mathcal{G}'_1}$ is a Diff-algebra, and it is different from the Diff-algebra spanned by $(\mathcal{R}_{\mathcal{G}})_1$. In fact
$(\mathcal{R}_{\mathcal{G}})_1$ is generated by $\{Y^3W\}$.
\end{example}

\begin{parrafo} {\em The situation in \ref{exf2}, in which in which
$(\mathcal{R}_{\mathcal{G}})_1$ and $\mathcal{R}_{\mathcal{G}'_1}$
do not span the same Diff-algebra, can only occur when the
invariant $e_0$ is not zero. In this case $p=2$, and $e_0=1$ at
the origin, which is a singular point of Diff-algebra $\mathcal{G}=\bigoplus I_k\cdot
W^k (\subset R[W])$.

An interesting case, also computable from our invariants, is that
of $Z^3+X^{13}Z+X^{16}$ in characteristic 3. The invariant $e_0$
at the origin is 1. This curve is analytically irreducible, the
singularity is resolved with five quadratic transformations.

In this case $\mathcal{R}_{\mathcal{G}'_i}$ spans the same
Diff-algebra as $(\mathcal{R}_{\mathcal{G}})_i$ for the first two
transformations (for $i=0,1$), but not for the next 3 quadratic
transformations.}
\end{parrafo}

\end{document}